\documentclass{article}
\usepackage{amsmath}
\usepackage{amssymb}
\usepackage{amsfonts}
\usepackage{graphicx}
\providecommand{\U}[1]{\protect\rule{.1in}{.1in}}
\newtheorem{theorem}{Theorem}

\newtheorem{definition}[theorem]{Definition}

\begin{document}

\title{Introduction to the method of multiple scales}
\author{Per Kristen Jakobsen}
\maketitle
\tableofcontents

\section{Introduction}

Perturbation methods are aimed at finding approximate analytic solutions to
problems whose exact analytic solutions can not be found. The setting where
perturbation methods are applicable is where there is a family of equations,
 $\mathcal{P}(\varepsilon)$, depending on a parameter $\varepsilon<<1$, and
where $\mathcal{P}(0)$ has a known solution. Perturbation methods are designed
to construct solutions to $\mathcal{P}(\varepsilon)$ by adding small
corrections to known solutions of $\mathcal{P}(0)$. The singular aim of
perturbation methods is to calculate corrections to solutions of
$\mathcal{P}(0)$. Perturbation methods do not seek to prove that a solution of
$\mathcal{P}(0)$, with corrections added, is close to a solution of
$\mathcal{P}(\varepsilon)$ for $\varepsilon$ in some finite range with respect
to some measure of error. It's sole aim is to compute corrections and to make
sure that the first correction is small with respect to the choosen solution
of $\mathcal{P}(0)$, that the second correction is small with respect to the
first correction and so on, all in the limit when $\varepsilon$ approaches
zero. This formal nature and limited aim of is why we prefer to call it
\textit{perturbation methods} rather than \textit{perturbation theory}. A
mathematical theory is a description of proven mathematical relations among
well defined objects of human thought. Perturbation methods does not amount to
a mathematical theory in this sense. It is more like a very large bag of
tricks, whose elements have a somewhat vague domain of applicability, and where
the logical relations between the tricks are not altogether clear, to put it nicely.

After all this negative press you might ask why we should bother with this
subject at all, and why we should not rather stay with real, honest to God,
mathematics. The reason is simply this: If you want analytic solutions to
complex problems, it is the only game in town. In fact, for quantum theory,
which is arguably our best description of reality so far, perturbation methods
is almost always the first tool we reach for. For the quantum theory of
fields, like quantum electrodynamics, perturbation methods are essentially the
only tools available. These theories are typically only known in terms of
perturbation expansions. You could say that we don't actually know what the
mathematical description of these very fundamental theories is. But at the
same time, quantum theory of fields give some of the most accurate,
experimentally verified, predictions in all of science.

So clearly, even if perturbation methods are somewhat lacking in mathematical
justification, they work pretty well. And in the end that is the only thing
that really counts.

These lecture notes are not meant to be a general introduction to the wide
spectrum of perturbation methods that are used  all across science. Many textbooks
exists whose aim is to give such a broad overview, an overview that includes the most
commonly used perturbation methods\cite{Hinch},\cite{Nayfeh},\cite{Holmes}%
,\cite{Murdoc}. Our aim is more limited; we focus on one such method, which 
is widely used in many areas of applied science. This is the \textit{method of
multiple scales}. The method of multiple scales is described in all
respectable books on perturbation methods and there are also more specialized
books on singular perturbation methods where the method of multiple scales has
a prominent place\cite{Kevorkian},\cite{Johnson}. There are, however, quite
different views on how the method is to be applied, and what its limitations
are. Therefore, descriptions of the method appears quite different in the
various sources, depending on the views of the authors. In these lecture notes
we decribe the method in a way that is different from most textbooks, but
which is very effective and makes it possible to take the perturbation
expansions to higher order in the small perturbation parameter that would otherwise be possible. The source that
is closest to our approach is \cite{NLO}.

We do not assume that the reader has had any previous exposure to perturbation
methods. These lecture notes therefore starts off by introducing the
basic ideas of asymptotic expansions and illustrate them  using
algebraic equations. The lecture notes then proceeds by introducing regular perturbation expansions for 
single ODEs, study the breakdown of these expansions, and show how to avoid the breakdown using the
method of multiple scales.  The method of multiple scales is then generalized to systems of ODEs, boundary layer problems for ODEs  and to PDEs.
In the last section we illustrate the method of multiple scales by applying it to the Maxwells equations; showing how the Nonlinear Schrodinger equation appears
as an approximation to the Maxwell equations in a situation where dispersion and nonlinearity balances. 
 Several exercises involving multiple
scales for ODEs and PDEs are included in the lecture notes.

\section{Regular and singular problems, applications to algebraic equations.}

In this section we will introduce perturbation methods in the context of
algebraic equations. One of the main goals of this section is to introduce the
all-important distinction between regular and singular perturbation problems,
but we also use the opportunity to introduce the notion of a
\textit{perturbation hierarchy} and describe some of its general properties.

\subsection*{Example 1: A regularly perturbed quaderatic equation}

Consider the polynomial equation
\begin{equation}
x^{2}-x+\varepsilon=0.\label{Eq2.1}%
\end{equation}
This is our \textit{perturbed problem }$\mathcal{P}(\varepsilon)$. The
\textit{unperturbed problem }$\mathcal{P}(0)$, is
\begin{equation}
x^{2}-x=0.\label{eq2.2}%
\end{equation}
This unperturbed problem is very easy to solve%
\begin{align}
x^{2}-x  &  =0,\nonumber\\
&  \Updownarrow\nonumber\\
x_{0}  &  =0,\nonumber\\
x_{1}  &  =1.\label{Eq2.3}
\end{align}
Let us focus on $x_{1}$ and let us assume that the \textit{perturbed problem}
has a solution in the form of a \textit{perturbation expansion}%
\begin{equation}
x(\varepsilon)=a_{0}+\varepsilon a_{1}+\varepsilon^{2}a_{2}+...\;\;.\label{Eq2.4}%
\end{equation}
where $a_{0}=1$. Our goal is to find the unknown numbers $a_{1},a_{2},..$\;.
These numbers should have a size of order $1$. This will ensure that
$\varepsilon a_{1}$ is a small correction to $a_{0}$, that $\varepsilon
^{2}a_{2}$ is a small correction to $\varepsilon a_{1}$and so on, all in the
limit of small $\varepsilon$. As we have stressed before, maintaining the
ordering of the perturbation expansion is the one and only unbreakable rule
when we do perturbation calculations. \ The perturbation method now proceeds
by inserting the expansion (\ref{Eq2.4}) into equation (\ref{Eq2.1}) and
collecting terms containing the same order of $\varepsilon$.%
\begin{gather}
(a_{0}+\varepsilon a_{1}+\varepsilon^{2}a_{2}+...)^{2}-(a_{0}+\varepsilon
a_{1}+\varepsilon^{2}a_{2}+...)+\varepsilon=0,\nonumber\\
\Downarrow\nonumber\\
a_{0}^{2}+2\varepsilon a_{0}a_{1}+\varepsilon^{2}(a_{1}^{2}+2a_{0}a_{2}%
)-a_{0}-\varepsilon a_{1}-\varepsilon^{2}a_{2}-..+\varepsilon=0,\nonumber\\
\Downarrow\nonumber\\
a_{0}^{2}-a_0+\varepsilon(2a_{0}a_{1}-a_{1}+1)+\varepsilon^{2}(2a_{0}a_{2}%
+a_{1}^{2}-a_{2})+...=0.\label{Eq2.5}%
\end{gather}
Since $a_{1},a_{2},..$ are all assumed to be of order $1$ this equation will
hold in the limit when $\varepsilon$ approach zero only if
\begin{align}
a_{0}^{2}-a_0  &  =0,\nonumber\\
2a_{0}a_{1}-a_{1}+1  &  =0\nonumber\\
2a_{0}a_{2}+a_{1}^{2}-a_{2}  &  =0.\label{eq2.6}
\end{align}
We started with one nonlinear equation for $x$, and have ended up with three
coupled nonlinear equations for $a_{0}$, $a_{1}$ and $a_{2}$. Why should we
consider this to be progress? It seems like we have rather substituted one
complicated problem with one that is even more complicated!

The reason why this is progress, is that the coupled system of nonlinear
equations has a very special structure. We can rewrite it in the form%
\begin{align}
a_{0}(a_{0}-1)  &  =0,\nonumber\\
(2a_{0}-1)a_{1}  &  =-1,\nonumber\\
(2a_{0}-1)a_{2}  &  =-a_{1}^{2}.\label{Eq2.7}
\end{align}
The first equation is nonlinear but simpler than the perturbed equation
(\ref{Eq2.1}), the second equation is linear in the variable $a_{1}$ and that
the third equation is linear in the variable $a_{2}$ when $a_{1}$ has been
found. Moreover, the linear equations are all determined by the same linear
operator\ $\mathcal{L}(\cdot)=(2a_{0}-1)(\cdot)$. This reduction to a simpler
nonlinear equation and a sequence of linear problems determined by the same
linear operator is what makes (\ref{Eq2.7}) essentially simpler than the
original equation (\ref{Eq2.1}), which does not have this special structure.
The system (\ref{Eq2.7}) is called a \textit{perturbation hierarchy} for
(\ref{Eq2.1}). The special structure of the perturbation hierarchy is key to
any successful application of perturbation methods, whether it is for
algebraic equations, ordinary differential equations or partial differential equations.

The perturbation hierarchy (\ref{Eq2.7}) is easy to solve and we find
\begin{align}
a_{0}  &  =1,\nonumber\\
a_{1}  &  =-1,\nonumber\\
a_{2}  &  =-1,\label{Eq2.8}
\end{align}
and our perturbation expansion to second order in $\varepsilon$ is
\begin{equation}
x(\varepsilon)=1-\varepsilon-\varepsilon^{2}+...\label{eq2.9}%
\end{equation}
For this simple case we can solve the unperturbed problem directly using the
solution formula for a quaderatic equation. Here are some numbers

\bigskip%

\begin{tabular}
[c]{lll}%
$\varepsilon$ & Exact solution & Perturbation solution\\
$0.001$ & $0.998999$ & $0.998999$\\
$0.01$ & $0.989898$ & $09989900$\\
$0.1$ & $0.887298$ & $0.890000$%
\end{tabular}

\bigskip

\noindent We see that our perturbation expansion is quite accurate even for
$\varepsilon$ as large as $0.1$.

Let us see if we can do better by finding an even more accurate approximation through
extention of the perturbation expansion to higher order in $\varepsilon$. In fact
let us take the perturbation expansion to infinite order in $\varepsilon$.
\begin{equation}
x(\varepsilon)=a_{0}+\epsilon a_{1}+\epsilon^{2}a_{2}+...=a_{0}+
{\displaystyle\sum_{n=1}^{\infty}}
\varepsilon^{n}a_{n}\label{Eq2.10}
\end{equation}
Inserting (\ref{Eq2.10}) into (\ref{Eq2.1}) and expanding we get
\begin{gather}
(a_{0}+
{\displaystyle\sum_{n=1}^{\infty}}
\varepsilon^{n}a_{n})(a_{0}+
{\displaystyle\sum_{m=1}^{\infty}}
\varepsilon^{m}a_{m})-a_{0}-
{\displaystyle\sum_{n=1}^{\infty}}
\varepsilon^{n}a_{n}+\varepsilon=0,\nonumber\\
\Downarrow\nonumber\\
a_{0}^{2}-a_{0}+
{\displaystyle\sum_{p=1}^{\infty}}
\varepsilon^{p}(2a_{0}-1)a_{p}+
{\displaystyle\sum_{p=2}^{\infty}}
\varepsilon^{p}\left(
{\displaystyle\sum_{m=1}^{p-1}}
a_{m}a_{p-m}\right)  +\varepsilon=0,\nonumber\\
\Downarrow\nonumber\\
a_{0}^{2}-a_{0}+\varepsilon\left(  (2a_{0}-1)a_{1}+1\right)  +
{\displaystyle\sum_{p=2}^{\infty}}
\varepsilon^{p}\left(  (2a_{0}-1)a_{p}+
{\displaystyle\sum_{m=1}^{p-1}}
a_{m}a_{p-m}\right)  =0.\label{Eq2.11}
\end{gather}
Therefore the complete perturbation hierarchy is
\begin{align}
a_{0}(a_{0}-1)  &  =0,\nonumber\\
(2a_{0}-1)a_{1}  &  =-1,\nonumber\\
(2a_{0}-1)a_{p}  &  =-{\displaystyle\sum_{m=1}^{p-1}}a_{m}a_{p-m},\text{ \ \ \ }p\geqq2.
\end{align}
The right-hand side of the equation for $a_{p}$ only depends on $a_{j}$ for
$j<p$. Thus the perturbation hierarchy is an infinite system of linear
equations that is coupled in such a special way that we can solve them one by
one. The perturbation hierarchy truncated at order $4$ is
\begin{align}
(2a_{0}-1)a_{1}  &  =-1,\nonumber\\
(2a_{0}-1)a_{2}  &  =-a_{1}^{2},\nonumber\\
(2a_{0}-1)a_{3}  &  =-2a_{1}a_{2},\nonumber\\
(2a_{0}-1)a_{4}  &  =-2a_{1}a_{3}-a_{2}^{2}.\label{Eq2.13}
\end{align}
Using $a_{0}=1$, the solution to the hierarchy is trivially found to be
\begin{align}
a_{1}  &  =-1,\nonumber\\
a_{2}  &  =-1,\nonumber\\
a_{3}  &  =-2,\nonumber\\
a_{4}  &  =-5.\label{Eq2.14}
\end{align}
For $\varepsilon=0.1$ the perturbation expansion gives
\begin{equation}
x(0.1)=0.8875...\;\;,\label{Eq2.15}%
\end{equation}
whereas the exact solution is
\begin{equation}
x(0.1)=0.8872...\;\;.\label{Eq2.16}%
\end{equation}
we are clearly getting closer. However we did not get all that much in return
for our added effort.

Of course we did not actually have to use perturbation methods to find
solutions to equation (\ref{Eq2.1}), since it is exactly solvable using the
formula for the quaderatic equation. The example, however, illustrate many
general features of perturbation calculations that will appear again and again
in different guises.

\subsection*{Example 2: A regularly perturbed quintic equation}

Let us consider the equation
\begin{equation}
x^{5}-2x+\varepsilon=0.\label{Eq2.17}
\end{equation}
This is our perturbed problem, $\mathcal{P}(\varepsilon)$. For this case
perturbation methods are neccessary, since there is no solution formula for
general polynomial equations of order higher than four. The unperturbed
problem, $\mathcal{P}(0)$, is
\begin{equation}
x^{5}-2x=0.\label{Eq2.18}
\end{equation}
It is easy to see that the unperturbed equation has a real solution
\begin{equation}
x=\sqrt[4]{2}\equiv a_{0}.\label{Eq2.19}
\end{equation}
We will now construct a perturbation expansion for a solution to
(\ref{Eq2.17}), starting with the solution $x=a_{0}$ of the unperturbed
equation (\ref{Eq2.18}). We therefore introduce the expansion
\begin{equation}
x(\varepsilon)=a_{0}+\varepsilon a_{1}+\varepsilon^{2}a_{2}+....\;\;.\label{Eq2.20}
\end{equation}
Inserting (\ref{Eq2.20}) into equation (\ref{Eq2.17}) and expanding we get
\begin{gather}
(a_{0}+\varepsilon a_{1}+\varepsilon^{2}a_{2}+..)^{5}\nonumber\\
-2(a_{0}+\varepsilon a_{1}+\varepsilon^{2}a_{2}+..)+\varepsilon=0,\nonumber\\
\Downarrow\nonumber\\
a_{0}^{5}+5a_{0}^{4}(\varepsilon a_{1}+\varepsilon^{2}a_{2}+...)+10a_{0}
^{3}(\varepsilon a_{1}+...)^{2}+..\nonumber\\
-2a_{0}-2\varepsilon a_{1}-2\varepsilon^{2}a_{2}-...+\varepsilon=0,\nonumber\\
\Downarrow\nonumber\\
a_{0}^{5}-2a_{0}+\varepsilon(1+5a_{0}^{4}a_{1}-2a_{1})+\varepsilon^{2}
(5a_{0}^{4}a_{2}+10a_{0}^{3}a_{1}^{2}-2a_{2})+...=0.\label{Eq2.21}
\end{gather}
Thus the perturbation hierarchy to order two in $\varepsilon$ is
\begin{align}
a_{0}^{5}-2a_{0} &  =0,\nonumber\\
(5a_{0}^{4}-2)a_{1} &  =-1, \nonumber\\
(5a_{0}^{4}-2)a_{2} &  =-10a_{0}^{3}a_{1}^{2}.\label{Eq2.24}
\end{align}
Observe that the first equation in the hierarchy for $a_{0}$ is nonlinear, 
whereas the equations for $a_{p}$ are linear in $a_{p}$ for $p>0$. All the
linear equations are defined in terms of the same linear operator
$\mathcal{L}(\cdot)=(5a_{0}^{4}-2)(\cdot)$. This is the same structure that we
saw in the previous example. If the unperturbed problem is linear, the first
equation in the hierarchy will also in general be linear.

  The perturbation hierarchy is easy to solve, and we find
\begin{align}
a_{1} &  =-\frac{1}{5a_{0}^{4}-2}=\frac{-1}{8},\nonumber\\
a_{2} &  =-\frac{10a_{0}^{3}a_{1}^{2}}{5a_{0}^{4}-2}=-\frac{5\sqrt[4]{8}}{256}.\label{Eq2.25}
\end{align}
The perturbation expansion to second order is then
\begin{equation}
x(\varepsilon)=\sqrt[4]{2}-\frac{1}{8}\varepsilon-\frac{5\sqrt[4]{8}}
{256}\varepsilon^{2}+...\;\;.\label{Eq2.26}
\end{equation}
Here are some numbers

\bigskip%

\begin{tabular}
[c]{lll}%
$\varepsilon$ & Exact solution & Perturbation solution\\
$0.001$ & $1.18908$ & $1.18908$\\
$0.01$ & $1.19795$ & $1.19795$\\
$0.1$ & $1.17636$ & $1.17638$%
\end{tabular}

\noindent Perturbation expansions for the other solutions to equation (\ref{Eq2.17}) can
be found by starting with the other four solutions of the equation
(\ref{Eq2.18}). In this way we get perturbation expansions for all the
solutions of (\ref{Eq2.17}), and the effort is not much larger than for the
quaderatic equation.

If we can find perturbation expansions for all the solutions of a problem
$\mathcal{P}(\varepsilon)$, by starting with solutions of the unperturbed
problem $\mathcal{P}(0)$, we say that $\mathcal{P}(\varepsilon)$ is a
\textit{regular} perturbation of $\mathcal{P}(0)$. \ If the perturbation is
not regular it is said to be \textit{singular}. This distinction applies to
all kinds of perturbation problems whether we are looking at algebraic
equations, ordinary differential equations or partial differential equations.
Clearly, for polynomial equations a neccessary condition for being a regular
perturbation problem is that $\mathcal{P}(\varepsilon)$ and $\mathcal{P}(0)$
have the same algebraic order. This is not always the case as the next example shows.

\subsection*{Example 3: A singularly perturbed quaderatic equation.}

Let us consider the following equation
\begin{equation}
\varepsilon x^{2}+x-1=0.\label{Eq2.28}
\end{equation}
This is our perturbed problem $\mathcal{P}(\varepsilon)$. The unperturbed
problem $\mathcal{P}(0)$, is
\begin{equation}
x-1=0.\label{Eq2.29}
\end{equation}
There is only one solution to the unperturbed problem
\begin{equation}
x=1\equiv a_{0}.\label{Eq2.30}
\end{equation}
Let us find a perturbation expansion for a solution to (\ref{Eq2.28}) starting
with the solution (\ref{Eq2.30}) of the unperturbed problem
\begin{equation}
x(\varepsilon)=a_{0}+\varepsilon a_{1}+\varepsilon^{2}a_{2}+...\;\;.\label{Eq2.31}
\end{equation}
Inserting (\ref{Eq2.31}) into equation (\ref{Eq2.28}) and expanding we get
\begin{gather}
\varepsilon(a_{0}+\varepsilon a_{1}+\varepsilon^{2}a_{2}+...)^{2}
+a_{0}+\varepsilon a_{1}+\varepsilon^{2}a_{2}+...-1=0,\nonumber\\
\Downarrow\nonumber\\
\varepsilon(a_{0}^{2}+2\varepsilon a_{0}a_{1}+...)+a_{0}+\varepsilon
a_{1}+\varepsilon^{2}a_{2}+...-1=0,\nonumber\\
\Downarrow\nonumber\\
a_{0}-1+\varepsilon(a_{1}+a_{0}^{2})+\varepsilon^{2}(a_{2}+2a_{0}
a_{1})+...=0.\label{Eq2.32}
\end{gather}
The perturbation hierarchy, up to second order in $\varepsilon$ is thus
\begin{align}
a_{0}  &  =1,\nonumber\\
a_{1}  &  =-a_{0}^{2},\nonumber\\
a_{2}  &  =-2a_{0}a_{1}.\label{Eq2.33}
\end{align}
The solution of the perturbation hierarchy is
\begin{align}
a_{0}  &  =1,\nonumber\\
a_{1}  &  =-1,\nonumber\\
a_{2}  &  =2,\label{Eq2.34}
\end{align}
and the perturbation expansion for the solution to (\ref{Eq2.28}) starting
from the solution $x=1$ to the unperturbed problem (\ref{Eq2.29}) is
\begin{equation}
x(\varepsilon)=1-\varepsilon+2\varepsilon^{2}+...\;\;.\label{Eq2.35}
\end{equation}
In order to find a perturbation expansion for the other solution to the
quaderatic equation (\ref{Eq2.28}), the unperturbed problem (\ref{Eq2.29}) is
of no help.

However, looking at equation (\ref{Eq2.28}) we learn something important: In
order for a solution different from $x=1$ to appear in the limit when
$\varepsilon$ approaches zero, the first term in (\ref{Eq2.28}) can not
approach zero. This is only possible if $x$ approaches infinity as
$\varepsilon$ approaches zero.

Inspired by this, let us introduce a change of variables
\begin{equation}
x=\varepsilon^{-p}y,\label{Eq2.36}
\end{equation}
where $p>0$. If $y$ is of order one, as $\varepsilon$ approaches zero, then $x $
will approach infinity in this limit, and will thus be the solution we lost in (\ref{Eq2.29}
). Inserting (\ref{Eq2.36}) into (\ref{Eq2.28}) gives us
\begin{align}
\varepsilon(\varepsilon^{-p}y)^{2}+\varepsilon^{-p}y-1 &  =0,\nonumber\\
&  \Downarrow\nonumber\\
\varepsilon^{1-2p}y^{2}+\varepsilon^{-p}y-1 &  =0,\nonumber\\
&  \Downarrow\nonumber\\
y^{2}+\varepsilon^{p-1}y-\varepsilon^{2p-1} &  =0.\label{Eq2.37}
\end{align}
The idea is now to pick a value for $p$, thereby defining a perturbed problem
$\mathcal{P}(\varepsilon)$, such that $\mathcal{P}(0)$ has a solution of order
one. For $p>1$ we get in the limit when $\varepsilon$ approches zero the
problem%
\begin{equation}
y^{2}=0,\label{Eq2.38}%
\end{equation}
which does not have any solution of order one. One might be inspired to choose
$p=\frac{1}{2}$. We then get the equation
\begin{equation}
\sqrt{\varepsilon}y^{2}+y-\sqrt{\varepsilon}=0,\label{Eq2.39}
\end{equation}
which in the limit when $\varepsilon$ approaches zero turns into
\begin{equation}
y=0.\label{Eq2.40}
\end{equation}
This equation clearly has no solution of order one. Another possibility is to
choose $p=1$. Then we get the equation
\begin{equation}
y^{2}+y-\varepsilon=0.\label{Eq2.41}
\end{equation}
In the limit when $\varepsilon$ approaches zero this equation turns into
\begin{equation}
y^{2}+y=0.\label{Eq2.42}
\end{equation}
This equation has a solution $y=-1$ which \textit{is} of order one. We
therefore proceed with this choise for $p$, and introduce a perturbation
expansion for the solution to (\ref{Eq2.41}) that starts at the solution
$y\equiv a_{0}=-1$ to the unperturbed equation (\ref{Eq2.42}).
\begin{equation}
y(\varepsilon)=a_{0}+\varepsilon a_{1}+\varepsilon^{2}a_{2}+...\;\;.\label{Eq2.43}
\end{equation}
Inserting the perturbation expansion (\ref{Eq2.43}) into equation
(\ref{Eq2.41}) and expanding we get
\begin{gather}
(a_{0}+\varepsilon a_{1}+\varepsilon^{2}a_{2}+...)^{2}+a_{0}+\varepsilon
a_{1}+\varepsilon^{2}a_{2}+...-\varepsilon=0,\nonumber\\
\Downarrow\nonumber\\
a_{0}^{2}+a_{0}+\varepsilon((2a_{0}+1)a_{1}-1)+\varepsilon^{2}((2a
_{0}+1)a_{2}+a_{1}^{2})+...=0.\label{Eq2.44}
\end{gather}
The perturbation hierarchy to second order in $\varepsilon$ is then
\begin{align}
a_{0}^{2}+a_{0} &  =0,\nonumber\\
(2a_{0}+1)a_{1} &  =1,\nonumber\\
(2a_{0}+1)a_{2} &  =-a_{1}^{2}.\label{Eq2.45}
\end{align}
We observe in passing, that the perturbation hierarchy has the special
structure we have seen earlier. The solution to the perturbation hierarchy is
\begin{align}
a_{1} &  =-1,\nonumber\\
a_{2} &  =1,\label{Eq2.46}
\end{align}
and the perturbation expansion to second order in $\varepsilon$ is
\begin{equation}
y(\varepsilon)=-1-\varepsilon+\varepsilon^{2}+...\;\;.\label{Eq2.47}
\end{equation}
Going back to the original coordinate $x$ we finally get
\begin{equation}
x(\varepsilon)=-\varepsilon^{-1}-1+\varepsilon+...\;\;.\label{Eq2.48}
\end{equation}
Even for $\varepsilon$ as large as $0.1$ the perturbation expansion and the
exact solution, $x_{E}(\varepsilon)$, are close
\begin{align}
x(\varepsilon) &  =-\varepsilon^{-1}-1+\varepsilon+...\approx
-10.900..\;\;,\nonumber\\
x_{E}(\varepsilon) &  =\frac{-1-\sqrt{1+4\varepsilon}}{2\varepsilon}
\approx-10.916...\;\;.\label{Eq2.49}
\end{align}

The perturbation problem we have discussed in this example is evidently a
singular problem. For singular problems, a coordinate transformation, like the
one defined by (\ref{Eq2.36}), must at some point be used to transform the
singular perturbation problem into a regular one.

At this point I need to be honest with you: There is really no general rule
for how to find the right transformations. Skill, experience, insight and
sometimes even dumb luck is needed to succeed. This is one of the reasons why
I prefer to call our subject perturbation methods and not perturbation theory.
Certain classes of commonly occuring singular perturbation problems have
however been studiet extensively and rules for finding the correct
transformations have been designed. In general, what one observe, is that some
kind of \textit{scaling transformation}, like in (\ref{Eq2.36}), is almost
always part of the mix.

\section{Asymptotic sequences and series.}

When using perturbation methods, our main task is to investigate the behaviour
of unknown functions $f(\varepsilon)$, in the limit when $\varepsilon$
approaches zero. This is what we did in examples one, two and three.

The way we approach this problem is to compare the unknown function $f\left(
\varepsilon\right)  $ to one or several known functions when $\varepsilon$
approaches zero. In example one and two we compared our unknow function to the known functions
$\{1,\varepsilon,\varepsilon^{2},...\}$ whereas in example three we used the
functions $\{\varepsilon^{-1},1,\varepsilon,...\}$. In order to facilitate
such comparisons, we introduce the "large-O" and "little-o" notation.

\begin{definition}
Let $f(\varepsilon)$ be a function of $\varepsilon$. Then

\begin{description}
\item[i)] $f(\varepsilon)=O(g(\varepsilon))$ $\ ,\,\ \ \,\varepsilon
\rightarrow0$ $\ \Leftrightarrow\lim_{\varepsilon\rightarrow0}\left\vert
\frac{f(\varepsilon)}{g(\varepsilon)}\right\vert \neq0,$

\item[ii)] $f(\varepsilon)=o(g(\varepsilon))\ ,\,\ \ \,\varepsilon
\rightarrow0\ \ \ \ \Leftrightarrow\lim_{\varepsilon\rightarrow0}\left\vert
\frac{f(\varepsilon)}{g(\varepsilon)}\right\vert =0.$
\end{description}
\end{definition}

\noindent Thus $\ f(\varepsilon)=O(g(\varepsilon))\,$\ means that $f(\varepsilon)$ and
$g(\varepsilon)$ are of roughly the same size when $\varepsilon$ approaches
zero and $f(\varepsilon)=o(g(\varepsilon))\,$\ means that $f(\varepsilon)$ is
much smaller than $g(\varepsilon)$ when $\varepsilon$ approaches zero.

We have for example that

\begin{enumerate}
\item $\sin(\varepsilon)=O(\varepsilon)$ $\ ,\,\ \ \,\varepsilon\rightarrow0,
$ \ \ \ because
\[
\lim_{\varepsilon\rightarrow0}\left\vert \frac{\sin(\varepsilon)}{\varepsilon
}\right\vert =1\neq0,
\]

\item $\sin(\varepsilon^{2})=o(\varepsilon)\ ,\,\ \ \,\varepsilon
\rightarrow0,$ \ \ because
\[
\lim_{\varepsilon\rightarrow0}\left\vert \frac{\sin(\varepsilon^{2}
)}{\varepsilon}\right\vert =\lim_{\varepsilon\rightarrow0}\left\vert
\frac{2\varepsilon\cos(\varepsilon^{2})}{1}\right\vert =0,
\]

\item $1-\cos(\varepsilon)=o(\varepsilon),\,\ \ \,\varepsilon\rightarrow0,$
\ \ because
\[
\lim_{\varepsilon\rightarrow0}\left\vert \frac{1-\cos(\varepsilon
)}{\varepsilon}\right\vert =\lim_{\varepsilon\rightarrow0}\left\vert
\frac{\sin(\varepsilon)}{1}\right\vert =0,
\]

\item $\ln(\varepsilon)=o(\varepsilon^{-1}),\,\ \ \,\varepsilon\rightarrow0,$ \ \ because

\[
\lim_{\varepsilon\rightarrow0}\left\vert \frac{\ln(\varepsilon)}
{\varepsilon^{-1}}\right\vert =\lim_{\varepsilon\rightarrow0}\left\vert
\frac{\varepsilon^{-1}}{\varepsilon^{-2}}\right\vert =\lim_{\varepsilon
\rightarrow0}\varepsilon=0.
\]
\end{enumerate}
When we apply perturbation methods, we usually use a whole sequence of
comparison functions. In examples one and two we used the sequence
\[
\{\delta_{n}(\varepsilon)=\varepsilon^{n}\}_{n=1}^{\infty},
\]
and in example three we used the sequence
\[
\{\delta_{n}(\varepsilon)=\varepsilon^{n}\}_{n=-1}^{\infty}.
\]
What is characteristic about these sequences is that
\begin{equation}
\delta_{n+1}(\varepsilon)=o(\delta_{n}(\varepsilon)),\,\ \ \,\varepsilon
\rightarrow0,\label{Eq3.7}
\end{equation}
for all $n$ in the range defining the sequences. Sequences of functions that
satisfy conditions (\ref{Eq3.7}) are calles \textit{asymptotic sequencs}.

Here are some asymptotic sequences

\begin{enumerate}
\item $\delta_{n}(\varepsilon)=\sin(\varepsilon)^{n}$,

\item $\delta_{n}(\varepsilon)=\ln(\varepsilon)^{-n}$,

\item $\delta_{n}(\varepsilon)=(\sqrt{\varepsilon})^{n}$.
\end{enumerate}

\noindent Using the notion of asymptotic sequences, we can define asymptotic expansion
analogous to the way inifinite series are defined in elementary calculus

\begin{definition}
Let $\{\delta_{n}(\varepsilon)\}$ be an asymptotic sequence. Then a formal
series
\begin{equation}
{\displaystyle\sum_{n=1}^{\infty}}a_{n}\delta_{n}(\varepsilon)\label{Eq3.9}
\end{equation}

\end{definition}

is an asymptotic expansion for a function $f(\varepsilon)$ as $\varepsilon$
approaches zero if
\begin{equation}
f(\varepsilon)-{\displaystyle\sum_{n=1}^{N}}a_{n}\delta_{n}(\varepsilon)=o(\delta_{N}(\varepsilon)),\text{ \ \ \ }\varepsilon\rightarrow0.\label{Eq3.10}
\end{equation}
Observe that
\begin{align}
f(\varepsilon)-a_{1}\delta_{1}(\varepsilon)  &  =o(\delta_{1}(\varepsilon
)),\text{ \ \ \ }\varepsilon\rightarrow0,\nonumber\\
&  \Downarrow\nonumber\\
\lim_{\varepsilon\rightarrow0}\left\vert \frac{f(\varepsilon)-a_{1}\delta
_{1}(\varepsilon)}{\delta_{1}(\varepsilon)}\right\vert  &  =0,\nonumber\\
&  \Downarrow\nonumber\\
\lim_{\varepsilon\rightarrow0}\left\vert a_{1}-\frac{f(\varepsilon)}
{\delta_{1}(\varepsilon)}\right\vert  &  =0,\nonumber\\
&  \Downarrow\nonumber\\
a_{1}  &  =\lim_{\varepsilon\rightarrow0}\frac{f(\varepsilon)}{\delta
_{1}(\varepsilon)}.\label{Eq3.12}
\end{align}
In an entirely similar way we find that for all $m\geqq1$ that%
\begin{equation}
a_{m}=\lim_{\varepsilon\rightarrow0}\left\vert \frac{f(\varepsilon)-{\displaystyle\sum_{n=1}^{m-1}}a_{n}\delta_{n}(\varepsilon)}{\delta_{m}(\varepsilon)}\right\vert.\label{Eq3.13}
\end{equation}
This shows that for a fixed asymptotic sequence, the coefficients of the
asymptotic expansion for a function $f(\varepsilon)$ are determined by taking
limits. Observe that formula (\ref{Eq3.13}) does not require
differentiability for $f(\varepsilon)$ at $\varepsilon=0$. This is very
different from Taylor expansions which requires that $f(\varepsilon)$ is
infinitely differentiable at $\varepsilon=0$.

This is a hint that asymptotic expansions are much more general than the usual
convergent expansions, for example power series, that we are familiar with from
elementary calculus. In fact, asymptotic expansions may well diverge, but this
does not make them less useful! The following example was first discussed by
Leonard Euler in 1754.

\subsection*{Eulers example}

Let a function $f(\varepsilon)$ be defined by the formula
\begin{equation}
f(\varepsilon)=\int_{0}^{\infty}dt\frac{e^{-t}}{1+\varepsilon t}.\label{Eq3.14}
\end{equation}
The integral defining $f(\varepsilon)$ converge very fast, and because of this
$f(\varepsilon)$ is a very smooth function, in fact it is infinitely smooth
and moreover analytic in the complex plane where the negative real axis has been removed.

Using the properties of telescoping series we observe that for all $m\geqq0$
\begin{equation}
\frac{1}{1+\varepsilon t}={\displaystyle\sum_{n=0}^{m}}(-\varepsilon t)^{n}+\frac{(-\varepsilon t)^{m+1}}{1+\varepsilon t}.\label{Eq3.15}
\end{equation}
Inserting (\ref{Eq3.15}) into (\ref{Eq3.14}) we find that
\begin{equation}
f(\varepsilon)=S_{m}(\varepsilon)+R_{m}(\varepsilon),\label{Eq3.16}
\end{equation}
where
\begin{align}
S_{m}(\varepsilon) &  =
{\displaystyle\sum_{n=0}^{m}}(-1)^{n}n!\varepsilon^{n},\nonumber\\
R_{m}(\varepsilon) &  =(-\varepsilon)^{m+1}\int_{0}^{\infty}dt\frac
{t^{m+1}e^{-t}}{1+\varepsilon t}.
\end{align}
For the quantity $R_{m}(\varepsilon)$ we have the estimate
\begin{equation}
\left\vert R_{m}(\varepsilon)\right\vert \leqq\varepsilon^{m+1}\int_{0}^{\infty}dt\frac{t^{m+1}e^{-t}}{1+\varepsilon t}\leqq\varepsilon^{m+1}
\int_{0}^{\infty}dt\smallskip t^{m+1}e^{-t}=(m+1)!\varepsilon^{m+1},\label{Eq3.18}
\end{equation}
from which it follows that
\begin{equation}
\lim_{\varepsilon\rightarrow0}\left\vert \frac{R_{m}(\varepsilon)}{\varepsilon^{m}}\right\vert \leqq\lim_{\varepsilon\rightarrow0}(m+1)!\varepsilon=0.\label{Eq3.19}
\end{equation}
Thus we have proved that an asymptotic expansion for $f(\varepsilon)$ is
\begin{equation}
f(\varepsilon)={\displaystyle\sum_{n=0}^{\infty}}(-1)^{n}n!\varepsilon^{n}.\label{Eq3.20}
\end{equation}
It is on the other hand trivial to verify that the formal power series
\begin{equation}
{\displaystyle\sum_{n=0}^{\infty}}(-1)^{n}n!\varepsilon^{n},\label{Eq3.21}
\end{equation}
diverge for all $\varepsilon\neq0$!

\noindent In figure 1, we compare the function $f(\varepsilon)$ with what we get from the
asymptotic expansion for a range of $\varepsilon$ and several truncation
levels for the expansion. From this example we make the following two
observations that are quite generic with regards to the convergence or
divergence of asymptotic expansions.%

\begin{figure}[htbp]
\centering
\includegraphics[
natheight=2.460400in,
natwidth=3.766300in,
height=2.501in,
width=3.8147in
]{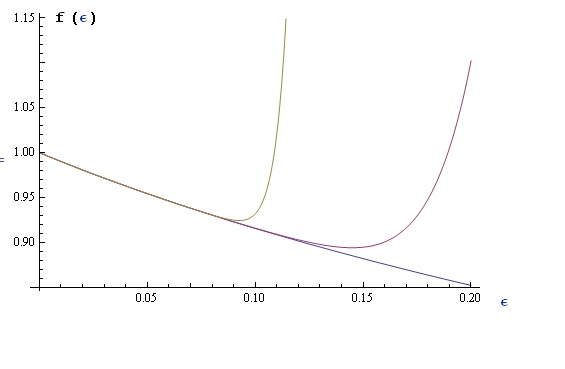}
\caption{Comparing the exact(blue) expression for $f(\varepsilon)$ with the
asymptotic expansion (\ref{Eq3.20}) containing ten(red) and twenty(yellow)
terms}
\label{fig1}
\end{figure}

\noindent Firstly, the asymptotic expansion (\ref{Eq3.20}) is an accurate representation
of $f(\varepsilon)$ in the limit when $\varepsilon$ approaches zero even if
the expansions is divergent. Secondly, adding more terms to the expansion for
a fixed value of $\varepsilon$ makes the expansion less accurate.

In reality we are most of the time, because of algebraic complexity, only
able to calculate a few terms of an asymptotic expansion. Thus convergence
properties of the expansion are most of the time unknown. As this example
shows, convergence properties are also not relevant for what we are trying to
achieve when we solve problems using perturbation methods.

\section{Regular perturbation expansions for ODEs.}

It is now finally time to start solving differential equations using
asymptotic expansions. Let us start with a simple boundary value problem for a
first order ordinary differential equation.

\subsection*{Example 1: A weakly nonlinear boundary value problem.}

Consider the following boundary value problem
\begin{align}
y^{\prime}(x)+y(x)+\varepsilon y^{2}(x)  &  =x,\text{ \ }0<x<1,\nonumber\\
y(1)  &  =1,\label{Eq3.22}
\end{align}
where $\varepsilon$ as usual is a small number. Since the differential
equation is nonlinear and nonseparable, this is a nontrivial problem. The
unperturbed problem is
\begin{align}
y^{\prime}(x)+y(x)  &  =x \text{ \ }0<x<1,\nonumber\\
y(1)  &  =1.\label{Eq3.23}
\end{align}
The unperturbed problem is easy to solve since the equation is a first order
linear equation. The general solution to the equation is
\begin{equation}
y(x)=x-1+Ae^{-x}.\label{Eq3.24}
\end{equation}
The arbitrary constant $A$ is determined from the boundary condition
\begin{align}
y(1)  &  =1,\nonumber\\
&  \Downarrow\nonumber\\
1-1+Ae^{-1}  &  =1,\nonumber\\
&  \Downarrow\nonumber\\
A  &  =e.\label{Eq3.25}
\end{align}
Thus the unique solution to the unperturbed problem is
\begin{equation}
y_{0}(x)=x-1+e^{1-x}.\label{Eq3.26}
\end{equation}
We now want to find an asymptotic expansion for the solution to the perturbed
problem (\ref{Eq3.22}), starting from the solution $y_{0}(x)$. We thus
postulate an expansion of the form
\begin{equation}
y(x;\varepsilon)=y_{0}(x)+\varepsilon y_{1}(x)+\varepsilon^{2}y_{2}
(x)+...\;\;.\label{Eq3.27}
\end{equation}
Inserting (\ref{Eq3.27}) into (\ref{Eq3.22}) and expanding we get
\begin{gather}
(y_{0}+\varepsilon y_{1}+\varepsilon^{2}y_{2}+...)^{\prime}+y_{0}+\varepsilon
y_{1}+\varepsilon^{2}y_{2}+...\nonumber\\
+\varepsilon(y_{0}+\varepsilon y_{1}+\varepsilon^{2}y_{2}+...)^{2}=x,\nonumber\\
\Downarrow\nonumber\\
y_{0}^{\prime}+\varepsilon y_{1}^{\prime}+\varepsilon^{2}y_{2}^{\prime
}+...+y_{0}+\varepsilon y_{1}+\varepsilon^{2}y_{2}+...\nonumber\\
\varepsilon(y_{0}^{2}+2\varepsilon y_{0}y_{1}+..)=x,\nonumber\\
\Downarrow\nonumber\\
y_{0}^{\prime}+y_{0}+\varepsilon(y_{1}^{\prime}+y_{1}+y_{0}^{2})+\varepsilon
^{2}(y_{2}^{\prime}+y_{2}+2y_{0}y_{1})+...=x.\label{Eq3.28}
\end{gather}
We must also expand the boundary condition
\begin{equation}
y_{0}(0)+\varepsilon y_{1}(0)+\varepsilon^{2}y_{2}(0)+...=1.\label{Eq3.29}
\end{equation}
From (\ref{Eq3.28}) and (\ref{Eq3.29}) we get the following perturbation
hierarchy
\begin{align}
y_{0}^{\prime}(x)+y_{0}(x)  &  =x,\nonumber\\
y_{0}(1)  &  =1,\nonumber\\
& \nonumber\\
y_{1}^{\prime}(x)+y_{1}(x)  &  =-y_{0}^{2}(x),\nonumber\\
y_{1}(1)  &  =0,\nonumber\\
& \nonumber\\
y_{2}^{\prime}(x)+y_{2}(x)  &  =-2y_{0}(x)y_{1}(x),\nonumber\\
y_{2}(1)  &  =0.\label{Eq3.30}
\end{align}
We observe that the perturbation hierarch has the special structure that we
have noted earlier. All equations in the hierarchy are determined by the
linear operator $\mathcal{L}=\frac{d}{dx}+1$. The first boundary value problem
in the hierarchy has already been solved. The second equation in the hierarchy
is
\begin{equation}
y_{1}^{\prime}(x)+y_{1}(x)=-y_{0}^{2}(x).\label{Eq3.31}
\end{equation}
Finding a special solution to this equation is simple
\begin{gather}
y_{p}^{\prime}(x)+y_{p}(x)=-y_{0}^{2}(x),\nonumber\\
\Downarrow\nonumber\\
(y_{p}(x)e^{x})^{\prime}=-y_{0}^{2}(x)e^{x},\nonumber\\
\Downarrow\nonumber\\
y_{p}(x)=-e^{-x}\int_{0}^{x}dx^{\prime}e^{x^{\prime}}y_{0}^{2}(x^{\prime
}).\label{Eq3.32}
\end{gather}
Adding a general solution to the homogenous equation, we get the general solution to
equation (\ref{Eq3.31}) in the form
\begin{equation}
y_{1}(x)=A_{1}e^{-x}-e^{-x}\int_{0}^{x}dx^{\prime}e^{x^{\prime}}y_{0}^{2}(x^{\prime
}).\label{Eq3.32.1}
\end{equation}
Inserting the expression for $y_0(x)$ from (\ref{Eq3.26}) into (\ref{Eq3.32.1}), expanding and doing the required integrals, we find that after applying the boundary condition, $y_1(1)=0$, we get 
\begin{equation}
y_{1}(x)=  -x^{2}+4x-5+(2x-x^{2})e^{1-x}+e^{2-2x}.
\label{Eq3.33}
\end{equation}
 We can thus conclude that the perturbation expansion to first order in
$\varepsilon$ is
\begin{equation}
y(x;\varepsilon)=x-1+e^{1-x}+\varepsilon\left(  -x^{2}+4x-5+(2x-x^{2}
)e^{1-x}+e^{2-2x}\right)  +...\;\;.\label{Eq3.34}
\end{equation}
The general solution to the third equation in the perturbation hierarchy is in
a similar way found to be
\begin{equation}
y_{2}(x)=A_{2}e^{-x}-2e^{-x}\int_{0}^{x}dx^{\prime}e^{x^{\prime}}
y_{0}(x^{\prime})y_{1}(x^{\prime}).\label{Eq3.35}
\end{equation}
 The integral in (\ref{Eq3.35})
will have fifteen terms that needs to be integrated. We thus see that even for this very simple
example the algebraic complexity grows quickly. \ 

Recall that we are only ensured that the correction $\varepsilon y_{1}(t)$ is
small with respect to the unperturbed solution $y_{0}(t)$ in the limit when
$\varepsilon$ approaches zero. The perturbation method does not say anything
about the accuracy for any finite value of $\varepsilon$. The hope is of
course that the perturbation expansion also gives a good approximation for
some range of $\varepsilon>0$.

\begin{figure}[htbp]
\centering
\includegraphics[
natheight=2.460400in,
natwidth=3.666800in,
height=2.501in,
width=3.7144in
]{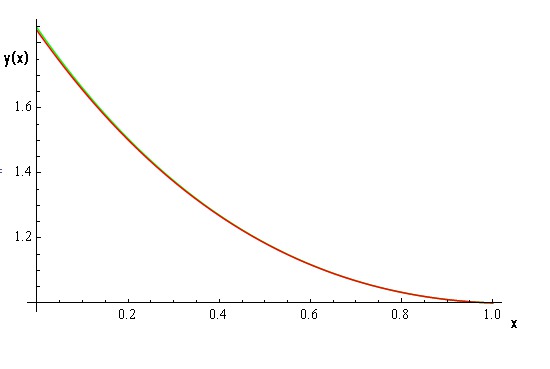}
\caption{Comparing direct perturbation expansion(red) to the exact solution(green),of the boundary value problem.}
\label{fig2}
\end{figure}

Our original equation (\ref{Eq3.22}) is a Ricatti equation and an exact solution to the boundary value problem can be found in terms of Airy functions.
In figure \ref{fig2} we compare our perturbation expansion to the exact solution in the domain $0<x<1$. We observe that even for
$\varepsilon$ as large as $0.05$ our perturbation expansion give a very
accurate representation of the solution over the whole domain.

\noindent In general, we will not have an exact solution that can be used to
investigate the accuracy of the perturbation expansion for finite values of
$\varepsilon$. For example, if our original equation contained $y^3$ instead of 
$y^2$, an exact solution can not be found. This is the normal situation when we apply perturbation
methods. The only way to get at the accuracy of the perturbation expansion is
to compare it to an approximate solution found by some other, independent,
approximation scheme. Often this involve numerical methods, but it could also
be another perturbation method.

As the next example show, things does not always work out as well as in the current example.

\subsection*{Example 2: A weakly damped linear oscillator.}

Consider the following initial value problem
\begin{align}
y^{\prime\prime}(t)+\varepsilon y^{\prime}(t)+y(t)  &  =0,\text{
\ }t>0,\nonumber\\
y(0)  &  =1,\nonumber\\
y^{\prime}(0)  &  =0.\label{Eq3.36}
\end{align}
This is our perturbed problem $\mathcal{P}(\varepsilon)$. The unperturbed
problem, $\mathcal{P}(0)$, is
\begin{align}
y^{\prime\prime}(t)+y(t)  &  =0,\nonumber\\
y(0)  &  =1,\nonumber\\
y^{\prime}(0)  &  =0.\label{Eq3.37}
\end{align}
The general solution to the unperturbed equation is evidently
\begin{equation}
y_{0}(t)=A_{0}e^{it}+A_{0}^{\ast}e^{-it},\label{Eq3.38}
\end{equation}
and the initial condition is satisfied if
\begin{align}
A_{0}+A_{0}^{\ast}  &  =1,\nonumber\\
iA_{0}-iA_{0}^{\ast}  &  =0,\label{Eq3.39}
\end{align}
which has the unique solution $A_{0}=\frac{1}{2}$. Thus the unique solution to
the unperturbed problem is
\begin{equation}
y_{0}(t)=\frac{1}{2}e^{it}+(\ast),\label{Eq3.40}
\end{equation}
where $z+(\ast)$ means $z+z^{\ast}$. This is a very common notation.

We want to find a perturbation expansion for the solution to the perturbed
problem, starting with the solution $y_{0}$ of the unperturbed problem. The
simplest approach is to use an expansion of the form
\begin{equation}
y(t)=y_{0}(t)+\varepsilon y_{1}(t)+\varepsilon^{2}y_{2}(t)...\;\;.\label{Eq3.41}
\end{equation}
We now, as usual, insert (\ref{Eq3.41}) into the perturbed equation
(\ref{Eq3.36}) and expand
\begin{gather}
(y_{0}+\varepsilon y_{1}+\varepsilon^{2}y_{2}+...)^{\prime\prime}\nonumber\\
+\varepsilon(y_{0}+\varepsilon y_{1}+\varepsilon^{2}y_{2}+...)^{\prime}
+y_{0}+\varepsilon y_{1}+\varepsilon^{2}y_{2}+...=0,\nonumber\\
\Downarrow\nonumber\\
y_{0}^{\prime\prime}+y_{0}+\varepsilon(y_{1}^{\prime\prime}+y_{1}
+y_{0}^{\prime})+\varepsilon^{2}(y_{2}^{\prime\prime}+y_{2}+y_{1}^{\prime
})+...=0.\label{Eq3.42}
\end{gather}
We must in a similar way expand the initial conditions
\begin{align}
y_{0}(0)+\varepsilon y_{1}(0)+\varepsilon^{2}y_{2}(0)+...  &  =1,\nonumber
\\
y_{0}^{\prime}(0)+\varepsilon y_{1}^{\prime}(0)+\varepsilon^{2}y_{2}^{\prime
}(t)+...  &  =0.\label{Eq3.43}
\end{align}
From equations (\ref{Eq3.42}) and (\ref{Eq3.43}) we get the following
perturbation hierarchy%
\begin{align}
y_{0}^{^{\prime\prime}}+y_{0}  &  =0,\text{ \ }t>0,\nonumber\\
y_{0}(0)  &  =1,\nonumber\\
y_{0}^{\prime}(0)  &  =0,\nonumber\\
& \nonumber\\
y_{1}^{\prime\prime}+y_{1}  &  =-y_{0}^{\prime},\text{ \ }t>0,\nonumber\\
y_{1}(0)  &  =0,\nonumber\\
y_{1}^{\prime}(0)  &  =0,\nonumber\\
& \nonumber\\
y_{2}^{\prime\prime}+y_{2}  &  =-y_{1}^{\prime},\text{ \ }t>0,\nonumber\\
y_{2}(0)  &  =0,\nonumber\\
y_{2}^{\prime}(0)  &  =0.\label{Eq3.44}
\end{align}
We note that the perturbation hierarchy has the special form discussed
earlier. Here the linear operator determining the hierarchy is $L=\frac{d^{2}
}{dt^{2}}+1$.

The first initial value problem in the hierarchy has already been solved. The
solution is (\ref{Eq3.40}). Inserting $y_{0}(t)$ into the second equation in
the hierarchy we get
\begin{equation}
y_{1}^{\prime\prime}+y_{1}=-\frac{i}{2}e^{it}+(\ast).\label{Eq3.45}
\end{equation}
Looking for particular solutions of the form
\[
y_{1}^{p}(t)=Ce^{it}+(\ast),
\]
will not work, here because the right-hand side of (\ref{Eq3.45}) is a solution
to the homogenous equation. In fact (\ref{Eq3.45}) is a harmonic oscillator
driven on ressonance. For such cases we must rather look for a special
solution of the form
\begin{equation}
y_{1}^{p}(t)=Cte^{it}+(\ast).\label{Eq3.46}
\end{equation}
By inserting (\ref{Eq3.46}) into (\ref{Eq3.45}) we find $C=-\frac{1}{4}$. The
general solution to equation (\ref{Eq3.45}) is then
\begin{equation}
y_{1}(t)=A_{1}e^{it}-\frac{1}{4}te^{it}+(\ast).\label{Eq3.47}
\end{equation}
Applying the initial condition for $y_{1}(t)$ we easily find that
$A_{1}=-\frac{i}{4}$. Thus the perturbation expansion to first order in
$\varepsilon$ is
\begin{equation}
y(t)=\frac{1}{2}e^{it}+\varepsilon\frac{1}{4}(i-t)e^{it}+(\ast).\label{Eq3.48}
\end{equation}
Let $y_{E}(t)$ be a high precision numerical solution to the perturbed problem
(\ref{Eq3.36}). For $\varepsilon=0.01$ we get for increasing time

\bigskip

$
\begin{array}
[c]{ccc}
t & y_{E} & y\\
4 & -0.6444 & -0.6367\\
40 & -0.5426 & -0.5372\\
400 & -0.0722 & 0.5295
\end{array}
$

\bigskip

\begin{figure}[ptb]
\centering
\includegraphics[
natheight=3.844100in,
natwidth=5.405900in,
height=3.4765in,
width=4.8801in
]{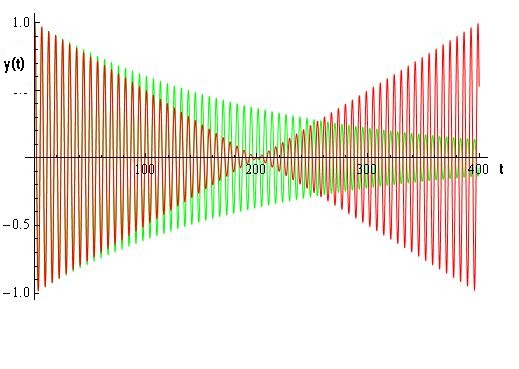}
\caption{Comparing the direct perturbation expansion(red) and a high precision
numerical solution(green)}
\label{fig3}
\end{figure}

\noindent The solution starts out by being quite accurate, but as $t$ increases, the
perturbation expansion eventually looses any relation to the exact solution.
The true extent of the disaster is seen in figure \ref{fig3}.

So what is going on, why is the perturbation expansion such a bad
approximation in this example?

Observe that $y_{1}$ contain a term that is proportional to $t$. Thus as $t$
grows the size of $y_{1}$ also grows and when
\begin{equation}
t\sim\frac{1}{\varepsilon}\label{Eq3.49}
\end{equation}
the second term in the perturbation expansion become as large as the first
term. The ordering of the expansion breaks down and the first correction,
$\varepsilon y_{1}$, is of the same size as the solution to the unperturbed
problem, $y_{0}$.

The reason why the growing term, $y_{1}$, is a problem here, but was not a
problem in the previous example, is that here the domain for the independent
variable is unbounded.

\noindent Let us at this point introduce some standard terminology. The last two
examples involved perturbation expansions where the coefficients depended on a
parameter. In general such expansions takes the form
\begin{equation}
f(\varepsilon;\mathbf{x})\sim{\displaystyle\sum_{n=1}^{\infty}}
a_{n}(\mathbf{x})\delta_{n}(\varepsilon),\text{ \ \ }\varepsilon\rightarrow0\label{Eq3.50}
\end{equation}
where the parameter, $\mathbf{x}$, ranges over some domain $V\subset
\mathbb{R}^{m}$ for some $m$. In example one, $V$ is the interval $[0,1]$ whereas in
example two, $V$ is the unbounded interval $(0,\infty)$.

With the introduction of a parameter dependence of the coefficients, a
breakdown of order in the expansion for some region(s) in $V$ becomes a
possibility. We saw how this came about for the case of the damped harmonic oscillator model (\ref{Eq3.36}).

And let me be clear about this: Breakdown of order in parameter dependent
perturbation expansions is not some weird, rarely occuring, event. On the
contrary it is very common.

Thus methods has to be invented to handle this fenomenon, which is called
\textit{nonuniformity} of asymptotic expansions. The multiple scale method is
design to do exactly this.

\section{The multiple scale method for weakly nonlinear scalar ODEs and systems of ODEs.}

In the previous section we saw that trying to represent the solution to the
problem
\begin{align}
y^{\prime\prime}(t)+\varepsilon y^{\prime}(t)+y(t)  &  =0,\text{
\ }t>0,\nonumber\\
y(0)  &  =1,\nonumber\\
y^{\prime}(0)  &  =0,\label{Eq4.1}
\end{align}
using a regular perturbation expansion
\begin{equation}
y(t)=y_{0}(t)+\varepsilon y_{1}(t)+\varepsilon^{2}y_{2}(t)...\;\;,\label{Eq4.2}
\end{equation}
leads to a nonuniform expansion where ordering of the terms broke down for
$t\sim\frac{1}{\varepsilon}$. In order to understand how to fix this, let us
have a look at the exact solution to (\ref{Eq4.1}). The exact solution can be
found using characteristic polynomials. We get
\begin{equation}
y(t)=Ce^{-\frac{1}{2}\varepsilon t}e^{i\sqrt{1-\frac{1}{4}\varepsilon^{2}}t}+(\ast),\label{Eq4.3}
\end{equation}
where
\begin{equation}
C=\frac{-\lambda^{\ast}}{\lambda-\lambda^{\ast}},\text{ \ \ }\lambda=-\frac
{1}{2}\varepsilon+i\sqrt{1-\frac{1}{4}\varepsilon^{2}}.\label{Eq4.4}
\end{equation}
If we expand the square root in the exponent with respect to $\varepsilon$, we
get
\begin{equation}
y(t)\approx Ce^{-\frac{1}{2}\varepsilon t}e^{it}e^{-\frac{i}{8}\varepsilon
^{2}t}+(\ast).\label{Eq4.5}
\end{equation}
Observe that if $f(\xi)$ is a function whose derivative is of order one, then
the function
\begin{equation}
g_{n}(t)=f(\varepsilon^{n}t),\label{Eq4.6}
\end{equation}
satisfy
\begin{equation}
\bigtriangleup g_{n}(t)=g_{n}(t+T)-g_{n}(t)\approx\varepsilon^{n}f^{\prime
}(\varepsilon^{n}t)T=O(1)\Longleftrightarrow T\sim\varepsilon^{-n}.\label{Eq4.7}
\end{equation}
We express this by saying that the function $g_{n}(t)$ \textit{vary on the
time scale} $t_{n}=\varepsilon^{-n}t$. If we now look at equation
(\ref{Eq4.5}), we see that the approximate solution (\ref{Eq4.5}) vary on three
separate time scales $t_{0}=\varepsilon^{0}t,t_{1}=\varepsilon^{-1}t$ and
$t_{2}=\varepsilon^{-2}t$. If we include more terms in the Taylor expansion
for the square root in (\ref{Eq4.3}) the resulting solution will depend on
even more time scales.

Inspired by this example we postulate the existence of a function
\begin{equation}
h=h(t_{0},t_{1},t_{2},...),\label{Eq4.8}
\end{equation}
such that%
\begin{equation}
y(t)=h(t_{0},t_{1},t_{2},...)|_{t_{j}=\varepsilon^{j}t},\label{Eq4.9}
\end{equation}
is a solution to problem (\ref{Eq4.1}). Using the chain rule we evidently
have
\[
\frac{dy}{dt}(t)=\left\{  (\partial_{t_{0}}+\varepsilon\partial_{t_{1}
}+\varepsilon^{2}\partial_{t_{2}}+...)h\right\}  |_{t_{j}=\varepsilon^{j}t},
\]
which we formally write as
\begin{equation}
\frac{d}{dt}=\partial_{t_{0}}+\varepsilon\partial_{t_{1}}+\varepsilon
^{2}\partial_{t_{2}}+...\;\;.\label{Eq4.10}
\end{equation}
The function $h$ is represented using a perturbation expansion of the form
\begin{equation}
h=h_{0}+\varepsilon h_{1}+\varepsilon^{2}h_{2}+...\;\;.\label{Eq4.11}
\end{equation}
The multiple scale method now proceed by substituting (\ref{Eq4.10}) and
(\ref{Eq4.11}) into the differential equation
\begin{equation}
y^{\prime\prime}(t)+\varepsilon y^{\prime}(t)+y(t)=0,\label{Eq4.12}
\end{equation}
and expanding everything in sight.
\begin{gather}
(\partial_{t_{0}}+\varepsilon\partial_{t_{1}}+\varepsilon^{2}\partial_{t_{2}
}+...)(\partial_{t_{0}}+\varepsilon\partial_{t_{1}}+\varepsilon^{2}
\partial_{t_{2}}+...)\nonumber\\
(h_{0}+\varepsilon h_{1}+\varepsilon^{2}h_{2}+...)+\varepsilon(\partial
_{t_{0}}+\varepsilon\partial_{t_{1}}+\varepsilon^{2}\partial_{t_{2}
}+...)\nonumber\\
(h_{0}+\varepsilon h_{1}+\varepsilon^{2}h_{2}+...)+h_{0}+\varepsilon
h_{1}+\varepsilon^{2}h_{2}+...=0,\nonumber\\
\Downarrow\nonumber\\
(\partial_{t_{0}t_{0}}+\varepsilon(\partial_{t_{0}t_{1}}+\partial_{t_{1}t_{0}
})+\varepsilon^{2}(\partial_{t_{0}t_{2}}+\partial_{t_{1}t_{1}}+\partial
_{t_{2}t_{0}})+...)\nonumber\\
(h_{0}+\varepsilon h_{1}+\varepsilon^{2}h_{2}+...)+\varepsilon(\partial
_{t_{0}}+\varepsilon\partial_{t_{1}}+\varepsilon^{2}\partial_{t_{2}
}+...)\nonumber\\
(h_{0}+\varepsilon h_{1}+\varepsilon^{2}h_{2}+...)++h_{0}+\varepsilon
h_{1}+\varepsilon^{2}h_{2}+...=0,\nonumber\\
\Downarrow\nonumber\\
\partial_{t_{0}t_{0}}h_{0}+h_{0}+\varepsilon(\partial_{t_{0}t_{0}}h_{1}
+h_{1}+\partial_{t_{0}t_{1}}h_{0}+\partial_{t_{1}t_{0}}h_{0}+\partial_{t_{0}
}h_{0})\nonumber\\
+\varepsilon^{2}(\partial_{t_{0}t_{0}}h_{2}+h_{2}+\partial_{t_{0}t_{1}}
h_{1}+\partial_{t_{1}t_{0}}h_{1}+\partial_{t_{0}t_{2}}h_{0}+\partial
_{t_{1}t_{1}}h_{0}\nonumber\\
+\partial_{t_{2}t_{0}}h_{0}+\partial_{t_{1}}h_{0}+\partial_{t_{0}}
h_{1})+...=0,\label{Eq4.13}
\end{gather}
which gives us the following perturbation hierarchy to second order in
$\varepsilon$
\begin{align}
\partial_{t_{0}t_{0}}h_{0}+h_{0}  &  =0,\nonumber\\
& \nonumber\\
\partial_{t_{0}t_{0}}h_{1}+h_{1}  &  =-\partial_{t_{0}t_{1}}h_{0}
-\partial_{t_{1}t_{0}}h_{0}-\partial_{t_{0}}h_{0},\nonumber\\
& \nonumber\\
\partial_{t_{0}t_{0}}h_{2}+h_{2}  &  =-\partial_{t_{0}t_{1}}h_{1}
-\partial_{t_{1}t_{0}}h_{1}-\partial_{t_{0}t_{2}}h_{0}\nonumber\\
&  -\partial_{t_{1}t_{1}}h_{0}-\partial_{t_{2}t_{0}}h_{0}-\partial_{t_{1}
}h_{0}-\partial_{t_{0}}h_{1}.\label{Eq4.14}
\end{align}
We observe, in passing, that the perturbation hierarchy has the special form we
have seen several times before. Here the common differential operator is
$L=\partial_{t_{0}t_{0}}+1$.

At this point a remark is in order. It is fair to say that there is not a full
agreement among the practitioners of the method of multiple scales about how
to perform these calculations. The question really hinges on whether to take
the multiple variable function $h(t_{0},t_{1},..)$ seriously or not. If you
do, you will be lead to a certain way of doing these calculation. This is the
point of view used in most textbooks on this subject. We will not follow this
path here. We will not take $h$ seriously as a multiple variable function and
never forget that what we actually want is not $h$, but rather $y$, which is
defined in terms of $h$ through equation (\ref{Eq4.9}). This point of view
will lead us to do multiple scale calculations in a different way from what
you see in most textbooks. This way is very efficient and will make it
possible to go to order $\varepsilon^{2}$ and beyond without being overwhelmed
by the amount of algebra that needs to be done.

  What I mean when I say that we
will not take $h$ seriously as a multiple variable function will become clear
as we proceed. One immediate consequence of this choise is already evident from the way I write the
perturbation hierarchy. Observe that I keep
\begin{equation}
\partial_{t_{i}t_{j}}h_{k}\;\;\text{and \ }\partial_{t_{j}t_{i}}h_{k}\;\;,\label{Eq4.15}
\end{equation}
as separate terms, I don't use the equality of cross derivatives to simplify
my expressions. This is the first rule we must follow when we do multiple
scale calculations in the way I am teaching you in these lecture notes. If we
took $h$ seriously as a multiple variable function we would put cross
derivatives equal. The second rule we must follow is to disregard the
initial values for the time being. We will fit the initial values at the very
end of our calculations rather than do it at each order in $\varepsilon$ like
in example 1 and example 2.

Let us now proceed to solve the equations in the perturbation hierarchy. At
order $\varepsilon^{0}$ we have the equation
\begin{equation}
\partial_{t_{0}t_{0}}h_{0}+h_{0}=0.\label{Eq4.16}
\end{equation}
When we are applying multiple scales to ordinary differential equations we
always use the general solution to the order $\varepsilon^{0}$ equation. \ For
partial differential equations this will not be so, as we will see later. The
general solution to (\ref{Eq4.16}) is evidently
\begin{equation}
h_{0}(t_{0},t_{1},..)=A_{0}(t_{1},t_{2},..)e^{it_{0}}+(\ast).\label{Eq4.17}
\end{equation}
Observe that the equation only determines how $h_{0}$ depends on the fastest
time scale $t_{0}$, the dependence on the other time scales $t_{1},t_{2},..$,
is arbitrary at this point and this is reflected in the fact that the  integration
"constant" $A_{0}$ is actually a function depending on $t_{1},t_{2},..$\;.

We have now solved the order $\varepsilon^{0}$ equation. Inserting the
expression for $h_{0}$ into the order $\varepsilon$ equation, we get after
some simple algebra
\begin{equation}
\partial_{t_{0}t_{0}}h_{1}+h_{1}=-2i(\partial_{t_{1}}A_{0}+\frac{1}{2}
A_{0})e^{it_{0}}+(\ast).\label{Eq4.18}
\end{equation}
We now need a particular solution to this equation. Observe that since $A_{0}
$ only depends on the slow time scales $t_{1},t_{2},..$ equation
(\ref{Eq4.18}) is in fact a harmonic oscillator driven on ressonance. It is
simple to verify that it has a particular solution of the form
\begin{equation}
h_{1}(t_{0},t_{1},..)=-t_{0}(\partial_{t_{1}}A_{0}+\frac{1}{2}A_{0})e^{it_{0}}.\label{Eq4.19}
\end{equation}
But this term is growing and will lead to breakdown of ordering for the
perturbation expansion (\ref{Eq4.11}) when $t_{0}\sim\varepsilon^{-1}$. This
breakdown was exactly what we tried to avoid using the multiple scales approach!

But everything is not lost, we now have freedom to remove the growing term by
postulating that
\begin{equation}
\partial_{t_{1}}A_{0}=-\frac{1}{2}A_{0}.\label{Eq4.20}
\end{equation}
With this choise, the order $\varepsilon$ equation simplifies into
\begin{equation}
\partial_{t_{0}t_{0}}h_{1}+h_{1}=0.\label{Eq4.21}
\end{equation}
Terms in equations leading to linear growth like in (\ref{Eq4.19}), are
traditionally called \textit{secular terms}. The name are derived from the
Latin word soeculum that means century and are used here because this kind of
nonuniformity was first observed on century time scales in planetary orbit calculations.

At this point we introduce the third rule for doing multiple scale
calculations in the particular way that I advocate in these lecture notes. The
rule is to disregard the general solution of the homogenous equation for all
equations in the perturbation hierarchy except the first. \ We therefore
choose $h_{1}=0$ and proceed to the order $\varepsilon^{2}$ equation using
this choise. The equation for $h_{2}$ then simplifies into
\begin{equation}
\partial_{t_{0}t_{0}}h_{2}+h_{2}=-2i(\partial_{t_{2}}A_{0}-\frac{i}{2}
\partial_{t_{1}t_{1}}A_{0}-\frac{i}{2}\partial_{t_{1}}A_{0})e^{it_{0}}
+(\ast).\label{Eq4.22}
\end{equation}
We have a new secular term and in order to remove it we must postulate that
\begin{equation}
\partial_{t_{2}}A_{0}=\frac{i}{2}\partial_{t_{1}t_{1}}A_{0}+\frac{i}
{2}\partial_{t_{1}}A_{0}.\label{Eq4.23}
\end{equation}
Using this choise, our order $\varepsilon^{2}$ equation simplifies into
\begin{equation}
\partial_{t_{0}t_{0}}h_{2}+h_{2}=0.\label{Eq4.24}
\end{equation}
For this equation we use, according to the rules of the game, the special
solution $h_{2}=0$.

What we have found so far is then
\begin{equation}
h(t_{0},t_{1},t_{2},..)=A_{0}(t_{1},t_{2},..)e^{it_{0}}+(\ast)+O(\varepsilon
^{3}),\label{Eq4.25}
\end{equation}
where
\begin{align}
\partial_{t_{1}}A_{0}  &  =-\frac{1}{2}A_{0},\label{Eq4.26}\\
\partial_{t_{2}}A_{0}  &  =\frac{i}{2}\partial_{t_{1}t_{1}}A_{0}+\frac{i}
{2}\partial_{t_{1}}A_{0}.\label{Eq4.27}
\end{align}
At this point you might ask if we actually have done something useful. Instead
of one ODE we have ended up with two coupled partial differential equations,
and clearly, if we want to go to higher order we will get even more partial
differential equations.

Observe that if we use (\ref{Eq4.26}) we can simplify equation (\ref{Eq4.27})
by removing the derivatives on the right hand side. Doing this we get the
system
\begin{align}
\partial_{t_{1}}A_{0}  &  =-\frac{1}{2}A_{0},\label{Eq4.28}\\
\partial_{t_{2}}A_{0}  &  =-\frac{i}{8}A_{0}.\label{Eq4.29}
\end{align}
The first thing that should come to mind when we see a system like
(\ref{Eq4.28}) and (\ref{Eq4.29}), is that the count is wrong. There is one
unknown function, $A_{0}$, and two equations. The system is
\textit{overdetermined} and will get more so, if we extend our calculations to
higher order in $\varepsilon$. Under normal circumstances, overdetermined
systems of equations has no solutions, which for our setting means that under
normal circumstances the function $h(t_{0},t_{1},t_{2},..)$ does not exist!
This is what I meant when I said that we will not take the functions $h$
seriously as a multiple variable function. For systems of first order partial
differential equations like (\ref{Eq4.28}), (\ref{Eq4.29}) there is a simple
test we can use to decide if a solution actually does  exist. This is the cross
derivative test you know from elementary calculus. Taking $\partial_{t_{2}}$
of equation (\ref{Eq4.28}) and $\partial_{t_{1}}$ of equation (\ref{Eq4.29})
we get
\begin{align}
\partial_{t_{2}t_{1}}A_{0}  &  =\partial_{t_{2}}\partial_{t_{1}}A_{0}
=-\frac{1}{2}\partial_{t_{2}}A_{0}=\frac{i}{16}A_{0},\nonumber\\
\partial_{t_{1}t_{2}}A_{0}  &  =\partial_{t_{1}}\partial_{t_{2}}A_{0}
=-\frac{i}{8}\partial_{t_{1}}A_{0}=\frac{i}{16}A_{0}.\label{Eq4.30}
\end{align}
According to the cross derivative test the overdetermined system
(\ref{Eq4.28}), (\ref{Eq4.29}) is solvable. Thus in this case the function $h$
exists, at least as a two variable function. To make sure that it exists as a
function of three variables we must derive and solve the perturbation
hierarchy to order $\varepsilon^{3}$, and then perform the cross derivative
test. For the current example we will never get into trouble, the many
variable function $h$ will exist as a function of however many variables we
want. But I want you to reflect on how special this must be. We will at order
$\varepsilon^{n}$ have a system of $n$ partial differential equations for only
one unknown function ! In general we will not be so lucky as in the current example,and
the function $h(t_0,t_1,...)$ will not exist. This fact   is the
reason why we can not take $h$ seriously as a many variable function.

So, should we be disturbed by the nonexistence of the solution to
the perturbation hierarchy in the general
case? Actually no, and the reason is that we do not care about $h(t_{0}
,t_{1},..)$. What we care about is $y(t)$. \ 

Inspired by this let us define an \textit{amplitude}, $A(t)$, by
\begin{equation}
A(t)=A_{0}(t_{1},t_{2},..)|_{t_{j}=\varepsilon^{j}t}.\label{Eq4.31}
\end{equation}
Using this and equations (\ref{Eq4.9}) and  (\ref{Eq4.25}), our perturbation
expansion for $y(t)$ is
\begin{equation}
y(t)=A(t)e^{it}+(\ast)+O(\varepsilon^{3}).\label{Eq4.32}
\end{equation}
For the amplitude $A(t)$ we have, using equations (\ref{Eq4.10}),(\ref{Eq4.28}
),(\ref{Eq4.29}) and (\ref{Eq4.31})
\begin{gather}
\frac{dA}{dt}(t)=\{(\partial_{t_{0}}+\varepsilon\partial_{t_{1}}
+\varepsilon^{2}\partial_{t_{2}}+...)A_{0}(t_{1},t_{2},...)\}|_{t_{j}
=\varepsilon^{j}t},\nonumber\\
\Downarrow\nonumber\\
\frac{dA}{dt}(t)=\{-\varepsilon\frac{1}{2}A_{0}(t_{1},t_{2},...)-\varepsilon
^{2}\frac{i}{8}A_{0}(t_{1},t_{2},...)\}|_{t_{j}=\varepsilon^{j}t},\nonumber\\
\Downarrow\nonumber\\
\frac{dA}{dt}=-\varepsilon\frac{1}{2}A-\varepsilon^{2}\frac{i}{8}
A.\label{Eq4.33}
\end{gather}

\noindent This equation is our first example of an \textit{amplitude equation}. The
amplitude equation determines, through equation (\ref{Eq4.31}), the perturbation
expansion for our solution to the original equation (\ref{Eq4.1}). The
amplitude equation is of course easy to solve and we get
\begin{equation}
y(t)=Ce^{-\frac{1}{2}\varepsilon t}e^{it}e^{-\frac{i}{8}\varepsilon^{2}
t}+(\ast)+O(\varepsilon^{3}).\label{Eq4.34}
\end{equation}
The constant $C$ can be fitted to the initial conditions. What we get is equal
to the exact solution up to second order in $\varepsilon$ as we see by
comparing with (\ref{Eq4.5}).

Let us next apply the multiple scale method to some weakly nonlinear ordinary
differential equations. For these cases no exact solution is known, so the
multiple scale method will actually be useful!

\subsection*{Example 1}

Consider the initial value problem
\begin{align}
\frac{d^{2}y}{dt^{2}}+y  &  =\varepsilon y^{3},\nonumber\\
y(0)  &  =1,\nonumber\\
\frac{dy}{dt}(0)  &  =0.\label{Eq4.35}
\end{align}
If we try do solve this problem using a regular perturbation expansion, we will
get secular terms that will lead to breakdown of ordering on a time scale
$t\sim\varepsilon^{-1}$. Let us therefore apply the multiple scale approach.
We introduce a function $h$ through
\begin{equation}
y(t)=h(t_{0},t_{1},t_{2},...)|_{t_{j}=\varepsilon^{j}t},\label{Eq4.36}
\end{equation}
and expansions
\begin{align}
\frac{d}{dt}  &  =\partial_{t_{0}}+\varepsilon\partial_{t_{1}}+\varepsilon
^{2}\partial_{t_{2}}+...\;\;,\nonumber\\
h  &  =h_{0}+\varepsilon h_{1}+\varepsilon^{2}h_{2}+...\;\;.\label{Eq4.38}
\end{align}
Inserting these expansions into (\ref{Eq4.35}), we get
\begin{gather}
(\partial_{t_{0}}+\varepsilon\partial_{t_{1}}+\varepsilon^{2}\partial_{t_{2}
}+...)(\partial_{t_{0}}+\varepsilon\partial_{t_{1}}+\varepsilon^{2}
\partial_{t_{2}}+...)\nonumber\\
(h_{0}+\varepsilon h_{1}+\varepsilon^{2}h_{2}+...)+h_{0}+\varepsilon
h_{1}+\varepsilon^{2}h_{2}+...\nonumber\\
=\varepsilon(h_{0}+\varepsilon h_{1}+\varepsilon^{2}h_{2}+...)^{3},\nonumber\\
\Downarrow\nonumber\\
(\partial_{t_{0}t_{0}}+\varepsilon(\partial_{t_{0}t_{1}}+\partial_{t_{1}t_{0}
})+\varepsilon^{2}(\partial_{t_{0}t_{2}}+\partial_{t_{1}t_{1}}+\partial
_{t_{2}t_{0}})+...)\nonumber\\
(h_{0}+\varepsilon h_{1}+\varepsilon^{2}h_{2}+...)+h_{0}+\varepsilon
h_{1}+\varepsilon^{2}h_{2}+...\nonumber\\
=\varepsilon h_{0}^{3}+3\varepsilon^{2}h_{0}^{2}h_{1}+...\;\;,\nonumber\\
\Downarrow\nonumber\\
\partial_{t_{0}t_{0}}h_{0}+h_{0}+\varepsilon(\partial_{t_{0}t_{0}}h_{1}
+h_{1}+\partial_{t_{0}t_{1}}h_{0}+\partial_{t_{1}t_{0}}h_{0})\nonumber\\
+\varepsilon^{2}(\partial_{t_{0}t_{0}}h_{2}+h_{2}+\partial_{t_{0}t_{1}}
h_{1}+\partial_{t_{1}t_{0}}h_{1}+\partial_{t_{0}t_{2}}h_{0}+\partial
_{t_{1}t_{1}}h_{0}\nonumber\\
+\partial_{t_{2}t_{0}}h_{0})+...=\varepsilon h_{0}^{3}+3\varepsilon^{2}
h_{0}^{2}h_{1}+...\;\;,\label{Eq4.39}
\end{gather}
which gives us the following perturbation hierarchy to second order in
$\varepsilon$
\begin{align}
\partial_{t_{0}t_{0}}h_{0}+h_{0}  &  =0,\nonumber\\
& \nonumber\\
\partial_{t_{0}t_{0}}h_{1}+h_{1}  &  =h_{0}^{3}-\partial_{t_{0}t_{1}}
h_{0}-\partial_{t_{1}t_{0}}h_{0},\nonumber\\
& \nonumber\\
\partial_{t_{0}t_{0}}h_{2}+h_{2}  &  =3h_{0}^{2}h_{1}-\partial_{t_{0}t_{1}
}h_{1}-\partial_{t_{1}t_{0}}h_{1}-\partial_{t_{0}t_{2}}h_{0}\nonumber\\
&  -\partial_{t_{1}t_{1}}h_{0}-\partial_{t_{2}t_{0}}h_{0}.\label{Eq4.40}
\end{align}
The general solution to the first equation in the perturbation hierarchy is
\begin{equation}
h_{0}=A_{0}(t_{1},t_{2},...)e^{it_{0}}+(\ast).\label{Eq4.41}
\end{equation}
Inserting this into the right hand side of the second equation in the
hierarchy and expanding, we get
\begin{equation}
\partial_{t_{0}t_{0}}h_{1}+h_{1}=(3|A_{0}|^{2}A_{0}-2i\partial_{t_{1}}
A_{0})e^{it}+A_{0}^{3}e^{3it}+(\ast).\label{Eq4.42}
\end{equation}
In order to remove secular terms we must postulate that
\begin{equation}
\partial_{t_{1}}A_{0}=-\frac{3i}{2}|A_{0}|^{2}A_{0}.\label{Eq4.43}
\end{equation}
This choise simplify the equation for $h_{1}$ into
\begin{equation}
\partial_{t_{0}t_{0}}h_{1}+h_{1}=A_{0}^{3}e^{3it_{0}}+(\ast).\label{Eq4.44}
\end{equation}
According to the rules of the game we now need a particular solution to this
equation. It is easy to verify that
\begin{equation}
h_{1}=-\frac{1}{8}A_{0}^{3}e^{3it_{0}}+(\ast),\label{Eq4.45}
\end{equation}
is such a particular solution.

We now insert $h_{0}$ and $h_{1}$ into the right hand side of the third
equation in the perturbation hierarchy and find
\begin{equation}
\partial_{t_{0}t_{0}}h_{2}+h_{2}=(-\frac{3}{8}|A_{0}|^{4}A_{0}-2i\partial
_{t_{2}}A_{0}-\partial_{t_{1}t_{1}}A_{0})e^{it_{0}}+(\ast)+NST,\label{Eq4.46}
\end{equation}
where $NST$ is an acronym for "nonsecular terms". Since we are not here
planning to go beyond second order in $\varepsilon$, we will at this order
only need the secular terms and group the rest into $NST$. In order to remove
the secular terms we must postulate that
\begin{equation}
\partial_{t_{2}}A_{0}=\frac{3i}{16}|A_{0}|^{4}A_{0}+\frac{i}{2}\partial
_{t_{1}t_{1}}A_{0}.\label{Eq4.47}
\end{equation}
As before, it make sense to simplify (\ref{Eq4.47}) using equation
(\ref{Eq4.43}). This leads to the following overdetermined system of equations
for $A_{0}$
\begin{align}
\partial_{t_{1}}A_{0}  &  =-\frac{3i}{2}|A_{0}|^{2}A_{0},\nonumber\\
\partial_{t_{2}}A_{0}  &  =-\frac{15i}{16}|A_{0}|^{4}A_{0}\label{Eq4.48}
\end{align}
Let us check solvability of this system using the cross derivative test
\begin{align*}
\partial_{t_{2}t_{1}}A_{0}  &  =-\frac{3i}{2}\partial_{t_{2}}(A_{0}^{2}
A_{0}^{\ast})\\
&  =-\frac{3i}{2}\left(  2A_{0}\partial_{t_{2}}A_{0}A_{0}^{\ast}+A_{0}
^{2}\partial_{t_{2}}A_{0}^{\ast}\right) \\
&  =-\frac{3i}{2}\left(  2A_{0}\left(  -\frac{15i}{16}|A_{0}|^{4}A_{0}\right)
A_{0}^{\ast}+A_{0}^{2}\left(  \frac{15i}{16}|A_{0}|^{4}A_{0}^{\ast}\right)
\right) \\
&  =-\frac{45}{32}|A_{0}|^{6}A_{0}.
\end{align*}

\begin{align*}
\partial_{t_{1}t_{2}}A_{0}  &  =-\frac{15i}{16}\partial_{t_{1}}\left(
A_{0}^{3}A_{0}^{\ast2}\right) \\
&  =-\frac{15i}{16}\left(  3A_{0}^{2}\partial_{t_{1}}A_{0}A_{0}^{\ast2}
+2A_{0}^{3}A_{0}^{\ast}\partial_{t_{1}}A_{0}^{\ast}\right) \\
&  =-\frac{15i}{16}\left(  3A_{0}^{2}\left(  -\frac{3i}{2}|A_{0}|^{2}
A_{0}\right)  A_{0}^{\ast2}+2A_{0}^{3}A_{0}^{\ast}\left(  \frac{3i}{2}
|A_{0}|^{2}A_{0}^{\ast}\right)  \right) \\
&  =-\frac{45}{32}|A_{0}|^{6}A_{0}.
\end{align*}
The system is compatible, and thus the function $h_{0}$ exists as a function of
two variables. Of course, whether or not $h_{0}$ exists is only of academic
interest for us since our only aim is to find the solution of the original
equation $y(t)$.

  Defining an amplitude, $A(t)$ by
\begin{equation}
A(t)=A_{0}(t_{1},t_{2},...)|_{t_{j}=\varepsilon^{j}t},\label{Eq4.49}
\end{equation}
we find that the solution is
\begin{equation}
y(t)=A(t)e^{it}-\varepsilon\frac{1}{8}A^{3}e^{3it}+(\ast)+O(\varepsilon
^{2}),\label{Eq4.50}
\end{equation}
where the amplitude satisfy the equation
\begin{gather}
\frac{dA}{dt}(t)=\{(\partial_{t_{0}}+\varepsilon\partial_{t_{1}}
+\varepsilon^{2}\partial_{t_{2}}+...)A_{0}(t_{1},t_{2},...)\}|_{t_{j}
=\varepsilon^{j}t},\nonumber\\
\Downarrow\nonumber\\
\frac{dA}{dt}(t)=\{-\varepsilon\frac{3i}{2}|A_{0}|^{2}A_{0}(t_{1}
,t_{2},...)-\varepsilon^{2}\frac{15i}{16}|A_{0}|^{4}A_{0}(t_{1},t_{2}
,...)\}|_{t_{j}=\varepsilon^{j}t},\nonumber\\
\Downarrow\nonumber\\
\frac{dA}{dt}=-\varepsilon\frac{3i}{2}|A|^{2}A-\varepsilon^{2}\frac{15i}
{16}|A|^{4}A.\label{Eq4.51}
\end{gather}
Observe that this equation has a unique solution for a given set of initial
conditions regardless of whether the overdetermined system (\ref{Eq4.48}) has
a solution or not. Thus doing the cross derivative test was only motivated by
intellectual curiosity, we did not have to do it.

In summary, (\ref{Eq4.50}) and (\ref{Eq4.51}), determines a perturbation
expansion for $y(t)$ that is uniform for $t\lesssim\varepsilon^{-3}$.

At this point it is reasonable to ask in which sense we have made progress. We
started with one second order nonlinear ODE for a real function $y(t)$ and
have ended up with one first order nonlinear ODE for a complex function
$A(t)$.

  This question actually has two different answers.
The first one is that it is possible to get an analytical solution for
(\ref{Eq4.51}) whereas this is not possible for the original equation
(\ref{Eq4.35}). This possibility might however easily get lost as we proceed
to higher order in $\varepsilon$, since this will add more terms to the
amplitude equation. But even if we can not solve the amplitude equation
exactly, it is a fact that amplitude equations with the \textit{same
}mathematical structure will arise when we apply the multiple scale method to
many \textit{different} equations. Thus any insight into an amplitude equation
derived by some mathematical analysis has relevance for many different
situations. This is clearly very useful.

There is however a second, more robust, answer to the question of whether we
have made progress or not. From a numerical point of view, there is an
important difference between (\ref{Eq4.35}) and (\ref{Eq4.51}). If we solve
(\ref{Eq4.35}) numerically, the time step is constrained by the oscillation
period of the linearized system
\begin{equation}
\frac{d^{2}y}{dt^{2}}+y=0.\label{Eq4.52}
\end{equation}

\noindent which is of order $T\sim1$, whereas if we solve (\ref{Eq4.51}), numerically the
timestep is constrained by the period $T\sim\varepsilon^{-1}$. Therefore, if we
want to propagate out to a time $t\sim\varepsilon^{-2}$, we must take on the
order of $\varepsilon^{-2}$ time steps if we use (\ref{Eq4.52}) whereas we
only need on the order of $\varepsilon^{-1}$ time steps using (\ref{Eq4.51}).
If $\varepsilon$ is very small the difference in the number of time steps can
be highly significant. From this point of view, the multiple scale method is a
\textit{reformulation} that is the key element in a fast \textit{numerical} method
for solving weakly nonlinear ordinary and partial differential equation.

Let us next turn to the problem of fitting the initial conditions. Using
equation (\ref{Eq4.50}) and the initial conditions from (\ref{Eq4.35}) we get,
truncating at order $\varepsilon$, the following equations
\begin{align}
A(0)-\varepsilon\frac{1}{8}A^{3}(0)+(\ast)  &  =1,\nonumber\\
iA(0)-\varepsilon(\frac{3i}{2}|A(0)|^{2}A(0)+\frac{3i}{8}A^{3}(0))+(\ast)  &
=0.\label{Eq4.53}
\end{align}
The solution for $\varepsilon=0$ is
\begin{equation}
A(0)=\frac{1}{2}.\label{Eq4.54}
\end{equation}
For $\varepsilon>0$ we solve the equation by Newton iteration starting with
the solution for $\varepsilon=0$. This will give us the initial condition for
the amplitude equation correct to this order in $\varepsilon$

\begin{figure}[htbp]
\centering
\includegraphics[
natheight=2.460400in,
natwidth=3.880400in,
height=2.501in,
width=3.9288in
]{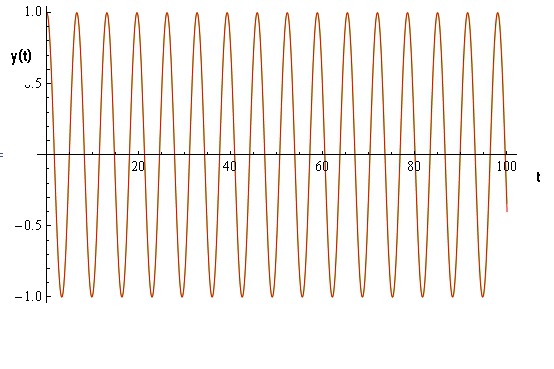}
\caption{Comparing the multiple scale solution, while keeping only the first
term in the amplitude equation(red), to a numerical solution(green) for
$t\lesssim\varepsilon^{-2}$}
\label{fig4}
\end{figure}

\begin{figure}[htbp]
\centering
\includegraphics[
natheight=2.460400in,
natwidth=3.999800in,
height=2.501in,
width=4.0491in
]{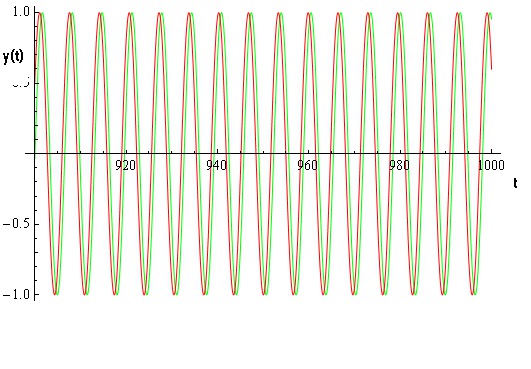}
\caption{Comparing the multiple scale solution, while keeping only the first
term in the amplitude equation(red), to a numerical solution(green) for
$t\lesssim\varepsilon^{-3}$}
\label{fig5}
\end{figure}

\begin{figure}[ptb]
\centering
\includegraphics[
natheight=2.460400in,
natwidth=3.594200in,
height=2.501in,
width=3.6409in
]{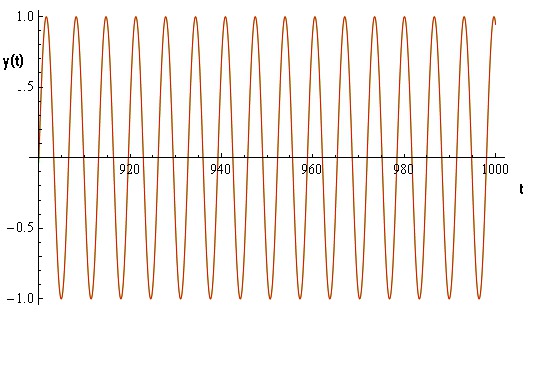}
\caption{Comparing the multiscale solution, while keeping both terms in the
amplitude equation(red), to a numerical solution(green) for
$t\lesssim\varepsilon^{-3}$}
\label{fig6}
\end{figure}

In figure (\ref{fig4}) we compare the multiple scale solution, keeping only
the first term in the amplitude equation, to a high precision numerical
solution for $\varepsilon=0.1$ for
$t\lesssim\varepsilon^{-2}$. We see that the
perturbation solution is very accurate for this range of $t$. In figure
(\ref{fig5}) we do the same comparison as in figure (\ref{fig4}) but now for
$t\lesssim\varepsilon^{-3}$. As expected the multiple scale solution and the
numerical solution starts to deviate for this range of $t$. In figure
(\ref{fig6}) we make the same comparison as in figure (\ref{fig5}), but now
include both terms in the amplitude equation. We see that high accuracy is
restored for the multiple scale solution for $t\lesssim\varepsilon^{-3}$.

\subsection*{Example 2}

Let us consider the weakly nonlinear equation
\begin{equation}
\frac{d^{2}y}{dt^{2}}+\frac{dy}{dt}+\varepsilon y^{2}=0,\text{ \ \ }
t>0.\label{Eq4.55}
\end{equation}
We want to apply the multiple scale method, and introduce a function
$h(t_{0},t_{1},t_{2},..)$ such that
\begin{equation}
y(t)=h(t_{0},t_{1},t_{2},..)|_{t_{j}=\varepsilon^{j}t},\label{Eq4.56}
\end{equation}
is a solution to equation (\ref{Eq4.55}). As usual we have the formal
expansions
\begin{align}
\frac{d}{dt}  &  =\partial_{t_{0}}+\varepsilon\partial_{t_{1}}+\varepsilon
^{2}\partial_{t_{2}}+...\;\;,\label{Eq4.57}\\
h  &  =h_{0}+\varepsilon h_{1}+\varepsilon^{2}h_{2}+...\;\;.\label{Eq4.58}
\end{align}
Inserting (\ref{Eq4.56}),(\ref{Eq4.57}) and (\ref{Eq4.58}) into equation
(\ref{Eq4.55}) and expanding, we get
\begin{gather}
(\partial_{t_{0}}+\varepsilon\partial_{t_{1}}+\varepsilon^{2}\partial_{t_{2}
}+...)(\partial_{t_{0}}+\varepsilon\partial_{t_{1}}+\varepsilon^{2}
\partial_{t_{2}}+...)\nonumber\\
(h_{0}+\varepsilon h_{1}+\varepsilon^{2}h_{2}+...)+(\partial_{t_{0}
}+\varepsilon\partial_{t_{1}}+\varepsilon^{2}\partial_{t_{2}}+...)\nonumber\\
(h_{0}+\varepsilon h_{1}+\varepsilon^{2}h_{2}+...)\nonumber\\
=-\varepsilon(h_{0}+\varepsilon h_{1}+\varepsilon^{2}h_{2}+...)^{2},\nonumber\\
\Downarrow\nonumber\\
(\partial_{t_{0}t_{0}}+\varepsilon(\partial_{t_{0}t_{1}}+\partial_{t_{1}t_{0}
})+\varepsilon^{2}(\partial_{t_{0}t_{2}}+\partial_{t_{1}t_{1}}+\partial
_{t_{2}t_{0}})+...)\nonumber\\
(h_{0}+\varepsilon h_{1}+\varepsilon^{2}h_{2}+...)+(\partial_{t_{0}
}+\varepsilon\partial_{t_{1}}+\varepsilon^{2}\partial_{t_{2}}+...)\nonumber\\
(h_{0}+\varepsilon h_{1}+\varepsilon^{2}h_{2}+...)=-\varepsilon h_{0}
^{2}-\varepsilon^{2}2h_{0}h_{1}+...\;\;,\nonumber\\
\Downarrow\nonumber\\
\partial_{t_{0}t_{0}}h_{0}+\partial_{t_{0}}h_{0}+\varepsilon(\partial
_{t_{0}t_{0}}h_{1}+\partial_{t_{0}}h_{1}+\partial_{t_{0}t_{1}}h_{0}
+\partial_{t_{1}t_{0}}h_{0}+\partial_{t_{1}}h_{0})\nonumber\\
+\varepsilon^{2}(\partial_{t_{0}t_{0}}h_{2}+\partial_{t_{0}}h_{2}
+\partial_{t_{0}t_{1}}h_{1}+\partial_{t_{1}t_{0}}h_{1}+\partial_{t_{0}t_{2}
}h_{0}+\partial_{t_{1}t_{1}}h_{0}\nonumber\\
+\partial_{t_{2}t_{0}}h_{0}+\partial_{t_{1}}h_{1}+\partial_{t_{2}}
h_{0})+...=-\varepsilon h_{0}^{2}-\varepsilon^{2}2h_{0}h_{1}+...\;\;,\label{Eq4.59}
\end{gather}
which gives us the perturbation hierarchy
\begin{align}
\partial_{t_{0}t_{0}}h_{0}+\partial_{t_{0}}h_{0}  &  =0,\nonumber\\
& \nonumber\\
\partial_{t_{0}t_{0}}h_{1}+\partial_{t_{0}}h_{1}  &  =-h_{0}^{2}
-\partial_{t_{0}t_{1}}h_{0}-\partial_{t_{1}t_{0}}h_{0}-\partial_{t_{1}}
h_{0},\nonumber\\
& \nonumber\\
\partial_{t_{0}t_{0}}h_{2}+\partial_{t_{0}}h_{2}  &  =-2h_{0}h_{1}
-\partial_{t_{0}t_{1}}h_{1}-\partial_{t_{1}t_{0}}h_{1}-\partial_{t_{0}t_{2}
}h_{0}\nonumber\\
&  -\partial_{t_{1}t_{1}}h_{0}-\partial_{t_{2}t_{0}}h_{0}-\partial_{t_{1}
}h_{1}-\partial_{t_{2}}h_{0}.\label{Eq4.60}
\end{align}
The general solution to the first equation in the perturbation hierarchy is
\begin{equation}
h_{0}(t_{0},t_{1},t_{2},...)=A_{0}(t_{1},t_{2},..)+B_{0}(t_{1},t_{2}
,...)e^{-t_{0}},\label{Eq4.61}
\end{equation}
where $A_{0}$ and $B_{0}$ are real functions of their arguments. Inserting
$h_{0}$ into the second equation in the hierarchy we get
\begin{equation}
\partial_{t_{0}t_{0}}h_{1}+\partial_{t_{0}}h_{1}=-\partial_{t_{1}}A_{0}
-A_{0}^{2}+(\partial_{t_{1}}B_{0}-2A_{0}B_{0})e^{-t_{0}}-B_{0}^{2}e^{-2t_{0}
}.\label{Eq4.62}
\end{equation}
In order to remove secular terms we must postulate that
\begin{align}
\partial_{t_{1}}A_{0}  &  =-A_{0}^{2},\nonumber\\
\partial_{t_{1}}B_{0}  &  =2A_{0}B_{0}.\label{Eq4.63}
\end{align}
Equation (\ref{Eq4.62}) simplifies into
\begin{equation}
\partial_{t_{0}t_{0}}h_{1}+\partial_{t_{0}}h_{1}=-B_{0}^{2}e^{-2t_{0}
},\label{Eq4.64}
\end{equation}
which has a special solution
\begin{equation}
h_{1}(t_{0},t_{1},..)=-\frac{1}{2}B_{0}^{2}e^{-2t_{0}}.\label{Eq4.65}
\end{equation}
Inserting (\ref{Eq4.61}) and (\ref{Eq4.65}) into the third equation in the
perturbation hierarchy we get
\begin{equation}
\partial_{t_{0}t_{0}}h_{2}+\partial_{t_{0}}h_{2}=-\partial_{t_{2}}
A_{0}-\partial_{t_{1}t_{1}}A_{0}+(\partial_{t_{2}}B_{0}-\partial_{t_{1}t_{1}
}B_{0})e^{-t_{0}}+NST.\label{Eq4.66}
\end{equation}
In order to remove secular terms we must postulate that%
\begin{align}
\partial_{t_{2}}A_{0}  &  =-\partial_{t_{1}t_{1}}A_{0},\nonumber\\
\partial_{t_{2}}B_{0}  &  =\partial_{t_{1}t_{1}}B_{0}.\label{Eq4.67}
\end{align}
We can as usual use (\ref{Eq4.63}) to simplify (\ref{Eq4.67}). We are thus
lead to the following overdetermined system for $A_{0}$ and $B_{0}$.
\begin{align}
\partial_{t_{1}}A_{0}  &  =-A_{0}^{2},\nonumber\\
\partial_{t_{1}}B_{0}  &  =2A_{0}B_{0},\nonumber\\
\partial_{t_{2}}A_{0}  &  =-2A_{0}^{3},\nonumber\\
\partial_{t_{2}}B_{0}  &  =2A_{0}^{2}B_{0}.\label{Eq4.68}
\end{align}
In order to satisfy our academic curiosity, let us do the cross derivative
test for solvability of (\ref{Eq4.68}).
\begin{align}
\partial_{t_{1}t_{2}}A_{0}  &  =-2\partial_{t_{1}}A_{0}^{3}=-6A_{0}
^{2}\partial_{t_{1}}A_{0}=6A_{0}^{4},\nonumber\\
\partial_{t_{2}t_{1}}A_{0}  &  =-\partial_{t_{2}}A_{0}^{2}=-2A_{0}
\partial_{t_{2}}A_{0}=4A_{0}^{4},\nonumber\\
& \nonumber\\
\partial_{t_{1}t_{2}}B_{0}  &  =2\partial_{t_{1}}(A_{0}^{2}B_{0}
)=4A_{0}\partial_{t_{1}}A_{0}B_{0}+2A_{0}^{2}\partial_{t_{1}}B_{0}
=0,\nonumber\\
\partial_{t_{2}t_{1}}B_{0}  &  =2\partial_{t_{2}}(A_{0}B_{0})=2\partial
_{t_{2}}A_{0}B_{0}+2A_{0}\partial_{t_{2}}B_{0}=0.\nonumber
\end{align}
We see that the test fails, so the system (\ref{Eq4.68}) has no solutions.
However the multiple scale method does \textit{not} fail since we are not
actually interested in the functions $A_{0}$ and $B_{0}$ that defines $h_{0}$,
but is rather interested in the function $y(t)$. Define two amplitudes $A(t)$
and $B(t)$ by%
\begin{align}
A(t)  &  =A_{0}(t_{1},t_{2},...)|_{t_{j}=\varepsilon^{j}t},\nonumber\\
B(t)  &  =B_{0}(t_{1},t_{2},...)|_{t_{j}=\varepsilon^{j}t},\label{Eq4.70}
\end{align}
then the solution to (\ref{Eq4.55}) is
\begin{equation}
y(t)=A(t)+B(t)e^{-t}-\varepsilon\frac{1}{2}B^{2}(t)e^{-2t}+O(\varepsilon
^{2}),\label{Eq4.71}
\end{equation}
where the amplitudes $A(t)$ and $B(t)$ satisfy the equations
\begin{align}
\frac{dA}{dt}(t)  &  =\{(\partial_{t_{0}}+\varepsilon\partial_{t_{1}
}+\varepsilon^{2}\partial_{t_{2}}+...)A_{0}(t_{1},t_{2},...)\}|_{t_{j}
=\varepsilon^{j}t},\nonumber\\
&  \Downarrow\nonumber\\
\frac{dA}{dt}(t)  &  =\{-\varepsilon A^{2}(t_{1},t_{2},...)-2\varepsilon
^{2}A(t_{1},t_{2},...)^{3}\}|_{t_{j}=\varepsilon^{j}t},\nonumber\\
&  \Downarrow\nonumber\\
\frac{dA}{dt}  &  =-\varepsilon A^{2}-2\varepsilon^{2}A^{3}.\label{Eq4.72}
\end{align}
and
\begin{align}
\frac{dB}{dt}(t)  &  =\{(\partial_{t_{0}}+\varepsilon\partial_{t_{1}
}+\varepsilon^{2}\partial_{t_{2}}+...)B_{0}(t_{1},t_{2},...)\}|_{t_{j}
=\varepsilon^{j}t},\nonumber\\
&  \Downarrow\nonumber\\
\frac{dB}{dt}(t)  &  =\{(2\varepsilon A_{0}(t_{1},t_{2},...)B_{0}(t_{1}
,t_{2},...)+2\varepsilon^{2}A_{0}^{2}(t_{1},t_{2},...)B_{0}(t_{1}
,t_{2},......)\}|_{t_{j}=\varepsilon^{j}t},\nonumber\\
&  \Downarrow\nonumber\\
\frac{dB}{dt}  &  =2\varepsilon AB+2\varepsilon^{2}A^{2}B.\label{Eq4.73}
\end{align}
Given the initial conditions for $A$ and $B$, equations (\ref{Eq4.72}) and
(\ref{Eq4.73}) clearly has a unique solution and our multiple scale method
will ensure that the perturbation expansion (\ref{Eq4.71}) will stay uniform
for $t\lesssim\varepsilon^{-3}$. As for the previous example, the initial
conditions $A(0)$ and $B(0)$ are calculated from the initial conditions for
(\ref{Eq4.55}) by a Newton iteration. Thus we see again that the existence or
not of $h(t_{0},..)$ is irrelevant for constructing a uniform perturbation expansion.

The system (\ref{Eq4.72}) and (\ref{Eq4.73}) can be solved analytically in
terms of implicit functions. However, as we have discussed before,
\ analytical solvability is nice, but not robust. If we take the expansion to
order $\varepsilon^{3}$, more terms are added to the amplitude equations and
the property of analytic solvability can easily be lost. What \textit{is}
robust is that the presense of $\varepsilon$ in the amplitude equations makes
(\ref{Eq4.72}) and (\ref{Eq4.73}) together with (\ref{Eq4.71}) into a fast
numerical scheme for solving the ordinary differential equation (\ref{Eq4.55}%
). This property does \textit{not} go away if we take the perturbation
expansion to higher order in $\varepsilon$.

\subsection*{Example 3}

So far, we have only been applying the method of multiple scales to scalar
ODEs. This is not a limitation on the method, it may equally well be applied
to systems of ordinary differential equations. The mechanics of the method for
systems of equations is very similar to what we have seen for scalar
equations. The only major difference is how we decide which terms are secular
and must be removed. For systems, this problem is solved by using the Fredholm
Alternative theorem, this is in fact one of the major areas of application for
this theorem in applied mathematics.

Let us consider the following system of two coupled second order ODEs.
\begin{align}
x^{\prime\prime}+2x-y  & =\varepsilon xy^{2},\nonumber\\
y^{\prime\prime}+3y-2x  & =\varepsilon yx^{2},\label{Eq4.74}
\end{align}
where $\varepsilon\ll1$. We will solve the system using the method of multiple
scales and introduce therefore two functions $h=h(t_{0},t_{1},...)$ and
$k=k(t_{0},t_{1},...)$ such that
\begin{align}
x(t)  & =h(t_{0},t_{1},...)|_{t_{j}=\varepsilon^{j}t},\nonumber\\
y(t)  & =k(t_{0},t_{1},...)|_{t_{j}=\varepsilon^{j}t},\label{Eq4.75}
\end{align}
is a solution to equation (\ref{Eq4.74}). As usual we have
\begin{equation}
\frac{d}{dt}=\partial_{t_{0}}+\varepsilon\partial_{t_{1}}+...\;\;,\label{Eq4.76}
\end{equation}
and for $h$ and $k$ we introduce the expansions
\begin{align}
h  & =h_{0}+\varepsilon h_{1}+...\,\,,\nonumber\\
k  & =k_{0}+\varepsilon k_{1}+...\;\;.\label{Eq4.77}
\end{align}
Inserting (\ref{Eq4.75}),(\ref{Eq4.76}) and (\ref{Eq4.77}) into equation
(\ref{Eq4.74}), and expanding everything in sight to first order in
$\varepsilon$ we get, after some tedious algebra, the following perturbation
hierarchy
\begin{align}
\partial_{t_{0}t_{0}}h_{0}+2h_{0}-k_{0}  & =0,\nonumber\\
\partial_{t_{0}t_{0}}k_{0}+3k_{0}-2h_{0}  & =0,\label{Eq4.78}\\
& \nonumber\\
\partial_{t_{0}t_{0}}h_{1}+2h_{1}-k_{1}  & =-\partial_{t_{0}t_{1}}
h_{0}-\partial_{t_{1}t_{0}}h_{0}+h_{0}k_{0}^{2},\nonumber\\
\partial_{t_{0}t_{0}}k_{1}+3k_{1}-2h_{1}  & =-\partial_{t_{0}t_{1}}
k_{0}-\partial_{t_{1}t_{0}}k_{0}+k_{0}h_{0}^{2}.\label{Eq4.79}
\end{align}
Let us start by finding the general solution to the order $\varepsilon^{0} $\newline
equations (\ref{Eq4.78}). They can be written as the following linear system
\begin{equation}
\partial_{t_{0}t_{0}}\left(
\begin{array}
[c]{c}
h_{0}\\
k_{0}
\end{array}
\right)  =\left(
\begin{array}
[c]{cc}
-2 & 1\\
2 & -3
\end{array}
\right)  \left(
\begin{array}
[c]{c}
h_{0}\\
k_{0}
\end{array}
\right),\label{Eq4.80}
\end{equation}
Let us look for a solution of the form
\begin{equation}
\left(
\begin{array}
[c]{c}
h_{0}\\
k_{0}
\end{array}
\right)  =\mathbf{\alpha}e^{i\omega t_{0}},\label{Eq4.81}
\end{equation}
where $\mathbf{\alpha}$ is a unknown vector and $\omega$ an unknown real
number. Inserting (\ref{Eq4.81}) into the system (\ref{Eq4.80}) and cancelling
a common factor we get the the following linear algebraic equation
\begin{equation}
\left(
\begin{array}
[c]{cc}
-2+\omega^{2} & 1\\
2 & -3+\omega^{2}
\end{array}
\right)  \mathbf{\alpha}=0.\label{Eq4.82}
\end{equation}
For there to be a nontrivial solution, the determinant of the matrix has to be
zero. This condition leads to the following polynomial equation for $\omega$
\begin{equation}
\omega^{4}-5\omega^{2}+4=0,\label{Eq4.83}
\end{equation}
which has four real solutions
\[
\omega_{1}=1,\omega_{2}=-1,\omega_{3}=2,\omega_{4}=-2.
\]
A basis for the solution space of (\ref{Eq4.82}) corresponding to
$\omega=\omega_{1},\omega_{2}$ is
\begin{equation}
\mathbf{\alpha}=\left(
\begin{array}
[c]{c}
1\\
1
\end{array}
\right), \label{Eq4.84}
\end{equation}
and a basis corresponding to $\omega=\omega_{3},\omega_{4}$ is
\begin{equation}
\mathbf{\beta}=\left(
\begin{array}
[c]{c}
1\\
-2
\end{array}
\right). \label{Eq4.85}
\end{equation}
It is then clear that a basis for the solution space for the linear system
(\ref{Eq4.80}) is
\[
\mathbf{\alpha e}^{\pm it_{0}},\mathbf{\beta}e^{\pm2it_{0}}.
\]
Therefore a general complex solution to (\ref{Eq4.80}) is
\begin{equation}
\left(
\begin{array}
[c]{c}
h_{0}\\
k_{0}
\end{array}
\right)  =A_{1}\mathbf{\alpha}e^{it_{0}}+A_{2}\mathbf{\alpha}e^{-it_{0}}%
+B_{1}\mathbf{\beta}e^{2it_{0}}+B_{2}\mathbf{\beta}e^{-2it_{0}}.\label{Eq4.86}
\end{equation}
However, we are looking for real solutions to the original system
(\ref{Eq4.74}), and in order to ensure reality for (\ref{Eq4.86}) we must
choose%
\begin{align*}
A_{1}  & =A_{0}^{\ast},\text{ \ \ }A_{2}=A_{0},\text{ \ \ \ }A_{0}=A_{0}
(t_{1},t_{2},...),\\
B_{1}  & =B_{0}^{\ast},\text{ \ \ }B_{2}=B_{0},\text{ \ \ \ }B_{0}=B_{0}
(t_{1},t_{2},...).
\end{align*}
Thus, a general solution to (\ref{Eq4.80}) is
\begin{equation}
\left(
\begin{array}
[c]{c}
h_{0}\\
k_{0}
\end{array}
\right)  =A_{0}\mathbf{\alpha}e^{-it_{0}}+B_{0}\mathbf{\beta}e^{-2it_{0}
}+(\ast).\label{Eq4.87}
\end{equation}
In component form, the general real solution is
\begin{align}
h_{0}  & =A_{0}e^{-it_{0}}+B_{0}e^{-2it_{0}}+(\ast),\nonumber\\
k_{0}  & =A_{0}e^{it_{0}}-2B_{0}e^{-2it_{0}}+(\ast).\label{Eq4.88}
\end{align}
We now insert the expressions (\ref{Eq4.88}) into the order $\varepsilon$
equations (\ref{Eq4.79}). After a large amount of tedious algebra, Mathematica
can be useful here, we find that the order $\varepsilon$ equations can be
written in the form
\begin{gather}
\partial_{t_{0}t_{0}}\left(
\begin{array}
[c]{c}
h_{1}\\
k_{1}
\end{array}
\right)  +\left(
\begin{array}
[c]{cc}
2 & -1\\
-2 & 3
\end{array}
\right)  \left(
\begin{array}
[c]{c}
h_{1}\\
k_{1}
\end{array}
\right)  =\left(
\begin{array}
[c]{c}
4B_{0}^{3}\\
-2B_{0}^{3}
\end{array}
\right)  e^{-6it_{0}}+\left(
\begin{array}
[c]{c}
0\\
-3A_{0}B_{0}^{2}
\end{array}
\right)  e^{-5it_{0}}\nonumber\\
-\left(
\begin{array}
[c]{c}
3A_{0}^{\ast}A_{0}^{\ast}B_{0}\\
0
\end{array}
\right)  e^{-4it_{0}}+\left(
\begin{array}
[c]{c}
A_{0}^{2}\\
A_{0}^{3}-3A_{0}^{\ast}B_{0}^{2}
\end{array}
\right)  e^{-3it_{0}}\nonumber\\
+\left(
\begin{array}
[c]{c}
4i\partial_{t_{1}}B_{0}+12|B_{0}|^{2}B_{0}-6|A_{0}|^{2}B_{0}\\
-8i\partial_{t_{1}}B_{0}-6|B_{0}|^{2}B_{0}
\end{array}
\right)  e^{-2it_{0}}\nonumber\\
+\left(
\begin{array}
[c]{c}
2i\partial_{t_{1}}A_{0}+3|A_{0}|^{2}A_{0}\\
2i\partial_{t_{1}}A_{0}+3|A_{0}|^{2}A_{0}-6|B_{0}|^{2}A_{0}
\end{array}
\right)  e^{-it_{0}}+\left(
\begin{array}
[c]{c}
-\frac{3}{2}A_{0}^{\ast}A_{0}^{\ast}B_{0}\\
0
\end{array}
\right)  +(\ast).\label{Eq4.89}
\end{gather}
We are not going to go beyond order $\varepsilon$ so we don't actually need to
solve this equation. What we need to do, however, is to ensure that the
solution is bounded in $t_{0}$. We only need a special solution to
(\ref{Eq4.89}), and because it is a linear equation, such a special solution
can be constructed as a sum of solutions where each solution in the sum
corresponds to a single term from the righhand side of (\ref{Eq4.89}). What we
mean by this is that
\begin{equation}
\left(
\begin{array}
[c]{c}
h_{1}\\
k_{1}
\end{array}
\right)  =\sum_{n=1}^{7}\left(
\begin{array}
[c]{c}
u_{n}\\
v_{n}
\end{array}
\right)  +(\ast),\label{Eq4.90}
\end{equation}
where for example
\begin{align}
\partial_{t_{0}t_{0}}\left(
\begin{array}
[c]{c}
u_{1}\\
v_{1}
\end{array}
\right)  +\left(
\begin{array}
[c]{cc}
2 & -1\\
-2 & 3
\end{array}
\right)  \left(
\begin{array}
[c]{c}
u_{1}\\
v_{1}
\end{array}
\right)   & =\left(
\begin{array}
[c]{c}
4B_{0}^{3}\\
-2B_{0}^{3}
\end{array}
\right)  e^{-6it_{0}},\nonumber\\
\partial_{t_{0}t_{0}}\left(
\begin{array}
[c]{c}
u_{2}\\
v_{2}
\end{array}
\right)  +\left(
\begin{array}
[c]{cc}
2 & -1\\
-2 & 3
\end{array}
\right)  \left(
\begin{array}
[c]{c}
u_{2}\\
v_{2}
\end{array}
\right)   & =\left(
\begin{array}
[c]{c}
0\\
-3A_{0}B_{0}^{2}
\end{array}
\right)  e^{-5it_{0}},\label{Eq4.91}
\end{align}
and so on. For the first equation we look for a solution of the form
\begin{equation}
\left(
\begin{array}
[c]{c}
u_{1}(t_{0})\\
v_{1}(t_{0})
\end{array}
\right)  =\mathbf{\xi}e^{-6it_{0}},\label{Eq4.92}
\end{equation}
where $\mathbf{\xi}$ is a constant vector. Observe that any solution of the
form (\ref{Eq4.92}), is bounded in $t_{0}$. If we insert (\ref{Eq4.92}) into
the first equation in (\ref{Eq4.91}) and cancel the common exponential factor
we find that the unknown vector $\xi$ has to be a solution of the following
linear algebraic system
\[
\left(
\begin{array}
[c]{cc}
-34 & -1\\
-2 & -33
\end{array}
\right)  \mathbf{\xi}=\left(
\begin{array}
[c]{c}
4B_{0}^{3}\\
-2B_{0}^{3}
\end{array}
\right).
\]
The matrix of this system is clearly nonsingular and the solution is
\[
\mathbf{\xi}=\frac{1}{560}\left(
\begin{array}
[c]{c}
-67B_{0}^{3}\\
38B_{0}^{3}
\end{array}
\right),
\]
which gives us the following bounded solution
\[
\left(
\begin{array}
[c]{c}
u_{1}(t_{0})\\
v_{1}(t_{0})
\end{array}
\right)  =\frac{1}{560}\left(
\begin{array}
[c]{c}
-67B_{0}^{3}\\
38B_{0}^{3}
\end{array}
\right)  e^{-6it_{0}}.
\]
A similar approach works for all but the fifth and the sixth term on the
righthand side of equation (\ref{Eq4.89}). For these two terms we run into
trouble. For the fifth term we must solve the equation
\begin{eqnarray}
\partial_{t_{0}t_{0}}\left(
\begin{array}
[c]{c}
u_{5}\\
v_{5}
\end{array}
\right)  +\left(
\begin{array}
[c]{cc}
2 & -1\\
-2 & 3
\end{array}
\right)  \left(
\begin{array}
[c]{c}
u_{5}\\
v_{5}
\end{array}
\right)=\nonumber\\
\left(
\begin{array}
[c]{c}
4i\partial_{t_{1}}B_{0}+12|B_{0}|^{2}B_{0}-6|A_{0}|^{2}B_{0}\\
-8i\partial_{t_{1}}B_{0}-6|B_{0}|^{2}B_{0}
\end{array}
\right)  e^{-2it_{0}}.\label{Eq4.93}
\end{eqnarray}
A bounded trial solution of the form
\begin{equation}
\left(
\begin{array}
[c]{c}
u_{5}(t_{0})\\
v_{5}(t_{0})
\end{array}
\right)  =\mathbf{\xi}e^{-2it_{0}},\label{Eq4.94}
\end{equation}
leads to the following algebraic equation for $\mathbf{\xi}$
\begin{equation}
\left(
\begin{array}
[c]{cc}
-2 & -1\\
-2 & -1
\end{array}
\right)  \mathbf{\xi}=\left(
\begin{array}
[c]{c}
4i\partial_{t_{1}}B_{0}+12|B_{0}|^{2}B_{0}-6|A_{0}|^{2}B_{0}\\
-8i\partial_{t_{1}}B_{0}-6|B_{0}|^{2}B_{0}
\end{array}
\right). \label{Eq4.95}
\end{equation}
The matrix for this linear system is singular, and the system will in general
have no solution. It will only have a solution, which will lead to a bounded
solution for (\ref{Eq4.93}), if the righthand side of (\ref{Eq4.95}) satisfy a
certain constraint. This constraint we get from the Fredholm Alternative
Theorem. Recall that this theorem say that a linear system
\[
M\mathbf{x}=b\mathbf{_{0}},
\]
has a solution only if
\[
(\mathbf{f},\mathbf{b_{0}})=0,
\]
for all vectors $\mathbf{f}$ such that
\[
M^{\ast}\mathbf{f}=0,
\]
where $M^{\ast}$ is the adjoint of the matrix $M$. For a real matric, like the
one we have, $M^{\ast}$ is just the transpose of $M$. For the matric of the
system (\ref{Eq4.95}) we get
\[
\left(
\begin{array}
[c]{cc}
-2 & -2\\
-1 & -1
\end{array}
\right)  \mathbf{f}=0.
\]
A basis for the solution space of this homogenous system can be taken to be
\[
\mathbf{f}=(1,-1).
\]
Thus in order to ensure solvability of the system (\ref{Eq4.95}) we must have
\begin{gather}
(1,-1)\cdot\left(
\begin{array}
[c]{c}
4i\partial_{t_{1}}B_{0}+12|B_{0}|^{2}B_{0}-6|A_{0}|^{2}B_{0}\\
-8i\partial_{t_{1}}B_{0}-6|B_{0}|^{2}B_{0}
\end{array}
\right)  =0,\nonumber\\
\Updownarrow\nonumber\\
4i\partial_{t_{1}}B_{0}+12|B_{0}|^{2}B_{0}-6|A_{0}|^{2}B_{0}+8i\partial
_{t_{1}}B_{0}+6|B_{0}|^{2}B_{0}=0,\nonumber\\
\Updownarrow\nonumber\\
\partial_{t_{1}}B_{0}=\frac{i}{2}(3|B_{0}|^{2}-|A_{0}|^{2})B_{0}.\label{Eq4.96}
\end{gather}
If this condition on the amplitudes is imposed on the original system
(\ref{Eq4.93}), it has a bounded solution. The sixth term in the sum
(\ref{Eq4.90}) must be treated in the same way. The equation that we must
solve is
\begin{eqnarray}
\partial_{t_{0}t_{0}}\left(
\begin{array}
[c]{c}
u_{6}\\
v_{6}
\end{array}
\right)  +\left(
\begin{array}
[c]{cc}
2 & -1\\
-2 & 3
\end{array}
\right)  \left(
\begin{array}
[c]{c}
u_{6}\\
v_{6}
\end{array}
\right)  =\nonumber\\
\left(
\begin{array}
[c]{c}
2i\partial_{t_{1}}A_{0}+3|A_{0}|^{2}A_{0}\\
2i\partial_{t_{1}}A_{0}+3|A_{0}|^{2}A_{0}-6|B_{0}|^{2}A_{0}
\end{array}
\right)  e^{-it_{0}}.
\end{eqnarray}
Using a bounded trial solution of the form
\begin{equation}
\left(
\begin{array}
[c]{c}
u_{6}(t_{0})\\
v_{6}(t_{0})
\end{array}
\right)  =\mathbf{\xi}e^{-it_{0}},\label{Eq4.99}
\end{equation}
leads to the following singular linear system
\begin{equation}
\left(
\begin{array}
[c]{cc}
1 & -1\\
-2 & 2
\end{array}
\right)  \mathbf{\xi}=\left(
\begin{array}
[c]{c}
2i\partial_{t_{1}}A_{0}+3|A_{0}|^{2}A_{0}\\
2i\partial_{t_{1}}A_{0}+3|A_{0}|^{2}A_{0}-6|B_{0}|^{2}A_{0}
\end{array}
\right). \label{Eq4.101}
\end{equation}
For this case we find that the null space of the transpose of the matrix in
(\ref{Eq4.101}) is spanned by the vector
\[
\mathbf{f}=(2,1),
\]
and the Fredholm Alternative gives us the solvability condition
\begin{gather}
(2.1)\cdot\left(
\begin{array}
[c]{c}
2i\partial_{t_{1}}A_{0}+3|A_{0}|^{2}A_{0}\\
2i\partial_{t_{1}}A_{0}+3|A_{0}|^{2}A_{0}-6|B_{0}|^{2}A_{0}
\end{array}
\right)  =0,\nonumber\\
\Updownarrow\nonumber\\
4i\partial_{t_{1}}A_{0}+6|A_{0}|^{2}A_{0}+2i\partial_{t_{1}}A_{0}
+3|A_{0}|^{2}A_{0}-6|B_{0}|^{2}A_{0}=0,\nonumber\\
\Updownarrow\nonumber\\
\partial_{t_{1}}A_{0}=\frac{i}{2}(3|A_{0}|^{2}-2|B_{0}|^{2})A_{0}.
\label{Eq4.102}
\end{gather}
The solutions $h_{1}$ and $k_{1}$ are thus bounded if we impose the following
two conditions on the amplitudes $A_{0}$ and $B_{0}$
\begin{align}
\partial_{t_{1}}A_{0}  & =\frac{i}{2}(3|A_{0}|^{2}-2|B_{0}|^{2})A_{0},\nonumber\\
\partial_{t_{1}}B_{0}  & =\frac{i}{2}(3|B_{0}|^{2}-|A_{0}|^{2})B_{0}.\label{Eq4.103}
\end{align}
Returning to the original variables $x(t)$ and $y(t)$ in the usual way, we
have thus found the following approximate solution to our system
(\ref{Eq4.74})
\begin{align}
x(t)  & =A(t)e^{-it}+B(t)e^{-2it}+\mathcal{O}(\varepsilon),\nonumber\\
y(t)  & =A(t)e^{-it}-2B(t)e^{-2it}+\mathcal{O}(\varepsilon),\label{Eq4.104}
\end{align}
where the amplitudes $A(t)$ and $B(t)$ are defined by
\begin{align*}
A(t)  & =A_{0}(t_{1},t_{2},...)|_{t_{j}=\varepsilon^{j}t},\\
B(t)  & =B_{0}(t_{1},t_{2},...)|_{t_{j}=\varepsilon^{j}t},
\end{align*}
and satisfy the following equations
\begin{align}
\partial_{t}A  & =\varepsilon\frac{i}{2}(3|A|^{2}-2|B|^{2})A,\label{Eq4.105}\\
\partial_{t}B  & =\varepsilon\frac{i}{2}(3|B|^{2}-|A|^{2})B.\nonumber
\end{align}
The expansions (\ref{Eq4.104}) are uniform for $t\lesssim\varepsilon^{-2}$.

The amplitude equations (\ref{Eq4.105}) looks complicated, but they are
special in the sense that they can be solved exactly. We have noted before
that the amplitude equations that appears when we use the method of multiple
scale tends to be special. We will see more of this later when we apply the
method to partial differential equations.

Observe that
\begin{align*}
\partial_{t}|A|^{2}  & =\partial_{t}(AA^{\ast})=A^{\ast}\partial
_{t}A+A\partial_{t}A^{\ast}\\
& =A^{\ast}(\frac{i}{2}(3|A|^{2}-2|B|^{2})A)+A(-\frac{i}{2}(3|A|^{2}
-2|B|^{2})A^{\ast})\\
& =\frac{i}{2}(3|A|^{4}-2|B|^{2}|A|^{2}-3|a|^{4}+2|B|^{2}|A|^{2})=0.
\end{align*}
Thus $|A(t)|=|A(0)|$ for all $t$. In a similar way we find that
$|B(t)|=|B(0)|$. Therefore the amplitude equations can be written as
\begin{align*}
\partial_{t}A  & =\frac{i}{2}(3|A(0)|^{2}-2|B(0)|^{2})A,\\
\partial_{t}B  & =\frac{i}{2}(3|B(0)|^{2}-|A(0)|^{2})B.
\end{align*}
and this system is trivial to solve. We find
\begin{align}
A(t)  & =A(0)e^{\frac{i}{2}(3|A(0)|^{2}-2|B(0)|^{2})t},\nonumber\\
B(t)  & =B(0)e^{\frac{i}{2}(3|B(0)|^{2}-|A(0)|^{2})t}.\label{Eq4.106}
\end{align}
The formulas (\ref{Eq4.106}) together with the expansions (\ref{Eq4.104})
gives us an approximate analytic solution to the original system (\ref{Eq4.74}).

\section{Boundary layer problems for ODEs}

Boundary layer problems first appeared in the theory of fluids. However, boundary layer problems are in no
way limited to fluid theory, but occurs in all areas of science and technology.

In these lecture notes, we will not worry about the physical context for these problems, but
will focus on how to apply the multiple scale method to solve a given problem
of this type. As usual we learn by doing examples.

\subsection*{Example 1}

Let us consider the following linear boundary value problem
\begin{align}
\varepsilon y^{\prime\prime}(x)+y^{\prime}(x)-y(x)  & =0,\text{ \ \ }
0<x<1,\nonumber\\
y(0)  & =1,\nonumber\\
y(1)  & =0.\label{Eq6.1}
\end{align}
We will assume that $\varepsilon\ll1$, and try to solve this problem using a
perturbation methods. The unperturbed problem is clearly
\begin{align}
y^{\prime}(x)-y(x)  & =0,\text{ \ \ }0<x<1,\nonumber\\
y(0)  & =1,\nonumber\\
y(1)  & =0.\label{Eq6.2}
\end{align}
The general solution to the differential equation is
\[
y(x)=Ae^{x},
\]
and fitting the boundary condition at $x=0$ we find that
\[
y(x)=e^{x},
\]
but for this solution we have
\[
y(1)=e\neq0,
\]
so the unperturbed problem has no solution. Our perturbation approach fail at
the very first step; there is no unperturbed solution that we can start
calculating corrections to! What is going on?

What is going on is that equation (\ref{Eq6.1}) is a singular perturbation
problem. For $\varepsilon\neq0$, we have a second order ODE, whose general
solution has two free constants that can be fitted to the two boundary
conditions, whereas for $\varepsilon=0$ we have a first order ODE whose
general solution has only one free constant. This single constant can in
general not be fitted to two boundary conditions.

We have seen such singular perturbation problems before when we applied
perturbation methods to polynomial equations. For the polynomial case, the
unperturbed problem was of lower algebraic order than the perturbed problem.
Here the unperturbed problem is of lower differential order than the perturbed problem.

For the polynomial case we solved the singular perturbation problem by
transforming it into a regular perturbation problem using a change of
variables. We do the same here.

Let
\begin{equation}
x=\varepsilon^{p}\xi,\text{ \ \ \ }y(x)=u(\frac{x}{\epsilon^{p}}),\label{Eq6.3}
\end{equation}
then the function $u(\xi)$ satisfy the equation
\begin{equation}
u^{\prime\prime}(\xi)+\varepsilon^{p-1}u^{\prime}(\xi)-\varepsilon^{2p-1}
u(\xi)=0.\label{Eq6.31}
\end{equation}
This equation constitute a regular perturbation problem if we, for example,
choose $p=1$. We thus have the following regularly perturbed boundary value
problem
\begin{align}
u^{\prime\prime}(\xi)+u^{\prime}(\xi)-\varepsilon u(\xi)  & =0,\text{
\ \ }0<\xi<\frac{1}{\varepsilon},\nonumber\\
u(0)  & =1,\nonumber\\
u(\frac{1}{\varepsilon})  & =0.\label{Eq6.4}
\end{align}
Let us try to solve this problem using a perturbation expansion
\begin{equation}
u(\xi)=u_{0}(\xi)+\varepsilon u_{1}(\xi)+...\;\;.\label{Eq6.5}
\end{equation}
We will solve the problem by first finding $u_{0}$ and $u_{1}$ and then
fitting the boundary conditions. If we insert the perturbation expansion
(\ref{Eq6.5}) into the equation (\ref{Eq6.4}) we find the following
perturbation hierarchy to first order in $\varepsilon$
\begin{align}
u_{0}^{\prime\prime}+u_{0}^{\prime}  & =0,\label{Eq6.51}\\
u_{1}^{\prime\prime}+u_{1}^{\prime}  & =u_{0}.\label{Eq6.52}
\end{align}
The general solution to the first equation in the perturbation hierarchy (\ref{Eq6.51}), is
clearly
\begin{equation}
u_{0}(\xi)=A_{0}+B_{0}e^{-\xi}.\label{Eq6.6}
\end{equation}
If we insert the solution (\ref{Eq6.6}) into the second equation in the
perturbation hierarchy (\ref{Eq6.52}), we get
\begin{equation}
u_{1}^{\prime\prime}+u_{1}^{\prime}=A_{0}+B_{0}e^{-\xi}.\label{Eq6.7}
\end{equation}
Note, that we only need a particular solution to this equation, since the first
term in the perturbation expansion (\ref{Eq6.5}) already have two free
constants, and we only need two constants to fit the two boundary data.
Integrating equation (\ref{Eq6.7}) once we get
\[
u_{1}^{\prime}+u_{1}=A_{0}\xi-B_{0}e^{-\xi},
\]
and using an integrating factor we get the following particular solution
\[
u_{1}(\xi)=A_{0}(\xi-1)-B_{0}\xi e^{-\xi}.
\]
Thus our perturbation solution to first order in $\varepsilon$ is
\begin{equation}
u(\xi)=A_{0}+B_{0}e^{-\xi}+\varepsilon\left(  A_{0}(\xi-1)-B_{0}\xi e^{-\xi
}\right)  +...\;\;.
\end{equation}
The two constants are fitted to the boundary conditions using the following
two equations
\begin{align*}
u(0)  & =1\text{ \ \ }\Longleftrightarrow\text{ \ }A_{0}+B_{0}-\varepsilon
A_{0}=1,\text{\ }\\
u(\frac{1}{\varepsilon})  & =0\text{ \ \ }\Longleftrightarrow\text{ }
A_{0}+B_{0}e^{-\frac{1}{\varepsilon}}+\varepsilon\left(  A_{0}(\frac
{1}{\varepsilon}-1)-B_{0}\frac{1}{\varepsilon}e^{-\frac{1}{\varepsilon}
}\right)  =0.
\end{align*}
However at this point disaster strikes. When we evaluate the solution at the
right boundary $\xi=\frac{1}{\varepsilon}$, using the perturbation expansion,
the ordering of terms is violated. The first and the second term in the
expansion are of the same order. This can not be allowed. Our perturbation
method fails. The reason why the direct perturbation expansion (\ref{Eq6.5})
fails is similar to the reason why the direct perturbation expansion failed
for the weakly damped oscillator. In both cases the expansions failed because
they became nonuniform when we evaluated the respective functions at values of
the independent variable that was of order $\frac{1}{\varepsilon}$.

We will resolve the problem with the direct expansion (\ref{Eq6.5}) by using
the method of multiple scales to derive a perturbation expansion for the
solution to the equation (\ref{Eq6.4}) that is uniform for $\xi\lesssim
\frac{1}{\varepsilon}$ and then use this expansion to satisfy the boundary
conditions at $x=0$ and $x=\frac{1}{\varepsilon}$.

\begin{figure}[h]
\centering
\includegraphics[
natheight=3.270700in,
natwidth=4.968400in,
height=2.2987in,
width=3.4817in
]
{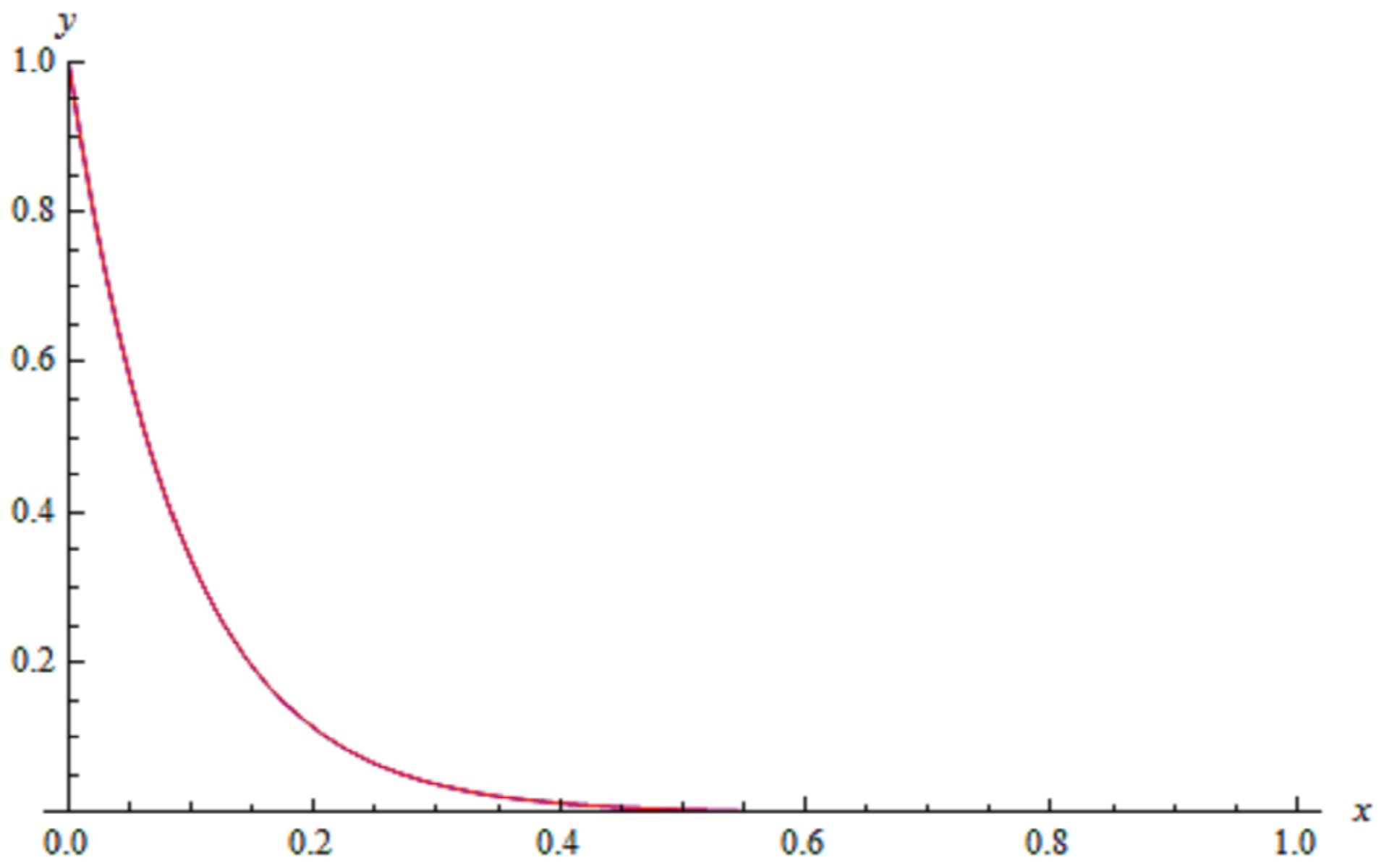}
\caption{Comparing the exact and approximate solution to the singularly
perturbed linear boundary value problem in example 1}
\label{fig7}
\end{figure}

We thus introduce a function $h=h(\xi_{0},\xi_{1},...)$, where $h$ is a function that will be designed 
to ensure that the function $u$,  defined by
\begin{equation}
u(\xi)=h(\xi_{0},\xi_{1},...)|_{\xi_{j}=\varepsilon^{j}\xi},\label{Eq6.8}
\end{equation}
is a solution to the equation (\ref{Eq6.4}). For the differential operator we
have in the usual way an expansion
\begin{equation}
\frac{d}{d\xi}=\partial_{\xi_{0}}+\varepsilon\partial_{\xi_{1}}+\varepsilon
^{2}\partial_{\xi_{2}}+...\;\;,\label{Eq6.9}
\end{equation}
and for the function $h$ we introduce the expansion
\begin{equation}
h=h_{0}+\varepsilon h_{1}+\varepsilon^{2}h_{2}+...\;\;.\label{Eq6.10}
\end{equation}
Inserting (\ref{Eq6.8}),(\ref{Eq6.9}) and (\ref{Eq6.10}) into the equation
(\ref{Eq6.4}) and expanding everything in sigh to second order in
$\varepsilon$, we get after a small amount of algebra the following
perturbation hierarchy
\begin{align}
\partial_{\xi_{0}\xi_{0}}h_{0}+\partial_{\xi_{0}}h_{0}  & =0,\nonumber\\
& \nonumber\\
\partial_{\xi_{0}\xi_{0}}h_{1}+\partial_{\xi_{0}}h_{1}  & =h_{0}-\partial
_{\xi_{0}\xi_{1}}h_{0}-\partial_{\xi_{1}\xi_{0}}h_{0}-\partial_{\xi_{1}}
h_{0},\nonumber\\
& \nonumber\\
\partial_{\xi_{0}\xi_{0}}h_{2}+\partial_{\xi_{0}}h_{2}  & =h_{1}-\partial
_{\xi_{0}\xi_{1}}h_{1}-\partial_{\xi_{1}\xi_{0}}h_{1}\nonumber\\
& -\partial_{\xi_{0}\xi_{2}}h_{0}-\partial_{\xi_{1}\xi_{1}}h_{0}-\partial
_{\xi_{2}\xi_{0}}h_{0}\nonumber\\
& -\partial_{\xi_{1}}h_{1}-\partial_{\xi_{2}}h_{0}.\label{Eq6.11}
\end{align}
The general solution to the first equation in the perturbation hierarchy
(\ref{Eq6.11}) is
\begin{equation}
h_{0}(\xi_{0},\xi_{1},\xi_{2},...)=A_{0}(\xi_{1},\xi_{2},...)+B_{0}(\xi
_{1},\xi_{2},...)e^{-\xi_{0}}.\label{Eq6.12}
\end{equation}
We now insert this solution into the righthand side of the second equation in
the perturbation hierarchy. Thus the order $\varepsilon$ equation is of the
form
\[
\partial_{\xi_{0}\xi_{0}}h_{1}+\partial_{\xi_{0}}h_{1}=A_{0}-\partial_{\xi
_{1}}A_{0}+(\partial_{\xi_{1}}B_{0}+B_{0})e^{-\xi_{0}}.
\]
Both terms on the righthand side of the equation are secular and in order to
avoid nonuniformity in our expansion we must enforce the conditions
\begin{align}
\partial_{\xi_{1}}A_{0}  & =A_{0},\nonumber\\
\partial_{\xi_{1}}B_{0}  & =-B_{0}.\label{Eq6.13}
\end{align}
With these conditions in place, the equation for $h_{1}$ simplify into
\[
\partial_{\xi_{0}\xi_{0}}h_{1}+\partial_{\xi_{0}}h_{1}=0.
\]
and for this equation we choose the special solution
\begin{equation}
h_{1}=0.\label{Eq6.14}
\end{equation}
Inserting (\ref{Eq6.12}) and (\ref{Eq6.14}) into the third equation in the
perturbation hierarchy (\ref{Eq6.11}) we get
\[
\partial_{\xi_{0}\xi_{0}}h_{2}+\partial_{\xi_{0}}h_{2}=-\partial_{\xi_{2}
}A_{0}-\partial_{\xi_{1}\xi_{1}}A_{0}+(\partial_{\xi_{2}}B_{0}-\partial
_{\xi_{1}\xi_{1}}B_{0})e^{-\xi_{0}}.
\]
In order to avoid secular terms we enforce the conditions
\begin{align}
\partial_{\xi_{2}}A_{0}  & =-\partial_{\xi_{1}\xi_{1}}A_{0},\nonumber\\
\partial_{\xi_{2}}B_{0}  & =\partial_{\xi_{1}\xi_{1}}B_{0},\label{Eq6.15}
\end{align}
and with this choise the equation for $h_{2}$ simplify into
\[
\partial_{\xi_{0}\xi_{0}}h_{2}+\partial_{\xi_{0}}h_{2}=0,
\]
and for this equation we choose the special solution
\[
h_{2}=0.
\]
Using (\ref{Eq6.13}), equations (\ref{Eq6.15}) can be simplified into
\begin{align}
\partial_{\xi_{2}}A_{0}  & =-A_{0},\nonumber\\
\partial_{\xi_{2}}B_{0}  & =B_{0}.\label{Eq6.15}
\end{align}
Returning to the original variable $u(\xi)$ in the usual way, we have an
approximate solution to the equation (\ref{Eq6.4}) of the form
\begin{equation}
u(\xi)=A(\xi)+B(\xi)e^{-\xi}+\mathcal{O}(\varepsilon^{3}),\label{Eq6.16}
\end{equation}
where the amplitudes $A$ and $B$ are defined by
\begin{align*}
A(\xi)  & =A_{0}(\xi_{1},\xi_{2},...)|_{\xi_{j}=\varepsilon^{j}\xi},\\
B(\xi)  & =B_{0}(\xi_{1},\xi_{2},...)|_{\xi_{j}=\varepsilon^{j}\xi},
\end{align*}
and satisfy the equations
\begin{align}
\frac{dA}{d\xi}  & =\varepsilon A-\varepsilon^{2}A,\nonumber\\
\frac{dB}{d\xi}  & =-\varepsilon B+\varepsilon^{2}B.\label{Eq6.17}
\end{align}
The amplitude equations (\ref{Eq6.17}) are easy to solve, the general solution
is
\begin{align}
A(\xi)  & =Ce^{(\varepsilon-\varepsilon^{2})\xi},\nonumber\\
B(\xi)  & =De^{(-\varepsilon+\varepsilon^{2})\xi},\label{Eq6.18}
\end{align}
where $C$ and $D$ are arbitrary real constants. If we insert the solution
(\ref{Eq6.18}) into (\ref{Eq6.16}) we get
\begin{equation}
u(\xi)=Ce^{(\varepsilon-\varepsilon^{2})\xi}+De^{(-\varepsilon+\varepsilon
^{2}-1)\xi}+\mathcal{O}(\varepsilon^{3}).\label{Eq6.19}
\end{equation}
We now determine the constants $C$ and $D$ such that (\ref{Eq6.19}) satisfy
the boundary conditions to order $\varepsilon^{2}$.
\begin{align*}
u(0)  & =1\text{ \ \ }\Longleftrightarrow\text{ \ }C+D=1,\\
u(\frac{1}{\varepsilon})  & =0\text{ \ \ }\Longleftrightarrow\text{
\ }Ce^{(1-\varepsilon)}+De^{(-1+\varepsilon-\frac{1}{\varepsilon})}=0.
\end{align*}
The linear system for $C$ and $D$ is easy to solve and we get
\begin{align*}
C  & =(1-e^{2-2\varepsilon+\frac{1}{\varepsilon}})^{-1},\\
D  & =(1-e^{-2+2\varepsilon-\frac{1}{\varepsilon}})^{-1},
\end{align*}
and the approximate solution to the original boundary value problem
(\ref{Eq6.1}) is
\begin{equation}
y(x)=(1-e^{2-2\varepsilon+\frac{1}{\varepsilon}})^{-1}e^{(1-\varepsilon
)x}+(1-e^{-2+2\varepsilon-\frac{1}{\varepsilon}})^{-1}e^{(-1+\varepsilon
-\frac{1}{\varepsilon})x}+\mathcal{O}(\varepsilon^{3}).\label{Eq6.20}
\end{equation}
In figure \ref{fig7} we compare a high precision numerical solution of (\ref{Eq6.1}) with the
approximate solution (\ref{Eq6.20}) for $\varepsilon=0.1$. The two solutions
are clearly very close over the whole domain.

\begin{figure}[h]
\centering
\includegraphics[
natheight=3.083100in,
natwidth=4.645800in,
height=2.3091in,
width=3.4696in
]
{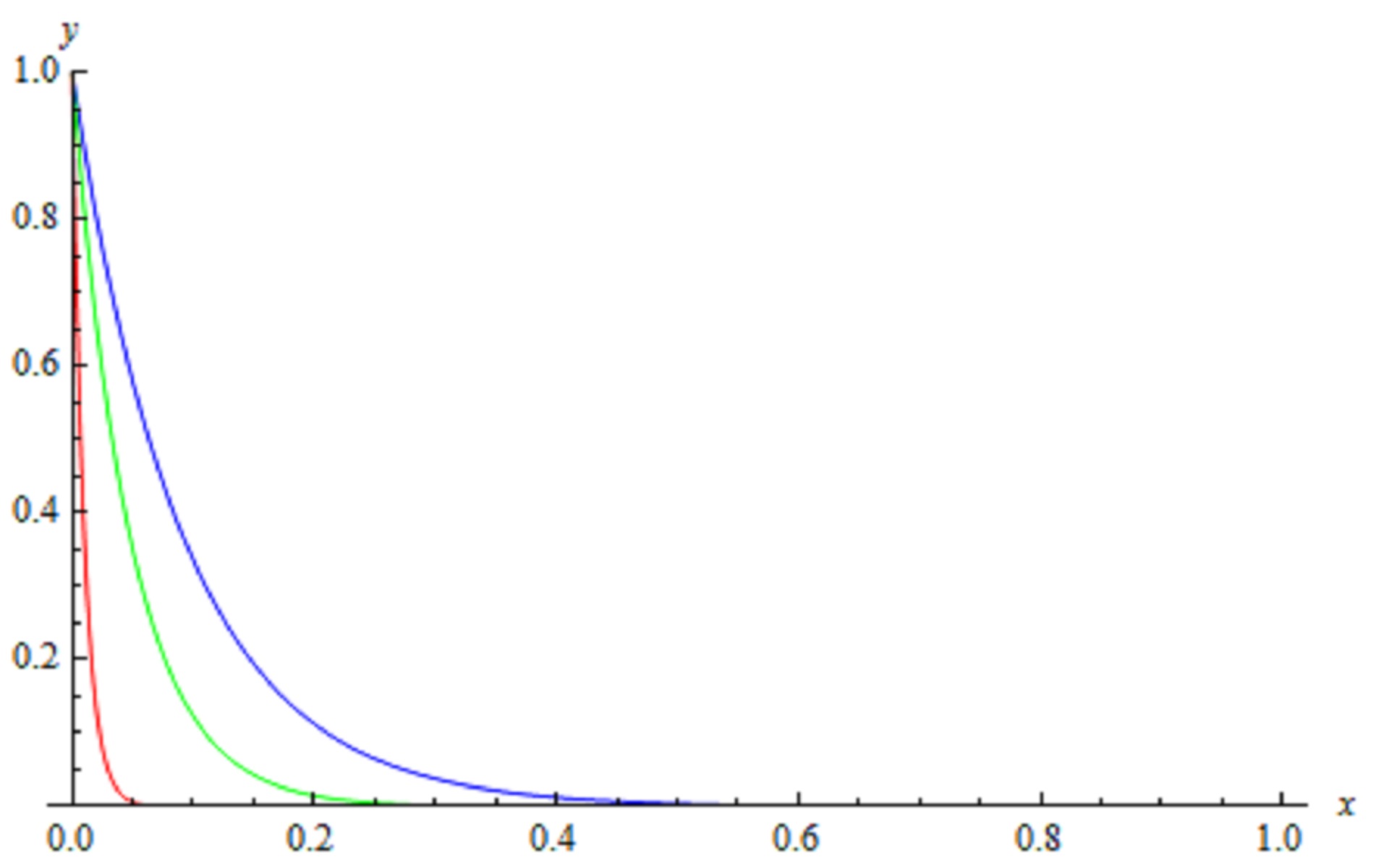}
\caption{A high precision numerical solution of the boundary value problem for
$\varepsilon=0.1$ (Blue), $\varepsilon=0.05$ (Green) and $\varepsilon=0.01$
(Red)}
\label{fig8}
\end{figure}

In figure \ref{fig8} we show a high precision numerical solution to the
boundary value problem (\ref{Eq6.1}) for $\varepsilon=0.1$ (Blue),
$\varepsilon=0.05$ (Green) and $\varepsilon=0.01$ (Red).

We observe that the solution is characterized by a very fast variation close
to $x=0$. The domain close to $x=0$, where $y(x)$ experience a fast variation
is called a \textit{boundary layer}. It's extent is of the order of
$\varepsilon$.

In the context of fluids, the boundary layer is the part of the fluid where
the viscosity plays a role. Away from the boundary layer, the dynamics of the
fluid is to a good approximation described by the Euler equation.

\subsection*{Example 2}

\begin{figure}[h]
\centering
\includegraphics[
natheight=1.656100in,
natwidth=5.083400in,
height=1.6924in,
width=5.1387in
]
{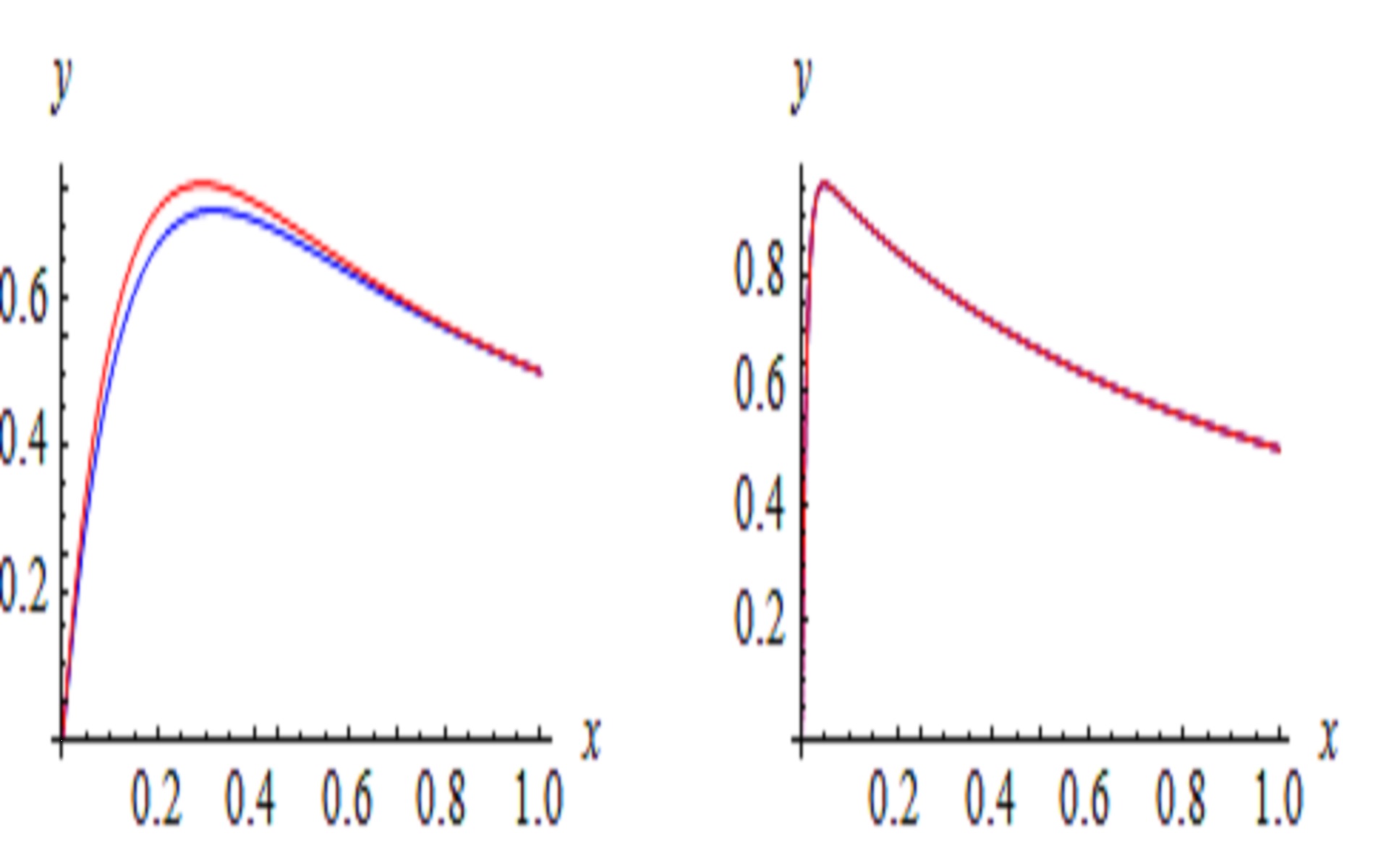}
\caption{Comparing a high precision numerical solution of the boundary value
problem (Red) with the approximate multiple scale solution (Blue) for
$\varepsilon=0.1\,$(Left panel) and $\varepsilon=0.01$ (Right panel)}
\label{fig9}
\end{figure}

Let us consider the following nonlinear boundary value problem
\begin{align}
\varepsilon y^{\prime\prime}+y^{\prime}+y^{2}  & =0,\text{ \ }
0<x<1,\nonumber\\
y(0)  & =0,\nonumber\\
y(1)  & =\frac{1}{2}.\label{Eq6.21}
\end{align}
We recognize that the differential equation in (\ref{Eq6.21}) is singularly
perturbed. The problem is transformed into a regularly perturbed problem using
the transformation
\begin{align}
x  & =\varepsilon\xi,\nonumber\\
y(x)  & =u(\frac{x}{\varepsilon}).\label{Eq6.22}
\end{align}
For the function $u(\xi)$ we get the following regularly perturbed boundary
value problem
\begin{align}
u^{\prime\prime}+u^{\prime}+\varepsilon u^{2}  & =0,\text{ \ \ }0<\xi<\frac
{1}{\varepsilon},\nonumber\\
u(0)  & =0,\nonumber\\
u(\frac{1}{\varepsilon})  & =\frac{1}{2}.\label{Eq6.23}
\end{align}
We have previously, in example 2 in section 5 constructed an approximate
solution to the equation in (\ref{Eq6.23}) that is uniform for $\xi<\frac
{1}{\varepsilon^{2}}$.
\begin{equation}
u(\xi)=A(\xi)+B(\xi)e^{-\xi}-\varepsilon\frac{1}{2}B^{2}(\xi)e^{-2\xi
}+O(\varepsilon^{2}),\label{Eq6.24}
\end{equation}
where the amplitudes $A(\xi)$ and $B(\xi)$ satisfy the equations
\begin{align}
\frac{dA}{d\xi}  & =-\varepsilon A^{2}-2\varepsilon^{2}A^{3},\nonumber\\
\frac{dB}{d\xi}  & =2\varepsilon AB+2\varepsilon^{2}A^{2}B.\label{Eq6.25}
\end{align}
From the boundary conditions on $u(\xi)$, we get
\begin{align}
u(0)  & =0,\text{ \ \ }\Longleftrightarrow\text{ \ \ }A(0)+B(0)-\varepsilon
\frac{1}{2}B^{2}(0)=0,\text{ }\nonumber\\
u(\frac{1}{\varepsilon})  & =\frac{1}{2},\text{ \ \ }\Longleftrightarrow
\text{\ \ \ }A(\frac{1}{\varepsilon})+B(\frac{1}{\varepsilon})e^{-\frac
{1}{\varepsilon}}-\varepsilon\frac{1}{2}B^{2}(\frac{1}{\varepsilon}
)e^{-\frac{2}{\varepsilon}}=\frac{1}{2}.\label{Eq6.26}
\end{align}
The equations (\ref{Eq6.24}),(\ref{Eq6.25}) and (\ref{Eq6.26}) can now be used
to design an efficient numerical algoritm for finding the solution to the
boundary value problem. We do this by defining a function $F(B_{0})$ by
\[
F(B_{0})=A(\frac{1}{\varepsilon})+B(\frac{1}{\varepsilon})e^{-\frac
{1}{\varepsilon}}-\varepsilon\frac{1}{2}B^{2}(\frac{1}{\varepsilon}
)e^{-\frac{2}{\varepsilon}}-\frac{1}{2},
\]
where the functions $A(\xi)$ and $B(\xi)$ are calculated by solving the system
(\ref{Eq6.25}) with initial conditions
\begin{align}
A(0)  & =-B_{0}+\varepsilon\frac{1}{2}B_{0}^{2},\nonumber\\
B(0)  & =B_{0}.\label{Eq6.27}
\end{align}
Using Newton iteration, we find a value of $B_{0}$ such that
\[
F(B_{0})=0.
\]
Inserting this value of $B_{0}$ into the formulas for the initial conditions
(\ref{Eq6.27}), calculating the amplitudes $A(\xi)$, $B(\xi)$ from
(\ref{Eq6.25}) and inserting $A(\xi)$ and $B(\xi)$ into the formula (\ref{Eq6.24}), gives
us a solution to the initial value problem (\ref{Eq6.21}). In figure \ref{fig9} we compare
a high precision numerical solution of (\ref{Eq6.21}) with our approximate
multiple scale solution for $\varepsilon=0.1$ and $\varepsilon=0.01$.

\noindent Apart from being able to use the amplitude equations to construct an
efficient, purely numerical, algoritm for solving the boundary value problem,
it is also possible do quite a lot of analytic work on the amplitude equations (\ref{Eq6.25}). It is
fairly easy to find an explicit formula for $B$ as a function of $A$, it
involves nothing more fancy than using partial fractions. It is also possible
to find an implicit solution for the function $A$, also using partial fractions.

\subsection*{Exercises}

For the following initial value problems, find asymptotic expansions that are
uniform for $t\lesssim\varepsilon^{-3}$. You thus need to take the expansions to
second order in $\varepsilon$. Compare your asymptotic solution to a high
precision numerical solution of the exact problem. Do the comparison for
several values of $\varepsilon$ and show that the asymptotic expansion and the
numerical solution of the exact problem deviates when $t\gtrsim\varepsilon^{-3}$.

\begin{description}
\item[Problem 1:]
\begin{align*}
\frac{d^{2}y}{dt^{2}}+y  &  =\varepsilon y^{2},\\
y(0)  &  =1\\
\frac{dy}{dt}(0)  &  =0
\end{align*}

\item[Problem 2:]
\begin{align*}
\frac{d^{2}y}{dt^{2}}+y  &  =\varepsilon(1-y^{2})\frac{dy}{dt}\\
y(0)  &  =1,\\
\frac{dy}{dt}(0)  &  =0.
\end{align*}

\item[Problem 3:]
\begin{align*}
\frac{d^{2}y}{dt^{2}}+y  &  =\varepsilon(y^{3}-2\frac{dy}{dt}),\\
y(0)  &  =1,\\
\frac{dy}{dt}(0)  &  =0.
\end{align*}

\item[Problem 4:] Let the initial value problem
\begin{align}
\frac{d^{2}y}{dt^{2}}+\frac{dy}{dt}+\varepsilon y^{2}  &  =0,\text{
\ \ }t>0,\nonumber\\
y(0)  &  =1,\nonumber\\
y^{\prime}(0)  &  =1,\label{Eq4.731}
\end{align}
be given. Design a numerical solution to this problem based on the amplitude
equations (\ref{Eq4.72}),(\ref{Eq4.73}) and (\ref{Eq4.71}). Compare this
numerical solution to a high precision numerical solution of (\ref{Eq4.731})
for $t\lesssim\varepsilon^{-3}$. Use several different values of $\varepsilon$ and
show that the multiple scale solution and the high precision solution starts to deviate
when $t\gtrsim\varepsilon^{-3} $.
\end{description}

\section{The multiple scale method for weakly nonlinear PDEs.}

It is now finally time to start applying the multiple scale method to partial
differential equations. The partial differential equations that are of
interest in the science of linear and nonlinear wave motion are almost always
hyperbolic, dispersive and weakly nonlinear. We will therefore focus all our
attention on such equations. 

\subsection*{Example 1}

Let us consider the equation
\begin{equation}
\partial_{tt}u-\partial_{xx}u+u=\varepsilon u^{2}.\label{Eq5.1}
\end{equation}
Inspired by our work on ordinary differential equations, we introduce a
function $h(x_{0},t_{0},x_{1},t_{1},...)$ such that
\begin{equation}
u(x,t)=h(x_{0},t_{0},x_{1},t_{1},...)|_{t_{j}=\varepsilon^{j}t,x_{j}
=\varepsilon^{j}x},\label{Eq5.2}
\end{equation}
is a solution of (\ref{Eq5.1}). The derivatives turns into
\begin{align}
\partial_{t}  &  =\partial_{t_{0}}+\varepsilon\partial_{t_{1}}+\varepsilon
^{2}\partial_{t_{2}}+...\,\,,\nonumber\\
\partial_{x}  &  =\partial_{x_{0}}+\varepsilon\partial_{x_{1}}+\varepsilon
^{2}\partial_{x_{2}}+...\;\;,\label{Eq5.3}
\end{align}
and for $h$ we use the expansion
\begin{equation}
h=h_{0}+\varepsilon h_{1}+\varepsilon^{2}h_{2}+...\;\;.\label{Eq5.4}
\end{equation}
Inserting (\ref{Eq5.2}),(\ref{Eq5.3}) and (\ref{Eq5.4}) and expanding
everything in sight, we get
\begin{gather}
(\partial_{t_{0}}+\varepsilon\partial_{t_{1}}+\varepsilon^{2}\partial_{t_{2}
}+...)(\partial_{t_{0}}+\varepsilon\partial_{t_{1}}+\varepsilon^{2}
\partial_{t_{2}}+...)\nonumber\\
(h_{0}+\varepsilon h_{1}+\varepsilon^{2}h_{2}+...)-\nonumber\\
(\partial_{x_{0}}+\varepsilon\partial_{x_{1}}+\varepsilon^{2}\partial_{x_{2}
}+...)(\partial_{x_{0}}+\varepsilon\partial_{x_{1}}+\varepsilon^{2}
\partial_{x_{2}}+...)\nonumber\\
(h_{0}+\varepsilon h_{1}+\varepsilon^{2}h_{2}+...)+(h_{0}+\varepsilon
h_{1}+\varepsilon^{2}h_{2}+...)\nonumber\\
=\varepsilon(h_{0}+\varepsilon h_{1}+\varepsilon^{2}h_{2}+...)^{2},\nonumber\\
\Downarrow\nonumber\\
(\partial_{t_{0}t_{0}}+\varepsilon(\partial_{t_{0}t_{1}}+\partial_{t_{1}t_{0}
})+\varepsilon^{2}(\partial_{t_{0}t_{2}}+\partial_{t_{1}t_{1}}+\partial
_{t_{2}t_{0}})+...)\nonumber\\
(h_{0}+\varepsilon h_{1}+\varepsilon^{2}h_{2}+...)-\nonumber\\
(\partial_{x_{0}x_{0}}+\varepsilon(\partial_{x_{0}x_{1}}+\partial_{x_{1}x_{0}
})+\varepsilon^{2}(\partial_{x_{0}x_{2}}+\partial_{x_{1}x_{1}}+\partial
_{x_{2}x_{0}})+...)\nonumber\\
(h_{0}+\varepsilon h_{1}+\varepsilon^{2}h_{2}+...)+(h_{0}+\varepsilon
h_{1}+\varepsilon^{2}h_{2}+...)\nonumber\\
=\varepsilon(h_{0}^{2}+2\varepsilon h_{0}h_{1}+...),\nonumber\\
\Downarrow\nonumber
\end{gather}
\begin{gather}
\partial_{t_{0}t_{0}}h_{0}+\varepsilon(\partial_{t_{0}t_{0}}h_{1}
+\partial_{t_{0}t_{1}}h_{0}+\partial_{t_{1}t_{0}}h_{0})+\nonumber\\
\varepsilon^{2}(\partial_{t_{0}t_{0}}h_{2}+\partial_{t_{0}t_{1}}h_{1}
+\partial_{t_{1}t_{0}}h_{1}+\partial_{t_{0}t_{2}}h_{0}+\partial_{t_{1}t_{1}
}h_{0}+\partial_{t_{2}t_{0}}h_{0})-...\nonumber\\
\partial_{x_{0}x_{0}}h_{0}-\varepsilon(\partial_{x_{0}x_{0}}h_{1}
+\partial_{x_{0}x_{1}}h_{0}+\partial_{x_{1}x_{0}}h_{0})-\nonumber\\
\varepsilon^{2}(\partial_{x_{0}x_{0}}h_{2}+\partial_{x_{0}x_{1}}h_{1}
+\partial_{x_{1}x_{0}}h_{1}+\partial_{x_{0}x_{2}}h_{0}+\partial_{x_{1}x_{1}
}h_{0}+\partial_{x_{2}x_{0}}h_{0})\nonumber\\
+h_{0}+\varepsilon h_{1}+\varepsilon^{2}h_{2}+...\nonumber\\
=\varepsilon h_{0}^{2}+2\varepsilon^{2}h_{0}h_{1}+...\;\;,\label{Eq5.5}
\end{gather}
which gives us the perturbation hierarchy
\begin{gather}
\partial_{t_{0}t_{0}}h_{0}-\partial_{x_{0}x_{0}}h_{0}+h_{0}=0,\label{Eq5.6}\\
\nonumber\\
\partial_{t_{0}t_{0}}h_{1}-\partial_{x_{0}x_{0}}h_{1}+h_{1}=h_{0}^{2}
-\partial_{t_{0}t_{1}}h_{0}-\partial_{t_{1}t_{0}}h_{0}\nonumber\\
+\partial_{x_{0}x_{1}}h_{0}+\partial_{x_{1}x_{0}}h_{0},\label{Eq5.7}\\
\nonumber\\
\partial_{t_{0}t_{0}}h_{2}-\partial_{x_{0}x_{0}}h_{2}+h_{2}=2h_{0}
h_{1}-\partial_{t_{0}t_{1}}h_{1}-\partial_{t_{1}t_{0}}h_{1}\nonumber\\
-\partial_{t_{0}t_{2}}h_{0}-\partial_{t_{1}t_{1}}h_{0}-\partial_{t_{2}t_{0}
}h_{0}+\partial_{x_{0}x_{1}}h_{1}+\partial_{x_{1}x_{0}}h_{1}\nonumber\\
+\partial_{x_{0}x_{2}}h_{0}+\partial_{x_{1}x_{1}}h_{0}+\partial_{x_{2}x_{0}
}h_{0}.\label{Eq5.8}
\end{gather}
For ordinary differential equations, we used the general solution to the order
$\varepsilon^{0}$ equation. For partial differential equations we can not do
this. We will rather use a finite sum of linear modes. The simplest
possibility is a single linear mode which we use here
\begin{equation}
h_{0}(t_{0},x_{0},t_{1},x_{1},...)=A_{0}(t_{1},x_{1},...)e^{i(kx_{0}-\omega
t_{0})}+(\ast).\label{Eq5.9}
\end{equation}
Since we are not using the general solution, we will in not be able to satisfy
arbitrary initial conditions. However, in the theory of waves this is
perfectly alright, since most of the time the relevant initial conditions are
in fact finite sums of wave packets or even a single wave packet. Such initial
conditions can be included in the multiple scale approach that we discuss in
this section. For (\ref{Eq5.9}) to actually be a solution to (\ref{Eq5.6}) we
must have
\begin{equation}
\omega=\omega(k)=\sqrt{1+k^{2}},\label{Eq5.10}
\end{equation}
which we of course recognize as the dispersion relation for the linearized
version of (\ref{Eq5.1}). With the choise of signs used here, (\ref{Eq5.9})
will represent a right-moving disturbance.

Inserting (\ref{Eq5.9}) into (\ref{Eq5.7}) we get
\begin{gather}
\partial_{t_{0}t_{0}}h_{1}-\partial_{x_{0}x_{0}}h_{1}+h_{1}=2|A_{0}
|^{2}\nonumber\\
+A_{0}^{2}e^{2i(kx_{0}-\omega t_{0})}+A_{0}^{\ast2}e^{-2i(kx_{0}-\omega
t_{0})}\nonumber\\
+(2i\omega\partial_{t_{1}}A_{0}+2ik\partial_{x_{1}}A_{0})e^{i(kx_{0}-\omega
t_{0})}\nonumber\\
-(2i\omega\partial_{t_{1}}A_{0}^{\ast}+2ik\partial_{x_{1}}A_{0}^{\ast
})e^{-i(kx_{0}-\omega t_{0})}.\label{Eq5.11}
\end{gather}
In order to remove secular terms, we must postulate that
\begin{gather}
2i\omega\partial_{t_{1}}A_{0}+2ik\partial_{x_{1}}A_{0}=0,\nonumber\\
\Updownarrow\nonumber\\
\partial_{t_{1}}A_{0}=-\frac{k}{\omega}\partial_{x_{1}}A_{0}.\label{Eq5.12}
\end{gather}
Here we assume that the terms
\[
e^{2i(kx_{0}-\omega t_{0})},e^{-2i(kx_{0}-\omega t_{0})}\;\;,
\]
are \textit{not} solutions to the homogenous equation
\[
\partial_{t_{0}t_{0}}h_{1}-\partial_{x_{0}x_{0}}h_{1}+h_{1}=0.
\]
For this to be true we must have
\begin{equation}
\omega(2k)\neq2\omega(k),\label{Eq5.12.1}
\end{equation}
and this is in fact true for all $k$. This is however not generally true for
dispersive wave equations. Whether it is true or not will depend on the exact
form of the dispersion relation for the system of interest. In the theory of
interacting waves, equality in (\ref{Eq5.12.1}), is called \textit{phase
matching}, and is of outmost importance.

The equation for $h_{1}$ now simplify into
\begin{equation}
\partial_{t_{0}t_{0}}h_{1}-\partial_{x_{0}x_{0}}h_{1}+h_{1}=2|A_{0}|^{2}
+A_{0}^{2}e^{2i(kx_{0}-\omega t_{0})}+A_{0}^{\ast2}e^{-2i(kx_{0}-\omega
t_{0})}.\label{Eq5.13}
\end{equation}
According to the rules of the game we need a special solution to this
equation. It is easy to verify that
\begin{equation}
h_{1}=2|A_{0}|^{2}-\frac{1}{3}A_{0}^{2}e^{2i(kx_{0}-\omega t_{0})}-\frac{1}
{3}A_{0}^{\ast2}e^{-2i(kx_{0}-\omega t_{0})},\label{Eq5.14}
\end{equation}
is such a special solution. Inserting (\ref{Eq5.9}) and (\ref{Eq5.14}) into
(\ref{Eq5.8}), we get
\begin{align}
\partial_{t_{0}t_{0}}h_{2}-\partial_{x_{0}x_{0}}h_{2}+h_{2}  &  =(2i\omega
\partial_{t_{2}}A_{0}+2ik\partial_{x_{2}}A_{0}-\partial_{t_{1}t_{1}}
A_{0}\label{Eq5.15}\\
&  +\partial_{x_{1}x_{1}}A_{0}+\frac{10}{3}|A_{0}|^{2}A_{0})e^{i(kx_{0}-\omega
t_{0})}+NST+(\ast).\nonumber
\end{align}
In order to remove secular terms we must postulate that
\begin{equation}
2i\omega\partial_{t_{2}}A_{0}+2ik\partial_{x_{2}}A_{0}-\partial_{t_{1}t_{1}
}A_{0}+\partial_{x_{1}x_{1}}A_{0}+\frac{10}{3}|A_{0}|^{2}A_{0}=0.\label{Eq5.16}
\end{equation}
Equations (\ref{Eq5.12}) and (\ref{Eq5.16}) is, as usual, an overdetermined
system. In general it is not an easy matter to verify that an overdetermined
system of partial differential equations is solvable and the methods that do
exist to adress such questions are mathematically quite sophisticated. For the
particular case discussed here it is however easy to verify that the system is
in fact solvable. But, as we have stressed several times in these lecture
notes, we are not really concerned with the solvability of the system
(\ref{Eq5.12}), (\ref{Eq5.16}) for the many variabable function $A_{0} $. We
are rather interested in the function $u(x,t)$ which is a solution to
(\ref{Eq5.1}). With that in mind, we define an amplitude
\begin{equation}
A(x,t)=A_{0}(t_{1},x_{1},...)|_{t_{j}=\varepsilon^{j}t,x_{j}=\varepsilon^{j}
x}.\label{Eq5.17}
\end{equation}
The solution to (\ref{Eq5.1}) is then
\begin{align}
u(x,t)  &  =A(x,t)e^{i(kx-\omega t)}+\varepsilon(2|A|^{2}(x,t)-\frac{1}
{3}A^{2}(x,t)e^{2i(kx_{0}-\omega t_{0})}\nonumber\\
&  -\frac{1}{3}A^{\ast2}(x,t)e^{-2i(kx_{0}-\omega t_{0})})+O(\varepsilon^{2}),\label{Eq5.18}
\end{align}
where $A(x,t)$ satisfy a certain amplitude equation that we will now derive.

Multiplying equation (\ref{Eq5.12}) by $\varepsilon$ , equation (\ref{Eq5.16})
by $\varepsilon^{2}$ and adding the two expressions, we get
\begin{gather}
\varepsilon(2i\omega\partial_{t_{1}}A_{0}+2ik\partial_{x_{1}}A_{0})\nonumber\\
+\varepsilon^{2}(2i\omega\partial_{t_{2}}A_{0}+2ik\partial_{x_{2}}
A_{0}-\partial_{t_{1}t_{1}}A_{0}+\partial_{x_{1}x_{1}}A_{0}+\frac{10}{3}
|A_{0}|^{2}A_{0})=0,\nonumber\\
\Downarrow\nonumber\\
2i\omega(\partial_{t_{0}}+\varepsilon\partial_{t_{1}}+\varepsilon^{2}
\partial_{t_{2}})A_{0}+2ik(\partial_{x_{0}}+\varepsilon\partial_{x_{1}
}+\varepsilon^{2}\partial_{x_{2}})A_{0}\nonumber\\
-(\partial_{t_{0}}+\varepsilon\partial_{t_{1}}+\varepsilon^{2}\partial_{t_{2}
})^{2}A_{0}+(\partial_{x_{0}}+\varepsilon\partial_{x_{1}}+\varepsilon
^{2}\partial_{x_{2}})^{2}A_{0}+\varepsilon^{2}\frac{10}{3}|A_{0}|^{2}
A_{0}=0,\label{Eq5.19}
\end{gather}
where we have used the fact that $A_{0}$ does not depend on $t_{0}$ and
$x_{0}$ and where the equation (\ref{Eq5.19}) is correct to second order in
$\varepsilon$. \ If we now evaluate (\ref{Eq5.19}) at $x_{j}=\varepsilon
^{j}x,t_{j}=\varepsilon^{j}t$, using (\ref{Eq5.3}) and (\ref{Eq5.17}), we
get the amplitude equation

\begin{gather}
2i\omega\partial_{t}A+2ik\partial_{x}A-\partial_{tt}A+\partial_{xx}
A+\varepsilon^{2}\frac{10}{3}|A|^{2}A=0,\nonumber\\
\Updownarrow\nonumber\\
\partial_{t}A=-\frac{k}{\omega}\partial_{x}A-\frac{i}{2\omega}\partial
_{tt}A+\frac{i}{2\omega}\partial_{xx}A+\varepsilon^{2}\frac{5i}{3\omega
}|A|^{2}A.\label{Eq5.21}
\end{gather}
This equation appears to have a problem since it contains a second
derivative with respect to time. The initial conditions for (\ref{Eq5.1}) is
only sufficient to determine $A(x,0)$. However, in order to be consistent with
the multiple scale procedure leading up to (\ref{Eq5.21}) we can only consider
solutions such that
\begin{gather}
\partial_{t}A\sim-\frac{k}{\omega}\partial_{x}A\sim\varepsilon,\nonumber\\
\Downarrow\nonumber\\
\partial_{tt}A\sim\left(  \frac{k}{\omega}\right)  ^{2}\partial_{xx}
A\sim\varepsilon^{2}.\label{Eq5.22}
\end{gather}
Thus we can, to second order in $\varepsilon$, rewrite the amplitude equation
as
\begin{equation}
\partial_{t}A=-\frac{k}{\omega}\partial_{x}A+\frac{i}{2\omega^{3}}
\partial_{xx}A+\varepsilon^{2}\frac{5i}{3\omega}|A|^{2}A.\label{Eq5.24}
\end{equation}
This is now first order in time and has a unique solution for a given initial
condition $A(x,0)$.

The multiple scale procedure demands that the amplitude $A(x,t)$ vary slowly
on scales $L=\frac{2\pi}{k},T=\frac{2\pi}{\omega}$. This means that
(\ref{Eq5.18}) and (\ref{Eq5.24}) can be thought of as a fast numerical scheme
for \textit{wavepackets} solutions to (\ref{Eq5.1}). If these are the kind of
solutions that we are interested in, and in the theory of waves this is often
the case, it is much more efficient to use (\ref{Eq5.18}) and (\ref{Eq5.24})
rather than having to resolve the scales $L$ and $T$ by integrating the
original equation (\ref{Eq5.1}).

The very same equation (\ref{Eq5.24}) appear as leading order amplitude
equation starting from a large set of nonlinear partial differential equations
describing a wide array of physical phenomena in fluid dynamics, clima
science, laser physics etc. The equation appeared for the first time more than
70 years ago, but it was not realized at the time that the Nonlinear
Schr\"{o}dinger equation (NLS), as it is called, is very special indeed.

V. Zakharov discovered in 1974 that NLS is in a certain sense completely
solvable. He discovered a nonlinear integral transform that decompose NLS into
an infinite system of uncoupled ODE's, that in many important cases are easy
to solve. This transform is called the \textit{Scattering Transform}.

Using this transform one can find explicit formulas for solutions of NLS that
acts like particles, they are localized disturbances in a wavefield that does
not disperse and they collide elastically just like particles do. The NLS equation
has a host of interesting and beautiful properties. It has for example
infinitely many quantities that are concerved under the time evolution and is
the continuum analog of a \textit{completely integrable} system of ODE's.

Many books and $\infty$- many papers have been written about this equation. In
the process of doing this, many other equations having similar wonderful
properties has been discovered. They \textit{all} appear through the use of
the method of multiple scales. However, all these wonderful properties,
however nice they are, are not robust. If we want to propagate our waves for
$t\lesssim\varepsilon^{-4}$, the multiple scale procedure must be extended to
order $\varepsilon^{3}$, and additional terms will appear in the amplitude
equation. These additional terms will destroy many of the wonderful
mathematical properties of the Nonlinear Schr\"{o}dinger equation but it will
\textit{not} destroy the fact that it is the key element in a fast numerical
scheme for wave packet solutions to (\ref{Eq5.1}).

\subsection*{Example 2}

Let us consider the equation
\begin{equation}
u_{tt}+u_{xx}+u_{xxxx}+u=\varepsilon u^{3}.\label{Eq5.25}
\end{equation}
Introducing the usual tools for the multiple scale method, we have
\begin{align}
u(x,t)  &  =h(x_{0},t_{0},x_{1},t_{1},...)|_{t_{j}=\varepsilon^{j}
t,x_{j}=\varepsilon^{j}x},\nonumber\\
\partial_{t}  &  =\partial_{t_{0}}+\varepsilon\partial_{t_{1}}+...\;\;,\nonumber\\
\partial_{x}  &  =\partial_{x_{0}}+\varepsilon\partial_{x_{1}}+...\;\;,\nonumber\\
h  &  =h_{0}+\varepsilon h_{1}+...\;\;.\label{Eq5.26}
\end{align}
Inserting these expressions into (\ref{Eq5.25}) and expanding we get
\begin{gather*}
(\partial_{t_{0}}+\varepsilon\partial_{t_{1}}+...)(\partial_{t_{0}
}+\varepsilon\partial_{t_{1}}+...)(h_{0}+\varepsilon h_{1}+...)+\\
(\partial_{x_{0}}+\varepsilon\partial_{x_{1}}+...)(\partial_{x_{0}
}+\varepsilon\partial_{x_{1}}+...)(h_{0}+\varepsilon h_{1}+...)+\\
(\partial_{x_{0}}+\varepsilon\partial_{x_{1}}+...)(\partial_{x_{0}
}+\varepsilon\partial_{x_{1}}+...)\\
(\partial_{x_{0}}+\varepsilon\partial_{x_{1}}+...)(\partial_{x_{0}
}+\varepsilon\partial_{x_{1}}+...)(h_{0}+\varepsilon h_{1}+...)\\
=\varepsilon(h_{0}+...)^{3},\\
\Downarrow\\
(\partial_{t_{0}t_{0}}+\varepsilon(\partial_{t_{0}t_{1}}+\partial_{t_{1}t_{0}
})+...)(h_{0}+\varepsilon h_{1}+...)+\\
(\partial_{x_{0}x_{0}}+\varepsilon(\partial_{x_{0}x_{1}}+\partial_{x_{1}x_{0}
})+...)(h_{0}+\varepsilon h_{1}+...)+\\
(\partial_{x_{0}x_{0}}+\varepsilon(\partial_{x_{0}x_{1}}+\partial_{x_{1}x_{0}
})+...)(\partial_{x_{0}x_{0}}+\varepsilon(\partial_{x_{0}x_{1}}+\partial
_{x_{1}x_{0}})+...)\\
(h_{0}+\varepsilon h_{1}+...)=\varepsilon h_{0}^{3}+...\;\;,\\
\Downarrow\\
\partial_{t_{0}t_{0}}h_{0}+\varepsilon(\partial_{t_{0}t_{0}}h_{1}
+\partial_{t_{0}t_{1}}h_{0}+\partial_{t_{1}t_{0}}h_{0})+\\
\partial_{x_{0}x_{0}}h_{0}+\varepsilon(\partial_{x_{0}x_{0}}h_{1}
+\partial_{x_{0}x_{1}}h_{0}+\partial_{x_{1}x_{0}}h_{0})+\\
\partial_{x_{0}x_{0}x_{0}x_{0}}h_{0}+\varepsilon(\partial_{x_{0}x_{0}
x_{0}x_{0}}h_{1}+\partial_{x_{0}x_{0}x_{0}x_{1}}h_{0}+\partial_{x_{0}
x_{0}x_{1}x_{0}}h_{0}\\
+\partial_{x_{0}x_{1}x_{0}x_{0}}h_{0}+\partial_{x_{1}x_{0}x_{0}x_{0}}
h_{0})+...\\
=\varepsilon h_{0}^{3}+...\;\;,
\end{gather*}
which gives us the perturbation hierarchy
\begin{gather}
\partial_{t_{0}t_{0}}h_{0}+\partial_{x_{0}x_{0}}h_{0}+\partial_{x_{0}
x_{0}x_{0}x_{0}}h_{0}=0,\label{Eq5.27}\\
\nonumber\\
\partial_{t_{0}t_{0}}h_{1}+\partial_{x_{0}x_{0}}h_{1}+\partial_{x_{0}
x_{0}x_{0}x_{0}}h_{1}=h_{0}^{3}\label{Eq5.28}\\
-\partial_{t_{0}t_{1}}h_{0}-\partial_{t_{1}t_{0}}h_{0}-\partial_{x_{0}x_{1}
}h_{0}-\partial_{x_{1}x_{0}}h_{0}\nonumber\\
-\partial_{x_{0}x_{0}x_{0}x_{1}}h_{0}-\partial_{x_{0}x_{0}x_{1}x_{0}}
h_{0}-\partial_{x_{0}x_{1}x_{0}x_{0}}h_{0}+\partial_{x_{1}x_{0}x_{0}x_{0}
}h_{0}.\nonumber
\end{gather}
For the order $\varepsilon^{0}$ equation, we choose a wave packet solution
\begin{equation}
h_{0}(x_{0},t_{0},x_{1},t_{1},...)=A_{0}(x_{1},t_{1},...)e^{i(kx_{0}-\omega
t_{0})}+(\ast),\label{Eq5.29}
\end{equation}
where the dispersion relation is
\begin{equation}
\omega=\sqrt{k^{4}-k^{2}+1}.\label{Eq5.30}
\end{equation}
Inserting (\ref{Eq5.29}) into (\ref{Eq5.28}), we get after a few algebraic
manipulations
\begin{gather}
\partial_{t_{0}t_{0}}h_{1}+\partial_{x_{0}x_{0}}h_{1}+\partial_{x_{0}
x_{0}x_{0}x_{0}}h_{1}=\nonumber\\
(2i\omega\partial_{t_{1}}A_{0}-2ik\partial_{x_{1}}A_{0}+4ik^{3}\partial
_{x_{1}}A_{0}+3|A_{0}|^{2}A_{0})e^{i(kx_{0}-\omega t_{0})}\nonumber\\
+A_{0}^{3}e^{3i(kx_{0}-\omega t_{0})}+(\ast).\label{Eq5.31}
\end{gather}
In order to remove secular terms we must postulate that
\begin{equation}
2i\omega\partial_{t_{1}}A_{0}-2ik\partial_{x_{1}}A_{0}+4ik^{3}\partial_{x_{1}
}A_{0}+3|A_{0}|^{2}A_{0}=0.\label{Eq5.32}
\end{equation}
But using the dispersion relation (\ref{Eq5.30}), we have
\[
-2ik+4ik^{3}=2i\omega\omega^{\prime},
\]
so that (\ref{Eq5.32}) simplifies into
\begin{equation}
2i\omega(\partial_{t_{1}}A_{0}+\omega^{\prime}\partial_{x_{1}}A_{0}
)+3|A_{0}|^{2}A_{0}=0.\label{Eq5.33}
\end{equation}
Introducing an amplitude
\[
A(x,t)=A_{0}(x_{1},t_{1},...)|_{x_{j}=e^{j}x,t_{j}=\varepsilon^{j}t},
\]
we get, following the approach from the previous example, the amplitude
equation
\begin{equation}
2i\omega(\partial_{t}A+\omega^{\prime}\partial_{x}A)=-3|A|^{2}
A.\label{Eq5.34}
\end{equation}
This equation together with the expansion
\begin{equation}
u(x,t)=A(t)e^{i(kx-\omega t)}+(\ast)+O(\varepsilon),\label{Eq5.35}
\end{equation}
constitute a fast numerical scheme for wave packet solutions to (\ref{Eq5.25})
for $t\lesssim\varepsilon^{-2}$. Of course, this particular amplitude equation
can be solved analytically, but as stressed earlier, this property is not
robust and can easily be lost if we take the expansion to higher order in
$\varepsilon$.

There is however one point in our derivation that we need to look more closely
into. We assumed that the term
\begin{equation}
A_{0}^{3}e^{3i(kx_{0}-\omega t_{0})},\label{Eq5.36}
\end{equation}
was \textit{not} a secular term. The term \textit{is} secular if
\begin{equation}
\omega(3k)=3\omega(k).\label{Eq5.37}
\end{equation}
Using the dispersion relation (\ref{Eq5.30}) we have
\begin{align}
\omega(3k) &  =3\omega(k),\nonumber\\
&  \Updownarrow\nonumber\\
\sqrt{81k^{4}-9k^{2}+1} &  =3\sqrt{k^{4}-k^{2}+1},\nonumber\\
&  \Updownarrow\nonumber\\
81k^{4}-9k^{2}+1 &  =9k^{4}-9k^{2}+9,\nonumber\\
&  \Updownarrow\nonumber\\
k &  =\pm\frac{1}{\sqrt{3}}.\label{Eq5.38}
\end{align}
Thus the term (\ref{Eq5.36}) \textit{can} be secular if the wave number of the
wave packet is given by (\ref{Eq5.38}). This is another example of the
fenomenon that we in the theory of interacting waves call phase matching. As
long as we stay away from the two particular values of the wave numbers given
in (\ref{Eq5.38}), our expansion (\ref{Eq5.34}) and (\ref{Eq5.35}) is uniform
for $t\lesssim\varepsilon^{-2}$. However if the wave number takes on one of
the two values in (\ref{Eq5.38}), nonuniformities will make the ordering of
the expansion break down for $t\sim\varepsilon^{-1}$. However this does not
mean that the multiple scale method breaks down. We only need to include a
second amplitude at order $\varepsilon^{0}$ that we can use to remove the
additional secular terms at order $\varepsilon^{1}$. We thus, instead of
(\ref{Eq5.29}), use the solution
\begin{align}
h_{0}(x_{0},t_{0},x_{1},t_{1},...) &  =A_{0}(x_{1},t_{1},...)e^{i(kx_{0}
-\omega t_{0})}\nonumber\\
&  +B_{0}(x_{1},t_{1},...)e^{3i(kx_{0}-\omega t_{0})}+(\ast),\label{Eq5.39}
\end{align}
where $k$ now is given by (\ref{Eq5.38}). Inserting this expression for
$h_{0}$ into the order $\varepsilon$ equation (\ref{Eq5.28}) we get, after a
fair amount of algebra, the equation
\begin{gather}
\partial_{t_{0}t_{0}}h_{1}+\partial_{x_{0}x_{0}}h_{1}+\partial_{x_{0}
x_{0}x_{0}x_{0}}h_{1}=\label{Eq5.40}\\
(2i\omega\partial_{t_{1}}A_{0}-2ik\partial_{x_{1}}A_{0}+4ik^{3}\partial
_{x_{1}}A_{0}\nonumber\\
+3|A_{0}|^{2}A_{0}+6|B_{0}|^{2}A_{0}+3A_{0}^{\ast2}B_{0})e^{i(kx_{0}-\omega
t_{0})}\nonumber\\
+(6i\omega\partial_{t_{1}}B_{0}-6ik\partial_{x_{1}}B_{0}+108ik^{3}
\partial_{x_{1}}B_{0}\nonumber\\
+3|B_{0}|^{2}B_{0}+6|A_{0}|^{2}B_{0}+A_{0}^{3})e^{3i(kx_{0}-\omega t_{0}
)}\nonumber\\
+NST+(\ast).\nonumber
\end{gather}
In order to remove secular terms we must postulate the two equations
\begin{gather}
2i\omega\partial_{t_{1}}A_{0}-2ik\partial_{x_{1}}A_{0}+4ik^{3}\partial_{x_{1}
}A_{0}\nonumber\\
+3|A_{0}|^{2}A_{0}+6|B_{0}|^{2}A_{0}+3A_{0}^{\ast2}B_{0}=0,\nonumber\\
\nonumber\\
6i\omega\partial_{t_{1}}B_{0}-6ik\partial_{x_{1}}B_{0}+108ik^{3}
\partial_{x_{1}}B_{0}\nonumber\\
+3|B_{0}|^{2}B_{0}+6|A_{0}|^{2}B_{0}+A_{0}^{3}=0.\label{Eq5.41}
\end{gather}
Using the dispersion relation we have
\begin{align*}
& \\
-6ik+108ik^{3}  & =2i\omega(3k)\omega^{\prime}(3k).
\end{align*}
Inserting this into the system (\ref{Eq5.41}), simplifies it into
\begin{align}
2i\omega(k)(\partial_{t_{1}}A_{0}+\omega^{\prime}(k)\partial_{x_{1}}A_{0}) &
=-3|A_{0}|^{2}A_{0}-6|B_{0}|^{2}A_{0}-3A_{0}^{\ast2}B_{0},\nonumber\\
2i\omega(3k)(\partial_{t_{1}}B_{0}+\omega^{\prime}(3k)\partial_{x_{1}}B_{0})
&  =-3|B_{0}|^{2}B_{0}-6|A_{0}|^{2}B_{0}-A_{0}^{3}.\label{Eq5.42}
\end{align}
Introducing amplitudes
\begin{align}
A(x,t) &  =A_{0}(x_{1},t_{1},...)|_{x_{j}=e^{j}x,t_{j}=\varepsilon^{j}
t},\nonumber\\
B(x,t) &  =B_{0}(x_{1},t_{1},...)|_{x_{j}=e^{j}x,t_{j}=\varepsilon^{j}
t},\label{Eq5.43}
\end{align}
the asymptotic expansion and corresponding amplitude equations for this case
are found to be
\begin{align}
u(x,t) &  =A(x,t)e^{i(kx-\omega t)}\nonumber\\
&  +B(x,t)e^{3i(kx-\omega t)}+(\ast)+O(\varepsilon),\nonumber\\
2i\omega(k)(\partial_{t}A+\omega^{\prime}(k)\partial_{x}A) &  =-3|A|^{2}
A-6|B|^{2}A-3A^{\ast2}B,\nonumber\\
2i\omega(3k)(\partial_{t}B+\omega^{\prime}(3k)\partial_{x}B) &  =-3|B|^{2}
B+6|A|^{2}B+A^{3}.\label{Eq5.44}
\end{align}
The same approach must be used to treat the case when we do not have exact
phase matching but we still have
\[
\omega(3k)\approx3\omega(k)
\]

\subsection*{Exercises}

In the following problems, use the methods from this section to find asymptotic
expansions that are uniform for $t\lesssim\varepsilon^{-2}$. Thus all
expansions must be taken to second order in $\varepsilon$.

\begin{description}
\item[Problem 1:]
\[
u_{tt}-u_{xx}+u=\varepsilon^{2}u^{3},
\]

\item[Problem 2:]
\[
u_{tt}-u_{xx}+u=\varepsilon(u^{2}+u_{x}^{2}),
\]

\item[Problem 3:]
\[
u_{tt}-u_{xx}+u=\varepsilon(uu_{xx}-u^{2}),
\]

\item[Problem 4:]
\[
u_{t}+u_{xxx}=\varepsilon u^{2}u_{x},
\]

\item[Problem 5:]
\[
u_{tt}-u_{xx}+u=\varepsilon(u_{x}^{2}-uu_{xx}).
\]

\end{description}

\section{The multiple scale method for Maxwell's equations}

In optics the equations of interest are of course Maxwell's equations. For a
situation without free carges and currents they are given by
\begin{align}
\partial_{t}\mathbf{B+\nabla\times E}  & =0,\nonumber\\
\partial_{t}\mathbf{D-\nabla\times H}  & =0,\nonumber\\
\nabla\cdot\mathbf{D}  & =0,\nonumber\\
\nabla\cdot\mathbf{B}  & =0.\label{Eq8.1}
\end{align}
At optical frequencies, materials of interest is almost always nonmagnetic so that we have
\begin{align}
\mathbf{H}  & =\frac{1}{\mu}\mathbf{B},\nonumber\\
\mathbf{D}  & =\varepsilon_{0}\mathbf{E}+\mathbf{P}.\label{Eq8.2}
\end{align}
The polarization is in general a sum of a terms that is linear in $\mathbf{E}
$ and one that is nonlinear in $\mathbf{E}$. We have
\begin{equation}
\mathbf{P}=\mathbf{P}_{L}+\mathbf{P}_{NL},\label{Eq8.3}
\end{equation}
where the term linear in $\mathbf{E}$ has the general form
\begin{equation}
\mathbf{P}_{L}(\mathbf{x},t)=\varepsilon_{0}\int_{-\infty}^{t}dt^{\prime}
\chi(t-t^{\prime})\mathbf{E}(\mathbf{x},t^{\prime}).\label{Eq8.4}
\end{equation}
Thus the polarization at a time $t$ depends on the electric field at all times
previous to $t$. This memory effect is what we in optics call \textit{temporal
dispersion}. The presence of dispersion in Maxwell equations spells trouble
for the intergration of the equations in time; we can not solve them as a
standard initial value problem. This is of course well known in optics and
various more or less ingenious methods has been designed for getting around
this problem. In optical pulse propagation one gets around the problem by
solving Maxwell's equations approximately as a boundary value problem rather
than as an initial value problem. A very general version of this approach is
the well known UPPE \cite{UPPE}\cite{per1} propagation scheme. In
these lecture notes we will, using the multiple scale method, derive
approximations to Maxwell's equations that can be solved as an initial value problem.

In the explicite calculations that we do we will assume that the nonlinear
polarization is generated by the Kerr effect. Thus we will assume that
\begin{equation}
\mathbf{P}_{NL}=\varepsilon_{0}\eta\mathbf{E}\cdot\mathbf{EE},\label{Eq8.5}
\end{equation}
where $\eta$ is the Kerr coefficient. This is a choise we make just to be
specific, the applicability of the multiple scale method to Maxwell's
equations in no way depend on this particular choise for the nonlinear response.

Before we proceed with the multiple scale method we will introduce a more
convenient representation of the dispersion. Observe that we have
\begin{align}
\mathbf{P}_{L}(\mathbf{x},t)  & =\varepsilon_{0}\int_{-\infty}^{t}dt^{\prime
}\chi(t-t^{\prime})\mathbf{E}(\mathbf{x},t^{\prime}),\nonumber\\
& =\varepsilon_{0}\int_{-\infty}^{\infty}d\omega\widehat{\chi}(\omega
)\widehat{\mathbf{E}}(\mathbf{x},\omega)e^{-i\omega t},\nonumber\\
& =\varepsilon_{0}\int_{-\infty}^{\infty}d\omega\left(
{\displaystyle\sum_{n=0}^{\infty}}
\frac{\widehat{\chi}^{(n)}(0)}{n!}\omega^{n}\right)  \widehat{\mathbf{E}
}(\mathbf{x},\omega)e^{-i\omega t},\nonumber\\
& =\varepsilon_{0}
{\displaystyle\sum_{n=0}^{\infty}}
\frac{\widehat{\chi}^{(n)}(0)}{n!}\left(  \int_{-\infty}^{\infty}d\omega
\omega^{n}\widehat{\mathbf{E}}(\mathbf{x},\omega)e^{-i\omega t}\right),
\nonumber\\
& =\varepsilon_{0}
{\displaystyle\sum_{n=0}^{\infty}}
\frac{\widehat{\chi}^{(n)}(0)}{n!}\left(  \int_{-\infty}^{\infty}
d\omega(i\partial_{t})^{n}\widehat{\mathbf{E}}(\mathbf{x},\omega)e^{-i\omega
t}\right), \nonumber\\
& =\varepsilon_{0}
{\displaystyle\sum_{n=0}^{\infty}}
\frac{\widehat{\chi}^{(n)}(0)}{n!}(i\partial_{t})^{n}\left(  \int_{-\infty
}^{\infty}d\omega\widehat{\mathbf{E}}(\mathbf{x},\omega)e^{-i\omega t}\right),
\nonumber\\
& =\widehat{\chi}(i\partial_{t})\mathbf{E}(\mathbf{x},t),\nonumber
\end{align}
where $\widehat{\chi}(\omega)$ is the fourier transform of $\chi(t)$. These
manipulations are of course purely formal; in order to make them into honest
mathematics we must dive into the theory of \textit{pseudo differential
operators}. In these lecture notes we will not do this as our focus is on
mathematical methods rather than mathematical theory.

Inserting (\ref{Eq8.2}),(\ref{Eq8.3}),(\ref{Eq8.4}) and (\ref{Eq8.5}) into
(\ref{Eq8.1}), we get Maxwell's equations in the form
\begin{align}
\partial_{t}\mathbf{B+\nabla\times E}  & =0,\nonumber\\
\partial_{t}\mathbf{E-c}^{2}\mathbf{\nabla\times B+\partial}_{t}\widehat{\chi
}(i\partial_{t})\mathbf{E}  & =-c^{2}\mu_{0}\partial_{t}\mathbf{P}
_{NL},\nonumber\\
\nabla\cdot\left(  \mathbf{E}+\widehat{\chi}(i\partial_{t})\mathbf{E}\right)
& =-\frac{1}{\varepsilon_{0}}\nabla\cdot\mathbf{P}_{NL},\nonumber\\
\nabla\cdot\mathbf{B}  & =0.\nonumber
\end{align}

\subsection*{TE scalar wave packets}

Let us first simplify the problem by only considering solutions of the form
\begin{align}
\mathbf{E}(x,y,z,t)  & =E(x,z,t)\mathbf{e}_{y},\nonumber\\
\mathbf{B}(x,y,z,t)  & =B_{1}(x,z,t)\mathbf{e}_{x}\mathbf{+}B_{2}
(x,z,t)\mathbf{e}_{z}.\label{Eq8.8}
\end{align}
For this simplified case, Maxwell's equations takes the form
\begin{align}
\partial_{t}B_{1}-\partial_{z}E  & =0,\nonumber\\
\partial_{t}B_{2}+\partial_{x}E  & =0,\nonumber\\
\partial_{t}E-c^{2}(\partial_{z}B_{1}-\partial_{x}B_{2})+\partial
_{t}\widehat{\chi}(i\partial_{t})E  & =-\partial_{t}P_{NL},\nonumber\\
\partial_{x}B_{1}+\partial_{z}B_{2}  & =0,\label{Eq8.9}
\end{align}
where
\begin{equation}
P_{NL}=\eta E^{3}.\label{Eq8.10}
\end{equation}
It is well known that this vector system is fully equivalent to the following
scalar equation
\begin{equation}
\partial_{tt}E-c^{2}\nabla^{2}E+\partial_{tt}\widehat{\chi}(i\partial
_{t})E=-\partial_{tt}P_{NL},\label{Eq8.11}
\end{equation}
where we have introduced the operator
\begin{equation}
\nabla^{2}=\partial_{xx}+\partial_{zz}.\label{Eq8.12}
\end{equation}
Equation (\ref{Eq8.11}) will be the staring point for our multiple scale
approach, but before that I will introduce the notion of a \textit{formal}
perturbation parameter. For some particular application of equation
(\ref{Eq8.11}) we will usually start by making the equation dimension-less by
picking some scales for space, time, and $E$ relevant for our particular
application. Here we don't want to tie our calculations to some particular
choise of scales and introduce therefore a formal perturbation parameter in
the equation multiplying the nonlinear polarization term. Thus we have
\begin{equation}
\partial_{tt}E-c^{2}\nabla^{2}E+\partial_{tt}\widehat{\chi}(i\partial
_{t})E=-\varepsilon^{2}\eta\partial_{tt}E^{3}.\label{Eq8.13}
\end{equation}
Starting with this equation we will proceed with our perturbation calculations
assuming that $\varepsilon<<1$ and in the end we will remove $\varepsilon$ by
setting it equal to $1$. What is going on here is that $\varepsilon$ is a
"place holder" for the actual small parameter that will appear in front of the
nonlinear term in the equation when we make a particular choise of scales.
Using such formal perturbation parameters is very common.

You might ask why I use $\varepsilon^{2}$ instead of $\varepsilon$ as formal
perturbation parameter? I will not answer this question here but will say
something about it at the very end of the lecture notes. We proceed with the
multiple scale method by introducing the expansions
\begin{align}
\partial_{t}  & =\partial_{t_{0}}+\varepsilon\partial_{t_{1}}+\varepsilon
^{2}\partial_{t_{2}}+...\;\;,\nonumber\\
\nabla & =\nabla_{0}+\varepsilon\nabla_{1}+\varepsilon^{2}\nabla
_{2}+...\;\;,\nonumber\\
e  & =e_{0}+\varepsilon e_{1}+\varepsilon^{2}e_{2}+...\;\;,\nonumber\\
E(\mathbf{x},t)  & =e(\mathbf{x}_{0},t_{0},\mathbf{x}_{1},t_{1},...)|_{t_{j}
=\varepsilon^{j}t,\mathbf{x}_{j}=\varepsilon^{j}\mathbf{x}},\label{Eq8.14}
\end{align}
where
\begin{equation}
\nabla_{j}=(\partial_{x_{j}},\partial_{z_{j}}),\label{Eq8.15}
\end{equation}
is the gradient with respect to $\mathbf{x}_{j}=(x_{j},z_{j})$. We now insert
(\ref{Eq8.14}) into (\ref{Eq8.13}) and expand everything in sight
\begin{align*}
& (\partial_{t_{0}}+\varepsilon\partial_{t_{1}}+\varepsilon^{2}\partial
_{t_{2}}+...)(\partial_{t_{0}}+\varepsilon\partial_{t_{1}}+\varepsilon
^{2}\partial_{t_{2}}+...)\\
& (e_{0}+\varepsilon e_{1}+\varepsilon^{2}e_{2}+...)-\\
& c^{2}(\nabla_{0}+\varepsilon\nabla_{1}+\varepsilon^{2}\nabla_{2}
+...)\cdot(\nabla_{0}+\varepsilon\nabla_{1}+\varepsilon^{2}\nabla_{2}+...)\\
& (e_{0}+\varepsilon e_{1}+\varepsilon^{2}e_{2}+...)+\\
& (\partial_{t_{0}}+\varepsilon\partial_{t_{1}}+\varepsilon^{2}\partial
_{t_{2}}+...)(\partial_{t_{0}}+\varepsilon\partial_{t_{1}}+\varepsilon
^{2}\partial_{t_{2}}+...)\\
& \widehat{\chi}(i\partial_{t_{0}}+i\varepsilon\partial_{t_{1}}+i\varepsilon
^{2}\partial_{t_{2}}+...)(e_{0}+\varepsilon e_{1}+\varepsilon^{2}e_{2}+...)\\
& =-\varepsilon^{2}\eta(\partial_{t_{0}}+\varepsilon\partial_{t_{1}
}+\varepsilon^{2}\partial_{t_{2}}+...)(\partial_{t_{0}}+\varepsilon
\partial_{t_{1}}+\varepsilon^{2}\partial_{t_{2}}+...)\\
& (e_{0}+\varepsilon e_{1}+\varepsilon^{2}e_{2}+...)^{3}\\
& \Downarrow\\
& (\partial_{t_{0}t_{0}}+\varepsilon(\partial_{t_{0}t_{1}}+\partial
_{t_{1}t_{0}})+\varepsilon^{2}(\partial_{t_{0}t_{2}}+\partial_{t_{1}t_{1}
}+\partial_{t_{2}t_{0}})+...)\\
& (e_{0}+\varepsilon e_{1}+\varepsilon^{2}e_{2}+...)-\\
& c^{2}(\nabla_{0}^{2}+\varepsilon(\nabla_{1}\cdot\nabla_{0}+\nabla_{0}
\cdot\nabla_{1})+\varepsilon^{2}(\nabla_{2}\cdot\nabla_{0}+\nabla_{1}
\cdot\nabla_{1}+\nabla_{0}\cdot\nabla_{2})+...)\\
& (e_{0}+\varepsilon e_{1}+\varepsilon^{2}e_{2}+...)+\\
& (\partial_{t_{0}t_{0}}+\varepsilon(\partial_{t_{0}t_{1}}+\partial
_{t_{1}t_{0}})+\varepsilon^{2}(\partial_{t_{0}t_{2}}+\partial_{t_{1}t_{1}
}+\partial_{t_{2}t_{0}})+...)\\
& (\widehat{\chi}(i\partial_{t_{0}})+\varepsilon\widehat{\chi}^{\prime
}(i\partial_{t_{0}})i\partial_{t_{1}}+\varepsilon^{2}(\widehat{\chi}^{\prime
}(i\partial_{t_{0}})i\partial_{t_{2}}-\frac{1}{2}\widehat{\chi}^{\prime\prime
}(i\partial_{t_{0}})\partial_{t_{1}t_{1}})+...)\\
& (e_{0}+\varepsilon e_{1}+\varepsilon^{2}e_{2}+...)\\
& =-\varepsilon^{2}\partial_{t_{0}t_{0}}e_{0}^{3}+...\;\;,\\
\nonumber\\
& \Downarrow
\end{align*}

\begin{align*}
& \partial_{t_{0}t_{0}}e_{0}+\varepsilon(\partial_{t_{0}t_{0}}e_{1}
+\partial_{t_{0}t_{1}}e_{0}+\partial_{t_{1}t_{0}}e_{0})\\
& +\varepsilon^{2}(\partial_{t_{0}t_{0}}e_{2}+\partial_{t_{0}t_{1}}
e_{1}+\partial_{t_{1}t_{0}}e_{1}+\partial_{t_{0}t_{2}}e_{0}+\partial
_{t_{1}t_{1}}e_{0}+\partial_{t_{2}t_{0}}e_{0})+...\\
& -c^{2}\nabla_{0}^{2}e_{0}-\varepsilon c^{2}(\nabla_{0}^{2}e_{1}+\nabla
_{1}\cdot\nabla_{0}e_{0}+\nabla_{0}\cdot\nabla_{1}e_{0})\\
& -\varepsilon^{2}c^{2}(\nabla_{0}^{2}e_{2}+\nabla_{1}\cdot\nabla_{0}
e_{1}+\nabla_{0}\cdot\nabla_{1}e_{1}\\
& +\nabla_{2}\cdot\nabla_{0}e_{0}+\nabla_{1}\cdot\nabla_{1}e_{0}+\nabla
_{0}\cdot\nabla_{2}e_{0})+...\\
& +\partial_{t_{0}t_{0}}\widehat{\chi}(i\partial_{t_{0}})e_{0}+\varepsilon
(\partial_{t_{0}t_{0}}\widehat{\chi}(i\partial_{t_{0}})e_{1}+\partial
_{t_{0}t_{0}}\widehat{\chi}^{\prime}(i\partial_{t_{0}})i\partial_{t_{1}}
e_{0}\\
& +\partial_{t_{0}t_{1}}\widehat{\chi}(i\partial_{t_{0}})e_{0}+\partial
_{t_{1}t_{0}}\widehat{\chi}(i\partial_{t_{0}})e_{0})+\varepsilon^{2}
(\partial_{t_{0}t_{0}}\widehat{\chi}(i\partial_{t_{0}})e_{2}\\
& +\partial_{t_{0}t_{0}}\widehat{\chi}^{\prime}(i\partial_{t_{0}}
)i\partial_{t_{1}}e_{1}+\partial_{t_{0}t_{1}}\widehat{\chi}(i\partial_{t_{0}
})e_{1}+\partial_{t_{1}t_{0}}\widehat{\chi}(i\partial_{t_{0}})e_{1}\\
& +\partial_{t_{0}t_{0}}\widehat{\chi}^{\prime}(i\partial_{t_{0}}
)i\partial_{t_{2}}e_{0}-\frac{1}{2}\partial_{t_{0}t_{0}}\widehat{\chi}
^{\prime\prime}(i\partial_{t_{0}})\partial_{t_{1}t_{1}}e_{0}+\partial
_{t_{1}t_{0}}\widehat{\chi}^{\prime}(i\partial_{t_{0}})i\partial_{t_{1}}
e_{0}\\
& +\partial_{t_{0}t_{1}}\widehat{\chi}^{\prime}(i\partial_{t_{0}}
)i\partial_{t_{1}}e_{0}+\partial_{t_{2}t_{0}}\widehat{\chi}(i\partial_{t_{0}
})e_{0}+\partial_{t_{1}t_{1}}\widehat{\chi}(i\partial_{t_{0}})e_{0}\\
& +\partial_{t_{0}t_{2}}\widehat{\chi}(i\partial_{t_{0}})e_{0})+...\\
& =-\varepsilon^{2}\partial_{t_{0}t_{0}}e_{0}^{3}+...\;\;,
\end{align*}
which gives us the perturbation hierarchy
\begin{gather}
\partial_{t_{0}t_{0}}e_{0}-c^{2}\nabla_{0}^{2}e_{0}+\partial_{t_{0}t_{0}
}\widehat{\chi}(i\partial_{t_{0}})e_{0}=0,\label{Eq8.16}\\
\nonumber\\
\partial_{t_{0}t_{0}}e_{1}-c^{2}\nabla_{0}^{2}e_{1}+\partial_{t_{0}t_{0}
}\widehat{\chi}(i\partial_{t_{0}})e_{1}=\nonumber\\
-\partial_{t_{0}t_{1}}e_{0}-\partial_{t_{1}t_{0}}e_{0}-c^{2}\nabla_{1}
\cdot\nabla_{0}e_{0}-c^{2}\nabla_{0}\cdot\nabla_{1}e_{0}\nonumber\\
-\partial_{t_{0}t_{0}}\widehat{\chi}^{\prime}(i\partial_{t_{0}})i\partial
_{t_{1}}e_{0}-\partial_{t_{0}t_{1}}\widehat{\chi}(i\partial_{t_{0}}
)e_{0}-\partial_{t_{1}t_{0}}\widehat{\chi}(i\partial_{t_{0}})e_{0},\label{Eq8.17}\\
\nonumber\\
\partial_{t_{0}t_{0}}e_{2}-c^{2}\nabla_{0}^{2}e_{2}+\partial_{t_{0}t_{0}
}\widehat{\chi}(i\partial_{t_{0}})e_{2}=\nonumber\\
-\partial_{t_{0}t_{1}}e_{1}-\partial_{t_{1}t_{0}}e_{1}-\partial_{t_{0}t_{2}
}e_{0}-\partial_{t_{1}t_{1}}e_{0}-\partial_{t_{2}t_{0}}e_{0}\nonumber\\
-c^{2}\nabla_{1}\cdot\nabla_{0}e_{1}-c^{2}\nabla_{0}\cdot\nabla_{1}e_{1}
-c^{2}\nabla_{2}\cdot\nabla_{0}e_{0}-c^{2}\nabla_{1}\cdot\nabla_{1}
e_{0}\nonumber\\
-c^{2}\nabla_{0}\cdot\nabla_{2}e_{0}-\partial_{t_{0}t_{0}}\widehat{\chi
}^{\prime}(i\partial_{t_{0}})i\partial_{t_{1}}e_{1}-\partial_{t_{0}t_{1}
}\widehat{\chi}(i\partial_{t_{0}})e_{1}\nonumber\\
-\partial_{t_{1}t_{0}}\widehat{\chi}(i\partial_{t_{0}})e_{1}-\partial
_{t_{0}t_{0}}\widehat{\chi}^{\prime}(i\partial_{t_{0}})i\partial_{t_{2}}
e_{0}+\frac{1}{2}\partial_{t_{0}t_{0}}\widehat{\chi}^{\prime\prime}
(i\partial_{t_{0}})\partial_{t_{1}t_{1}}e_{0}\nonumber\\
-\partial_{t_{1}t_{0}}\widehat{\chi}^{\prime}(i\partial_{t_{0}})i\partial
_{t_{1}}e_{0}-\partial_{t_{0}t_{1}}\widehat{\chi}^{\prime}(i\partial_{t_{0}
})i\partial_{t_{1}}e_{0}-\partial_{t_{2}t_{0}}\widehat{\chi}(i\partial_{t_{0}
})e_{0}\nonumber\\
-\partial_{t_{1}t_{1}}\widehat{\chi}(i\partial_{t_{0}})e_{0}-\partial
_{t_{0}t_{2}}\widehat{\chi}(i\partial_{t_{0}})e_{0}-\partial_{t_{0}t_{0}}
e_{0}^{3}.\label{Eq8.18}
\end{gather}
For the order $\varepsilon^{0}$ equation we choose the wave packet solution
\begin{equation}
e_{0}(\mathbf{x}_{0},t_{0},\mathbf{x}_{1},t_{1},..)=A_{0}(\mathbf{x}_{1}
,t_{1},...)e^{i\theta_{0}}+(\ast),\label{Eq8.19}
\end{equation}
where
\begin{align}
\mathbf{x}_{j}  & =(x_{j},z_{j}),\nonumber\\
\theta_{0}  & =\mathbf{k}\cdot\mathbf{x}_{0}-\omega t_{0},\label{Eq8.20}
\end{align}
and where $\mathbf{k}$ is a plane vector with components $\mathbf{k}=(\xi
,\eta)$. In (\ref{Eq8.20}), $\omega$, is a function of $k=||\mathbf{k}||$ that
satisfy the dispersion relation
\begin{equation}
\omega^{2}n^{2}(\omega)=c^{2}k^{2},\label{Eq8.21}
\end{equation}
where the refractive index, $n(\omega)$, is defined by
\begin{equation}
n^{2}(\omega)=1+\widehat{\chi}(\omega).\label{Eq8.22}
\end{equation}
We now must now calculate the right-hand side of the order $\varepsilon$ equation.

Observe that
\begin{align}
\partial_{t_{1}t_{0}}e_{0}  & =-i\omega\partial_{t_{1}}A_{0}e^{i\theta_{0}
}+(\ast),\nonumber\\
\partial_{t_{0}t_{1}}e_{0}  & =-i\omega\partial_{t_{1}}A_{0}e^{i\theta_{0}
}+(\ast),\nonumber\\
\nabla_{1}\cdot\nabla_{0}e_{0}  & =ik\nabla_{1}A_{0}\cdot\mathbf{u}
e^{i\theta_{0}}+(\ast),\nonumber\\
\nabla_{0}\cdot\nabla_{1}e_{0}  & =ik\nabla_{1}A_{0}\cdot\mathbf{u}
e^{i\theta_{0}}+(\ast),\nonumber\\
\partial_{t_{0}t_{0}}\widehat{\chi}^{\prime}(i\partial_{t_{0}})i\partial
_{t_{1}}e_{0}  & =-i\omega\widehat{\chi}^{\prime}(\omega)\partial_{t_{1}}
A_{0}e^{i\theta_{0}}+(\ast),\nonumber\\
\partial_{t_{0}t_{1}}\widehat{\chi}(i\partial_{t_{0}})e_{0}  & =-i\omega
\widehat{\chi}(\omega)\partial_{t_{1}}A_{0}e^{i\theta_{0}}+(\ast),\nonumber\\
\partial_{t_{1}t_{0}}\widehat{\chi}(i\partial_{t_{0}})e_{0}  & =-i\omega
\widehat{\chi}(\omega)\partial_{t_{1}}A_{0}e^{i\theta_{0}}+(\ast),\nonumber
\end{align}
where $\mathbf{u}$ is a unit vector in the direction of $\mathbf{k}$.
Inserting (\ref{Eq8.23}) into (\ref{Eq8.17}) we get
\begin{gather}
\partial_{t_{0}t_{0}}e_{1}-c^{2}\nabla_{0}^{2}e_{1}+\partial_{t_{0}t_{0}
}\widehat{\chi}(i\partial_{t_{0}})e_{1}=\nonumber\\
-\{-2i\omega\partial_{t_{1}}A_{0}-2ic^{2}k\mathbf{u}\cdot\nabla_{1}
A_{0}\nonumber\\
-i\omega^{2}\widehat{\chi}^{\prime}(\omega)\partial_{t_{1}}A_{0}
-2i\omega\widehat{\chi}(\omega)\partial_{t_{1}}A_{0}\}e^{i\theta_{0}}
+(\ast).\label{Eq8.24}
\end{gather}
In order to remove secular terms we must postulate that
\begin{gather}
-2i\omega\partial_{t_{1}}A_{0}-2ic^{2}k\mathbf{u}\cdot\nabla_{1}A_{0}
-i\omega^{2}\widehat{\chi}^{\prime}(\omega)\partial_{t_{1}}A_{0}
-2i\omega\widehat{\chi}(\omega)\partial_{t_{1}}A_{0}=0,\nonumber\\
\Updownarrow\nonumber\\
\omega(2n^{2}+\omega\widehat{\chi}^{\prime}(\omega))\partial_{t_{1}}
A_{0}-2ic^{2}k\mathbf{u}\cdot\nabla_{1}A_{0}=0.\label{Eq8.25}
\end{gather}
Observe that from the dispersion relation (\ref{Eq8.21}) we have
\begin{align*}
\omega^{2}n^{2}(\omega)  & =c^{2}k^{2},\\
& \Updownarrow\\
\omega^{2}(1+\widehat{\chi}(\omega))  & =c^{2}k^{2},\\
& \Downarrow\\
2\omega\omega^{\prime}n^{2}(\omega)+\omega^{2}\widehat{\chi}^{\prime}
(\omega)\omega^{\prime}  & =2c^{2}k,\\
& \Downarrow\\
\omega(2n^{2}+\omega\widehat{\chi}^{\prime}(\omega))\omega^{\prime}  &
=2c^{2}k.
\end{align*}
Thus (\ref{Eq8.25}) can be written in the form
\begin{equation}
\partial_{t_{1}}A_{0}+\mathbf{v}_{g}\cdot\nabla_{1}A_{0}=0,\label{Eq8.26}
\end{equation}
where $\mathbf{v}_{g}$ is the group velocity
\begin{equation}
\mathbf{v}_{g}=\omega^{\prime}(k)\mathbf{u}.\label{Eq8.27}
\end{equation}
The order $\varepsilon$ equation simplifies into
\begin{equation}
\partial_{t_{0}t_{0}}e_{1}-c^{2}\nabla_{0}^{2}e_{1}+\partial_{t_{0}t_{0}
}\widehat{\chi}(i\partial_{t_{0}})e_{1}=0.\label{Eq8.29}
\end{equation}
According to the rules of the game we choose the special solution
\begin{equation}
e_{1}=0,\label{Eq8.30}
\end{equation}
for (\ref{Eq8.29}). We now must compute the right-hand side of the order
$\varepsilon^{2}$ equation. Observe that

\begin{align}
\partial_{t_{2}t_{0}}e_{0}  & =-i\omega\partial_{t_{2}}A_{0}e^{i\theta_{0}
}+(\ast),\nonumber\\
\partial_{t_{1}t_{1}}e_{0}  & =\partial_{t_{1}t_{1}}A_{0}e^{i\theta_{0}}
+(\ast),\nonumber\\
\partial_{t_{0}t_{1}}e_{0}  & =-i\omega\partial_{t_{2}}A_{0}e^{i\theta_{0}
}+(\ast),\nonumber\\
\nabla_{2}\cdot\nabla_{0}e_{0}  & =ik\mathbf{u}\cdot\nabla_{2}A_{0}
e^{i\theta_{0}}+(\ast),\nonumber\\
\nabla_{1}\cdot\nabla_{1}e_{0}  & =\nabla_{1}^{2}A_{0}e^{i\theta_{0}}
+(\ast),\nonumber\\
\nabla_{0}\cdot\nabla_{2}e_{0}  & =ik\mathbf{u}\cdot\nabla_{2}A_{0}
e^{i\theta_{0}}+(\ast),\nonumber\\
\partial_{t_{0}t_{0}}\widehat{\chi}^{\prime}(i\partial_{t_{0}})i\partial
_{t_{2}}e_{0}  & =-i\omega^{2}\widehat{\chi}^{\prime}(\omega)\partial_{t_{2}
}A_{0}e^{i\theta_{0}}+(\ast),\nonumber\\
\frac{1}{2}\partial_{t_{0}t_{0}}\widehat{\chi}^{\prime\prime}(i\partial
_{t_{0}})\partial_{t_{1}t_{1}}e_{0}  & =-\frac{1}{2}\omega^{2}\widehat{\chi
}^{\prime\prime}(\omega)\partial_{t_{1}t_{1}}A_{0}e^{i\theta_{0}}
+(\ast),\nonumber\\
\partial_{t_{1}t_{0}}\widehat{\chi}^{\prime}(i\partial_{t_{0}})i\partial
_{t_{1}}e_{0}  & =\omega\widehat{\chi}^{\prime}(\omega)\partial_{t_{1}t_{1}
}A_{0}e^{i\theta_{0}}+(\ast),\nonumber\\
\partial_{t_{0}t_{1}}\widehat{\chi}^{\prime}(i\partial_{t_{0}})i\partial
_{t_{1}}e_{0}  & =\omega\widehat{\chi}^{\prime}(\omega)\partial_{t_{1}t_{1}
}A_{0}e^{i\theta_{0}}+(\ast),\nonumber\\
\partial_{t_{2}t_{0}}\widehat{\chi}(i\partial_{t_{0}})e_{0}  & =-i\omega
\widehat{\chi}(\omega)\partial_{t_{2}}A_{0}e^{i\theta_{0}}+(\ast),\nonumber\\
\partial_{t_{1}t_{1}}\widehat{\chi}(i\partial_{t_{0}})e_{0}  & =\widehat{\chi
}(\omega)\partial_{t_{1}t_{1}}A_{0}e^{i\theta_{0}}+(\ast),\nonumber\\
\partial_{t_{0}t_{2}}\widehat{\chi}(i\partial_{t_{0}})e_{0}  & =-i\omega
\widehat{\chi}(\omega)\partial_{t_{2}}A_{0}e^{i\theta_{0}}+(\ast),\nonumber\\
\partial_{t_{0}t_{0}}e_{0}^{3}  & =-3\omega^{2}\eta|A_{0}|^{2}A_{0}e^{i\theta
}+NST+(\ast).\label{Eq8.31}
\end{align}
Inserting (\ref{Eq8.30}) and (\ref{Eq8.31}) into the right-hand side of the
order $\varepsilon^{2}$ equation we get
\begin{gather}
\partial_{t_{0}t_{0}}e_{2}-c^{2}\nabla_{0}^{2}e_{2}+\partial_{t_{0}t_{0}
}\widehat{\chi}(i\partial_{t_{0}})e_{2}=\nonumber\\
-\{-2i\omega\partial_{t_{2}}A_{0}+\partial_{t_{1}t_{1}}A_{0}-2ic^{2}
k\mathbf{u}\cdot\nabla_{2}A_{0}-c^{2}\nabla_{1}^{2}A_{0}\nonumber\\
-i\omega^{2}\widehat{\chi}^{\prime}(\omega)\partial_{t_{2}}A_{0}+\frac{1}
{2}\omega^{2}\widehat{\chi}^{\prime\prime}(\omega)\partial_{t_{1}t_{1}}
A_{0}+2\omega\widehat{\chi}^{\prime}(\omega)\partial_{t_{1}t_{1}}
A_{0}\nonumber\\
-2i\omega\widehat{\chi}(\omega)\partial_{t_{2}}A_{0}+\widehat{\chi}
(\omega)\partial_{t_{1}t_{1}}A_{0}-3\omega^{2}\eta|A_{0}|^{2}\}e^{i\theta_{0}
}+NST+(\ast).\label{Eq8.32}
\end{gather}
In order to remove secular terms we must postulate that
\begin{gather}
-2i\omega\partial_{t_{2}}A_{0}+\partial_{t_{1}t_{1}}A_{0}-2ic^{2}
k\mathbf{u}\cdot\nabla_{2}A_{0}-c^{2}\nabla_{1}^{2}A_{0}-i\omega
^{2}\widehat{\chi}^{\prime}(\omega)\partial_{t_{2}}A_{0}\nonumber\\
+\frac{1}{2}\omega^{2}\widehat{\chi}^{\prime\prime}(\omega)\partial
_{t_{1}t_{1}}A_{0}+2\omega\widehat{\chi}^{\prime}(\omega)\partial_{t_{1}t_{1}
}A_{0}-2i\omega\widehat{\chi}(\omega)\partial_{t_{2}}A_{0}+\widehat{\chi
}(\omega)\partial_{t_{1}t_{1}}A_{0}\nonumber\\
-3\omega^{2}\eta|A_{0}|^{2}=0.\label{Eq8.33}
\end{gather}
Using the dispersion relation (\ref{Eq8.22}), equation (\ref{Eq8.33}) can be
simplified into
\begin{equation}
\partial_{t_{2}}A_{0}+\mathbf{v}_{g}\cdot\nabla_{2}A_{0}-i\beta\nabla_{1}
^{2}A_{0}+i\alpha\partial_{t_{1}t_{1}}A_{0}-i\gamma|A_{0}|^{2}A_{0}
=0,\label{Eq8.34}
\end{equation}
where
\begin{align}
\alpha & =\omega^{\prime}\frac{n^{2}+2\omega\widehat{\chi}^{\prime}
(\omega)+\frac{1}{2}\omega^{2}\widehat{\chi}^{\prime\prime}(\omega)}{2c^{2}
k},\nonumber\\
\beta & =\frac{\omega^{\prime}}{2k},\nonumber\\
\gamma & =\frac{3\eta\omega^{2}\omega^{\prime}}{2c^{2}k}.\nonumber
\end{align}
Defining an amplitude $A(\mathbf{x},t)$ by
\begin{equation}
A(\mathbf{x},t)=A_{0}(\mathbf{x}_{1},t_{1},...)|_{t_{j}=e^{j}t,\mathbf{x}
_{j}=\varepsilon^{j}\mathbf{x}},\label{Eq8.36}
\end{equation}
and proceeding in the usual way, using (\ref{Eq8.26}) and (\ref{Eq8.34}), we
get the following amplitude equation
\begin{equation}
\partial_{t}A+\mathbf{v}_{g}\cdot\nabla A-i\beta\nabla^{2}A+i\alpha
\partial_{tt}A-i\gamma|A|^{2}A=0,\label{Eq8.37}
\end{equation}
where we have put the formal perturbation parameter equal to $1$. From what we
have done it is evident that for
\begin{equation}
E(\mathbf{x},t)=A(\mathbf{x},t)e^{i(\mathbf{k}\cdot\mathbf{x}-\omega t)}
+(\ast),\label{Eq8.38}
\end{equation}
to be an approximate solution to (\ref{Eq8.13}) we must have
\begin{align}
\gamma|A|^{2}  & \sim\beta\nabla^{2}A\sim\alpha\partial_{tt}A\sim
O(\varepsilon^{2}),\nonumber\\
\partial_{t}A  & \sim\mathbf{v}_{g}\cdot\nabla A\sim O(\varepsilon),\label{Eq8.39}
\end{align}
where $\varepsilon$ is a number much smaller than $1$. Under these
circumstances (\ref{Eq8.37}),(\ref{Eq8.38}) is the key elements in a fast
numerical scheme for wave packet solutions to (\ref{Eq8.13}). Because of the
presence of the second derivative with respect to time, equation
(\ref{Eq8.37}) can not be solved as a standard initial value problem. However,
because of (\ref{Eq8.39}) we can remove the second derivative term by
iteration
\begin{align}
\partial_{t}A  & =-\mathbf{v}_{g}\cdot\nabla A\sim O(\varepsilon
),\nonumber\\
& \Downarrow\nonumber\\
\partial_{tt}A  & =(\mathbf{v}_{g}\cdot\nabla)^{2}A\sim O(\varepsilon
^{2}),\label{Eq8.40}
\end{align}
which leads to the equation
\begin{equation}
\partial_{t}A+\mathbf{v}_{g}\cdot\nabla A-i\beta\nabla^{2}A+i\alpha
(\mathbf{v}_{g}\cdot\nabla)^{2}A-i\gamma|A|^{2}A=0,\label{Eq8.41}
\end{equation}
which \textit{can} be solved as a standard initial value problem.

  In deriving
this equation we asssumed that the terms proportional to
\[
e^{\pm3i(\mathbf{k}\cdot\mathbf{x}-\omega t)},
\]
where nonsecular. For this to be true we must have
\begin{equation}
\omega(3k)\neq3\omega(k),\label{Eq8.42}
\end{equation}
where $\omega(k)$ is a solution to (\ref{Eq8.21}). If an equality holds in
(\ref{Eq8.42}) we have \textit{phase matching} and the multiple scale
calculation has to be redone, starting at (\ref{Eq8.19}), using a sum of two
wave packets with the appropriate center wave numbers and frequencies instead
of the single wavepacket we used in the calculation leading to (\ref{Eq8.37}).
It could also be the case that we are modelling a situation where several wave
packets are interacting in a Kerr medium. For such a case we would instead of
(\ref{Eq8.19}) use a finite sum of wave packets
\begin{equation}
e_{0}(\mathbf{x}_{0},t_{0},\mathbf{x}_{1},t_{1},..)=
{\displaystyle\sum_{j=0}^{N}}
A_{j}(\mathbf{x}_{1},t_{1},...)e^{i\theta_{j}}+(\ast).\label{Eq8.43}
\end{equation}
Calculations analogous to the ones leading up to equation (\ref{Eq8.37})
\ will now give a separate equation of the type (\ref{Eq8.37}) for each wave
packet, \textit{unless} we have phase matching. These phase matching
conditions appears from the nonlinear term in the order $\varepsilon^{2}$
equation and takes the familiar form
\begin{align}
\mathbf{k}_{j}  & =s_{1}\mathbf{k}_{j_{1}}+s_{2}\mathbf{k}_{j_{2}}
+s_{3}\mathbf{k}_{j_{3}},\nonumber\\
\omega(k_{j})  & =s_{1}\omega(k_{j_{1}})+s_{2}\omega(k_{j_{2}})+s_{3}
\omega(k_{j_{3}}),\label{Eq8.44}
\end{align}
where $s=\pm1$. The existence of phase matching leads to coupling of the
amplitude equations. If (\ref{Eq8.44}) holds, the amplitude equation for
$A_{j}$ will contain a coupling term proportional to
\begin{equation}
A_{j_{1}}^{s_{1}}A_{j_{2}}^{s_{2}}A_{j_{3}}^{s_{3}}\label{Eq8.46}%
\end{equation}
where by definition $A_{j}^{+1}=A_{j}$ and $A_{j}^{-1}=A_{j}^{\ast}$.

We have seeen that assuming a scaling of $\varepsilon$ for space and time
variables and $\varepsilon^{2}$ for the nonlinear term leads to an amplitude
equation where second derivatives of space and time appears at the same order
as the cubic nonlinearity. This amplitude equation can thus describe a
situation where diffraction, group velocity dispersion and nonlinearity are of
the same size. Other choises of scaling for space,time and nonlinearity will
lead to other amplitude equations where other physical effects are of the same
size. Thus, the choise of scaling is determined by what kind of physics we want to describe.

\subsection*{Linearly polarized vector wave packets}

Up til now all applications of the multiple scale method for PDEs  has involved scalar
equations. The multiple scale method is not limited to scalar equations and is
equally applicable to vector equations. However, for vector equations we need
to be more careful than for the scalar case when it comes to elliminating
secular terms. We will here use Maxwell's equations (\ref{Eq8.7}) to
illustrate how the method is applied to vector PDEs in general. Assuming,
as usual, a polarization response induced by the Kerr effect, our basic
equations are
\begin{align}
\partial_{t}\mathbf{B+\nabla\times E}  & =0,\nonumber\\
\partial_{t}\mathbf{E-c}^{2}\mathbf{\nabla\times B+\partial}_{t}\widehat{\chi
}(i\partial_{t})\mathbf{E}  & =-\varepsilon^{2}\eta\partial_{t}(\mathbf{E}
^{2}\mathbf{E)},\nonumber\\
\nabla\cdot\mathbf{B}  & =0,\nonumber\\
\nabla\cdot\mathbf{E}+\widehat{\chi}(i\partial_{t})\nabla\cdot\mathbf{E}  &
=-\varepsilon^{2}\eta\nabla\cdot(\mathbf{E}^{2}\mathbf{E)},\label{Eq8.47}
\end{align}
where we have introduced a formal perturbation parameter in front of the
nonlinear terms. We now introduce the usual machinery of the multiple scale method.

Let $\mathbf{e}(\mathbf{x}_{0},t_{0},\mathbf{x}_{1},t_{1},...)$ and
$\mathbf{b}(\mathbf{x}_{0},t_{0},\mathbf{x}_{1},t_{1},...)$ be functions such
that
\begin{align}
\mathbf{E}(x,t)  & =\mathbf{e}(\mathbf{x}_{0},t_{0},\mathbf{x}_{1}
,t_{1},...)|_{\mathbf{x}_{j}=\varepsilon^{j}\mathbf{x},t_{j}=\varepsilon^{j}
t},\nonumber\\
\mathbf{B}(x,t)  & =\mathbf{b}(\mathbf{x}_{0},t_{0},\mathbf{x}_{1}
,t_{1},...)|_{\mathbf{x}_{j}=\varepsilon^{j}\mathbf{x},t_{j}=\varepsilon^{j}
t},\label{Eq8.48}
\end{align}
and let
\begin{align}
\partial_{t}  & =\partial_{t_{0}}+\varepsilon\partial_{t_{1}}+\varepsilon
^{2}\partial_{t_{2}}+...\;\;,\nonumber\\
\nabla\times & =\nabla_{0}\times+\varepsilon\nabla_{1}\times+\varepsilon
^{2}\nabla_{2}\times+...\;\;,\nonumber\\
\nabla\cdot & =\nabla_{0}\cdot+\varepsilon\nabla_{1}\cdot+\varepsilon
^{2}\nabla_{2}\cdot+...\;\;,\nonumber\\
\mathbf{e}  & =\mathbf{e}_{0}+\varepsilon\mathbf{e}_{1}+\varepsilon
^{2}\mathbf{e}_{2}+...\;\;,\nonumber\\
\mathbf{b}  & =\mathbf{b}_{0}+\varepsilon\mathbf{b}_{1}+\varepsilon
^{2}\mathbf{b}_{2}+...\;\;.\label{Eq8.49}
\end{align}
We now insert (\ref{Eq8.49}) into (\ref{Eq8.47}) and expand everything in
sight to second order in $\varepsilon$. Putting each order of $\varepsilon$ to
zero separately gives us the perturbation hierarchy. At this point you should
be able to do this on your own so I will just write down the elements of the
perturbation hierarchy when they are needed.

The order $\varepsilon^{0}$ equations, which is the first element of the
perturbation hierarchy, is of course
\begin{align}
\partial_{t_{0}}\mathbf{b}_{0}\mathbf{+\nabla}_{0}\mathbf{\times e}_{0}  &
=0,\nonumber\\
\partial_{t_{0}}\mathbf{e}_{0}\mathbf{-c}^{2}\mathbf{\nabla}_{0}\mathbf{\times
b}_{0}\mathbf{+\partial}_{t_{0}}\widehat{\chi}(i\partial_{t_{0}}
)\mathbf{e}_{0}  & =0,\nonumber\\
\nabla_{0}\cdot\mathbf{b}_{0}  & =0,\nonumber\\
\nabla_{0}\cdot\mathbf{e}_{0}+\widehat{\chi}(i\partial_{t_{0}})\nabla_{0}
\cdot\mathbf{e}_{0}  & =0.\label{Eq8.50}
\end{align}
For the order $\varepsilon^{0}$ equations, we chose a linearly polarized wave
packet solution. It must be of the form
\begin{align}
\mathbf{e}_{0}(\mathbf{x}_{0},t_{0},\mathbf{x}_{1},t_{1},...)  & =\omega
A_{0}(\mathbf{x}_{1},t_{1},...)\mathbf{qe}^{i\theta_{0}}+(\ast),\nonumber
\\
\mathbf{b}_{0}(\mathbf{x}_{0},t_{0},\mathbf{x}_{1},t_{1},...)  &
=kA_{0}(\mathbf{x}_{1},t_{1},...)\mathbf{te}^{i\theta_{0}}+(\ast),\label{Eq8.51}
\end{align}
where
\begin{equation}
\theta_{0}=\mathbf{k}\cdot\mathbf{x}_{0}-\omega t_{0},\label{Eq8.52}
\end{equation}
and where
\[
\omega=\omega(k),
\]
is a solution to the dispersion relation
\begin{equation}
\omega^{2}n^{2}(\omega)=c^{2}k^{2}.\label{Eq8.53}
\end{equation}
The orthogonal unit vectors $\mathbf{q}$ and $\mathbf{t}$ span the space
transverse to $\mathbf{k}=k\mathbf{u}$, and the unit vectors $\{\mathbf{q}
,\mathbf{t},\mathbf{u}\}$ define a postively oriented fram for $\mathbb{R}^{3}$.

The order $\varepsilon$ equations are
\begin{gather}
\partial_{t_{0}}\mathbf{b}_{1}\mathbf{+\nabla}_{0}\mathbf{\times e}
_{1}=-\partial_{t_{1}}\mathbf{b}_{0}-\nabla_{1}\times\mathbf{e}_{0},\nonumber\\
\partial_{t_{0}}\mathbf{e}_{1}\mathbf{-c}^{2}\mathbf{\nabla}_{0}\mathbf{\times
b}_{1}\mathbf{+\partial}_{t_{0}}\widehat{\chi}(i\partial_{t_{0}}
)\mathbf{e}_{1}=\nonumber\\
-\partial_{t_{1}}\mathbf{e}_{0}+c^{2}\nabla_{1}\times\mathbf{b}_{0}
-\partial_{t_{1}}\widehat{\chi}(i\partial_{t_{0}})\mathbf{e}_{0}
-i\partial_{t_{0}}\widehat{\chi}^{\prime}(i\partial_{t_{0}})\partial_{t_{1}
}\mathbf{e}_{0},\nonumber\\
\nabla_{0}\cdot\mathbf{b}_{1}=-\nabla_{1}\cdot b_{0},\nonumber\\
\nabla_{0}\cdot\mathbf{e}_{1}+\widehat{\chi}(i\partial_{t_{0}})\nabla_{0}
\cdot\mathbf{e}_{1}=\nonumber\\
-\nabla_{1}\cdot\mathbf{e}_{0}-\widehat{\chi}(i\partial_{t_{0}})\nabla
_{1}\cdot\mathbf{e}_{0}-i\widehat{\chi}^{\prime}(i\partial_{t_{0}}
)\partial_{t_{1}}\nabla_{0}\cdot\mathbf{e}_{0}.\label{Eq8.54}
\end{gather}
Inserting (\ref{Eq8.51}) into (\ref{Eq8.54}) we get
\begin{gather}
\partial_{t_{0}}\mathbf{b}_{0}\mathbf{+\nabla}_{0}\mathbf{\times e}
_{0}=-\{k\partial_{t_{1}}A_{0}\mathbf{t}+\omega\nabla_{1}A_{0}\times
\mathbf{q\}}e^{i\theta_{0}}+(\ast),\nonumber\\
\partial_{t_{0}}\mathbf{e}_{0}\mathbf{-c}^{2}\mathbf{\nabla}_{0}\mathbf{\times
b}_{0}\mathbf{+\partial}_{t_{0}}\widehat{\chi}(i\partial_{t_{0}}
)\mathbf{e}_{0}=-\{(\omega n^{2}(\omega)+\omega^{2}\widehat{\chi}^{\prime
}(\omega))\partial_{t_{1}}A_{0}\mathbf{q}\nonumber\\
-c^{2}k\nabla_{1}A_{0}\times\mathbf{t}\}e^{i\theta_{0}}+(\ast),\nonumber\\
\nabla_{0}\cdot\mathbf{b}_{0}=-\{k\nabla_{1}A_{0}\cdot\mathbf{t\}}
e^{i\theta_{0}}+(\ast),\nonumber\\
\nabla_{0}\cdot\mathbf{e}_{0}+\widehat{\chi}(i\partial_{t_{0}})\nabla_{0}
\cdot\mathbf{e}_{0}=-\{\omega n^{2}(\omega)\nabla_{1}A_{0}\cdot\mathbf{q\}}
e^{i\theta_{0}}+(\ast).\label{Eq8.55}
\end{gather}  
If we can find a special solution to this system that is bounded, we will get a
perturbation expansion that is uniform for $t\lesssim\varepsilon^{-1}$. We
will look for solutions of the form
\begin{align}
\mathbf{e}_{1}  & =\mathbf{a}e^{i\theta_{0}}+(\ast),\nonumber\\
\mathbf{b}_{1}  & =\mathbf{b}e^{i\theta_{0}}+(\ast),\label{Eq8.56}
\end{align}
where $\mathbf{a}$ and $\mathbf{b}$ are constant vectors. Inserting
(\ref{Eq8.56}) into (\ref{Eq8.55}), we get the following linear algebraic
system of equations for the unknown vectors $\mathbf{a}$ and $\mathbf{b}$
\begin{gather}
-i\omega\mathbf{b}+ik\mathbf{u}\times\mathbf{a}=-\{k\partial_{t_{1}}
A_{0}\mathbf{t}+\omega\nabla_{1}A_{0}\times\mathbf{q\}},\label{Eq8.57}\\
-i\omega n^{2}(\omega)\mathbf{a}-ic^{2}k\mathbf{u}\times\mathbf{b}=-\{(\omega
n^{2}(\omega)+\omega^{2}\widehat{\chi}^{\prime}(\omega))\partial_{t_{1}}
A_{0}\mathbf{q}\nonumber\\
-c^{2}k\nabla_{1}A_{0}\times\mathbf{t}\},\label{Eq8.58}\\
ik\mathbf{u}\cdot\mathbf{b}=-k\nabla_{1}A_{0}\cdot\mathbf{t},\label{Eq8.59}\\
ikn^{2}(\omega)\mathbf{u}\cdot\mathbf{a}=-\omega n^{2}(\omega)\nabla_{1}
A_{0}\cdot\mathbf{q}\;.\label{Eq8.60}
\end{gather}
Introduce the longitudinal and transverse parts of $\mathbf{a}$ and
$\mathbf{b}$ through
\begin{align}
a_{\Vert}  & =(\mathbf{u}\cdot\mathbf{a)u},\text{ \ \ \ \ \ \ \ }a_{\bot
}=\mathbf{a}-a_{\Vert},\nonumber\\
b_{\Vert}  & =(\mathbf{u}\cdot\mathbf{b})\mathbf{u},\text{ \ \ \ \ \ \ \ }
b_{\bot}=\mathbf{b}-b_{\Vert}.\label{Eq8.61}
\end{align}
Then from (\ref{Eq8.59}) and (\ref{Eq8.60}) we get
\begin{align}
a_{\Vert}  & =(i\frac{\omega}{k}\nabla_{1}A_{0}\cdot\mathbf{q)u},
\label{Eq8.62}\\
b_{\Vert}  & =(i\nabla_{1}A_{0}\cdot\mathbf{t})\mathbf{u}.\label{Eq8.63}
\end{align}
However, the longitudinal part of (\ref{Eq8.57}) and (\ref{Eq8.58}) will also
determine $a_{\Vert}$ and $b_{\Vert}$. These values must be the same as the
ones just found in (\ref{Eq8.62}),(\ref{Eq8.63}). These are
\textit{solvability conditions}. Taking the longitudinal part of
(\ref{Eq8.57}) we get
\begin{gather}
-i\omega\mathbf{u}\cdot\mathbf{b}=-\omega\mathbf{u}\cdot(\nabla_{1}A_{0}
\times\mathbf{q}),\nonumber\\
\Updownarrow\nonumber\\
\mathbf{u}\cdot\mathbf{b}=i\nabla_{1}A_{0}\cdot\mathbf{t},\label{Eq8.64}
\end{gather}
which is consistent with (\ref{Eq8.63}). Thus this solvability condition is
automatically satisfied. Taking the longitudinal part of (\ref{Eq8.58}) we get
\begin{gather}
-i\omega n^{2}(\omega)\mathbf{u}\cdot\mathbf{a}\mathbf{=}c^{2}k\mathbf{u\cdot
(\nabla}_{1}A_{0}\times\mathbf{t}),\nonumber\\
\Updownarrow\nonumber\\
\mathbf{u}\cdot\mathbf{a}=i\frac{\omega}{k}\nabla_{1}A_{0}\cdot\mathbf{q},\label{Eq8.65}
\end{gather}
which is consistent with (\ref{Eq8.62}). Thus this solvability condition is
also automatically satisfied. 

The transversal part of (\ref{Eq8.57}) and
(\ref{Eq8.58}) are
\begin{align}
-i\omega b_{\bot}+ik\mathbf{u}\times a_{\bot}  & =-\{k\partial_{t_{1}}
A_{0}+\omega\nabla_{1}A_{0}\cdot\mathbf{u\}t},\label{Eq8.66}\\
-i\omega n^{2}(\omega)a_{\bot}-ic^{2}k\mathbf{u}\times b_{\bot}  &
=-\{\omega(n^{2}(\omega)+\omega\widehat{\chi}^{\prime}(\omega))\partial
_{t_{1}}A_{0}+c^{2}k\nabla_{1}A_{0}\cdot\mathbf{u}\}\mathbf{q},\nonumber
\end{align}
and this linear system is singular; the determinant is zero because of the
dispersion relation (\ref{Eq8.53}). It can therefore only be solved if the
right-hand side satisfy a certain solvability condition. The most effective
way to find this condition is to use the \textit{Fredholm Alternative}. It say
that a linear system
\[
A\mathbf{x}=\mathbf{c},
\]
has a solution if and only if
\[
\mathbf{f}\cdot\mathbf{c}=0,
\]
for all vectors $\mathbf{f}$, such that
\[
A^{\dag}\mathbf{f}=0,
\]
where $A^{\dag}$ is the adjoint of $A$.

The matrix for the system (\ref{Eq8.66}) is
\[
M=\left(
\begin{tabular}
[c]{ll}
$ik\mathbf{u}\times$ & $-i\omega$\\
$-i\omega n^{2}$ & $-ic^{2}k\mathbf{u}\times$
\end{tabular}
\right).
\]
The adjoint of this matrix is clearly

\begin{equation}
M^{\dag}=\left(
\begin{tabular}
[c]{ll}
$-ik\mathbf{u}\times$ & $-i\omega n^{2}$\\
$-i\omega$ & $ic^{2}k\mathbf{u}\times$
\end{tabular}
\right), \label{Eq8.67}
\end{equation}
and the null space of the adjoint is thus determined by
\begin{align}
-ik\mathbf{u}\times\mathbf{\alpha}-i\omega n^{2}\mathbf{\beta}  &
=0,\nonumber\\
-i\omega\mathbf{\alpha}+ic^{2}k\mathbf{u}\times\mathbf{\beta}  & =0.\label{Eq8.68}
\end{align}
A convenient basis for the null space is
\begin{equation}
\left\{  \left(
\begin{tabular}
[c]{l}
$-c^{2}k\mathbf{q}$\\
$\omega\mathbf{t}$
\end{tabular}
\right)  ,\left(
\begin{tabular}
[c]{l}
$c^{2}k\mathbf{t}$\\
$\omega\mathbf{q}$
\end{tabular}
\right)  \right\} \label{Eq8.69}
\end{equation}
The first basis vector gives a trivial solvability condition, whereas the
second one gives a nontrivial condition, which is
\begin{gather}
c^{2}k\{k\partial_{t_{1}}A_{0}+\omega\nabla_{1}A_{0}\cdot\mathbf{u\}+\omega
}\{\omega(n^{2}(\omega)+\omega\widehat{\chi}^{\prime}(\omega))\partial_{t_{1}
}A_{0}+c^{2}k\nabla_{1}A_{0}\cdot\mathbf{u}\}=0,\nonumber\\
\Updownarrow\nonumber\\
\omega^{2}(2n^{2}+\omega\widehat{\chi}^{\prime}(\omega))\partial_{t_{1}}
A_{0}+2c^{2}k\omega\mathbf{u}\cdot\nabla_{1}A_{0}=0.\label{Eq8.70}
\end{gather}
Observe that from the dispersion relation (\ref{Eq8.53}) we have
\begin{gather}
\omega^{2}n^{2}(\omega)=\omega^{2}(1+\widehat{\chi}(\omega))=c^{2}
k^{2},\nonumber\\
\Downarrow\nonumber\\
2\omega\omega^{\prime}n^{2}+\omega^{2}\widehat{\chi}^{\prime}(\omega
)\omega^{\prime}=2c^{2}k,\nonumber\\
\Downarrow\nonumber\\
\omega(2n^{2}+\omega\widehat{\chi}^{\prime}(\omega))\omega^{\prime}
=2c^{2}k.\label{Eq8.71}
\end{gather}
Using (\ref{Eq8.71}) in (\ref{Eq8.70}) the solvability condition can be
compactly written as
\begin{equation}
\partial_{t_{1}}A_{0}+\mathbf{v}_{g}\cdot\nabla_{1}A_{0}=0,\label{Eq8.72}
\end{equation}
where $\mathbf{v}_{g}$ is the \textit{group velocity}
\begin{equation}
\mathbf{v}_{g}=\frac{d\omega}{dk}\mathbf{u}.\label{Eq8.73}
\end{equation}
The system (\ref{Eq8.66}) is singular but consistent. We can therefore
disregard the second equation, and look for a special solution of the form
\begin{align}
a_{\bot}  & =a\mathbf{q},\nonumber\\
b_{\bot}  & =0.\label{Eq8.74}
\end{align}
Inserting (\ref{Eq8.74}) into the first equation in (\ref{Eq8.66}) we easily get
\begin{equation}
a_{\bot}=i\left\{  \partial_{t_{1}}A_{0}+\frac{\omega}{k}\mathbf{u}\cdot
\nabla_{1}A_{0}\right\}  \mathbf{q}.\label{Eq8.75}
\end{equation}
From (\ref{Eq8.62}),(\ref{Eq8.63}),(\ref{Eq8.74}) and (\ref{Eq8.75}), we get
the following bounded special solution to the order $\varepsilon$ equations
\begin{align}
\mathbf{e}_{1}  & =\{i(\partial_{t_{1}}A_{0}+\frac{\omega}{k}\mathbf{u}
\cdot\nabla_{1}A_{0})\mathbf{q}+i(\frac{\omega}{k}\mathbf{q}\cdot\nabla
_{1}A_{0})\mathbf{u\}e}^{i\theta_{0}}+(\ast),\nonumber\\
\mathbf{b}_{1}  & =\{i(\mathbf{t}\cdot\nabla_{1}A_{0})\mathbf{u\}e}
^{i\theta_{0}}+(\ast).\label{Eq8.76}
\end{align}
The order $\varepsilon^{2}$ equations are
\begin{gather}
\partial_{t_{0}}\mathbf{b}_{2}\mathbf{+\nabla}_{0}\mathbf{\times e}
_{2}=-\{\partial_{t_{1}}\mathbf{b}_{1}+\nabla_{1}\times\mathbf{e}_{1}
+\partial_{t_{2}}\mathbf{b}_{0}+\nabla_{2}\times\mathbf{e}_{0}\},\nonumber
\\
\nonumber\\
\partial_{t_{0}}\mathbf{e}_{2}\mathbf{-c}^{2}\mathbf{\nabla}_{0}\mathbf{\times
b}_{2}\mathbf{+\partial}_{t_{0}}\widehat{\chi}(i\partial_{t_{0}}
)\mathbf{e}_{2}=-\{\partial_{t_{1}}\mathbf{e}_{1}-c^{2}\nabla_{1}
\times\mathbf{b}_{1}+\partial_{t_{2}}\mathbf{e}_{0}\nonumber\\
-c_{2}\nabla_{2}\times\mathbf{b}_{0}+\partial_{t_{1}}\widehat{\chi}
(i\partial_{t_{0}})\mathbf{e}_{1}+i\partial_{t_{0}}\widehat{\chi}^{\prime
}(i\partial_{t_{0}})\partial_{t_{1}}\mathbf{e}_{1}\nonumber\\
\partial_{t_{2}}\widehat{\chi}(i\partial_{t_{0}})\mathbf{e}_{0}+i\partial
_{t_{1}}\widehat{\chi}^{\prime}(i\partial_{t_{0}})\partial_{t_{1}}
\mathbf{e}_{0}+i\partial_{t_{0}}\widehat{\chi}^{\prime}(i\partial_{t_{0}
})\partial_{t_{2}}\mathbf{e}_{0}\nonumber\\
-\frac{1}{2}\partial_{t_{0}}\widehat{\chi}^{\prime\prime}(i\partial_{t_{0}
})\partial_{t_{1}t_{1}}\mathbf{e}_{0}+\eta\partial_{t_{0}}\mathbf{e}_{0}
^{2}\mathbf{e}_{0}\},\nonumber\\
\nonumber\\
\nabla_{0}\cdot\mathbf{b}_{2}=-\{\nabla_{1}\cdot\mathbf{b}_{1}+\nabla_{2}
\cdot\mathbf{b}_{0}\},\nonumber\\
\nonumber\\
\nabla_{0}\cdot\mathbf{e}_{2}+\widehat{\chi}(i\partial_{t_{0}})\nabla_{0}
\cdot\mathbf{e}_{2}=-\{\nabla_{1}\cdot\mathbf{e}_{1}+\nabla_{2}\cdot
\mathbf{e}_{0}+\widehat{\chi}(i\partial_{t_{0}})\nabla_{1}\cdot\mathbf{e}
_{1}\nonumber\\
+i\widehat{\chi}^{\prime}(i\partial_{t_{0}})\partial_{t_{1}}\nabla_{0}
\cdot\mathbf{e}_{1}+\widehat{\chi}(i\partial_{t_{0}})\nabla_{2}\cdot
\mathbf{e}_{0}+i\widehat{\chi}^{\prime}(i\partial_{t_{0}})\partial_{t_{2}
}\nabla_{0}\cdot\mathbf{e}_{0}\nonumber\\
+i\widehat{\chi}^{\prime}(i\partial_{t_{0}})\partial_{t_{1}}\nabla_{1}
\cdot\mathbf{e}_{0}-\frac{1}{2}\widehat{\chi}^{\prime\prime}(i\partial_{t_{0}
})\partial_{t_{1}t_{1}}\nabla_{0}\cdot\mathbf{e}_{0}+\eta\nabla_{0}
\cdot(\mathbf{e}_{0}^{2}\mathbf{e}_{0})\}.\label{Eq8.77}
\end{gather}

We now insert (\ref{Eq8.51}) and (\ref{Eq8.76}) into (\ref{Eq8.77}). This
gives us
\begin{gather}
\partial_{t_{0}}\mathbf{b}_{2}\mathbf{+\nabla}_{0}\mathbf{\times e}
_{2}=-\{i(\partial_{t_{1}}\nabla_{1}A_{0}\cdot\mathbf{t)u}+i\mathbf{\nabla
}_{1}\partial_{t_{1}}A_{0}\times\mathbf{q},\nonumber\\
+i\frac{\omega}{k}\nabla_{1}(\nabla_{1}A_{0}\cdot\mathbf{u)\times q}
+i\frac{\omega}{k}\nabla_{1}(\nabla_{1}A_{0}\cdot\mathbf{q)\times
u}+k\mathbf{\partial}_{t_{2}}A_{0}\mathbf{t}\nonumber\\
\mathbf{+\omega\nabla}_{2}A_{0}\times\mathbf{q\}e}^{i\theta_{0}}
+(\ast),\nonumber\\
\nonumber\\
\partial_{t_{0}}\mathbf{e}_{2}\mathbf{-}c^{2}\mathbf{\nabla}_{0}\mathbf{\times
b}_{2}\mathbf{+\partial}_{t_{0}}\widehat{\chi}(i\partial_{t_{0}}
)\mathbf{e}_{2}=-\{iF(\omega)\partial_{t_{1}t_{1}}A_{0}\mathbf{q}\nonumber\\
+i\mathbf{G(\omega)(\partial}_{t_{1}}\nabla_{1}A_{0}\cdot\mathbf{u)q}
+iG(\omega)(\partial_{t_{1}}\nabla_{1}A_{0}\cdot\mathbf{q)u}-ic^{2}\nabla
_{1}(\nabla_{1}A_{0}\cdot\mathbf{t})\times\mathbf{u}\nonumber\\
-c^{2}k\nabla_{2}A_{0}\times\mathbf{t}+H\mathbf{(\omega)\partial}_{t_{2}}
A_{0}\mathbf{q}-3i\eta\omega^{4}|A_{0}|^{2}A_{0}\}e^{i\theta_{0}}
+(\ast)\nonumber,\\
\nonumber\\
\nabla_{0}\cdot\mathbf{b}_{2}=-\{i\nabla_{1}(\nabla_{1}A_{0}\cdot
\mathbf{t})\cdot\mathbf{u}+k\mathbf{\nabla}_{2}A_{0}\cdot\mathbf{t\}e}
^{i\theta_{0}}+(\ast),\nonumber\\
\nonumber\\
\nabla_{0}\cdot\mathbf{e}_{2}+\widehat{\chi}(i\partial_{t_{0}})\nabla_{0}
\cdot\mathbf{e}_{2}=-\{in^{2}\nabla_{1}\partial_{t_{1}}A_{0}\cdot
\mathbf{q}+in^{2}\frac{\omega}{k}\nabla_{1}(\nabla_{1}A_{0}\cdot
\mathbf{u})\cdot\mathbf{q}\nonumber\\
+in^{2}\frac{\omega}{k}\nabla_{1}(\nabla_{1}\cdot\mathbf{q})\cdot
\mathbf{u}+\omega n^{2}\nabla_{2}A_{0}\cdot\mathbf{q}\}e^{i\theta_{0}}
+(\ast),\label{Eq8.78}
\end{gather}
where we have defined
\begin{align}
F(\omega)  & =n^{2}+2\omega\widehat{\chi}^{\prime}(\omega)+\frac{1}{2}
\omega^{2}\widehat{\chi}^{\prime\prime}(\omega),\nonumber\\
G(\omega)  & =\frac{\omega}{k}(n^{2}+\omega\widehat{\chi}^{\prime}
(\omega)),\nonumber\\
H(\omega)  & =\omega(n^{2}+\omega\widehat{\chi}^{\prime}(\omega)).\label{Eq8.79}
\end{align}
Like for the order $\varepsilon$ equations, we will look for bounded solutions
of the form
\begin{align}
\mathbf{e}_{2}  & =\mathbf{a}e^{i\theta_{0}}+(\ast),\nonumber\\
\mathbf{b}_{2}  & =\mathbf{b}e^{i\theta_{0}}+(\ast).\label{Eq8.80}
\end{align}
Inserting (\ref{Eq8.80}) into (\ref{Eq8.78}) we get the following linear
system of equations for the constant vectors $\mathbf{a}$ and $\mathbf{b}$
\begin{gather}
-i\omega\mathbf{b}+ik\mathbf{u}\times\mathbf{a}=-\{i\mathbf{(}\partial_{t_{1}
}\nabla_{1}A_{0}\cdot\mathbf{t})\mathbf{u}+i\mathbf{\nabla}_{1}\partial
_{t_{1}}A_{0}\times\mathbf{q}\nonumber\\
+i\frac{\omega}{k}\nabla_{1}(\nabla_{1}A_{0}\cdot\mathbf{u)\times q}
+i\frac{\omega}{k}\nabla_{1}(\nabla_{1}A_{0}\cdot\mathbf{q)\times
u}+k\mathbf{\partial}_{t_{2}}A_{0}\mathbf{t}\nonumber\\
\mathbf{+\omega\nabla}_{2}A_{0}\times\mathbf{q\}},\label{Eq8.81}\\
\nonumber\\
-i\omega n^{2}(\omega)\mathbf{a}-ic^{2}k\mathbf{u}\times\mathbf{b}
=-\{iF(\omega)\partial_{t_{1}t_{1}}A_{0}\mathbf{q}\nonumber\\
+iG\mathbf{(\omega)(\partial}_{t_{1}}\nabla_{1}A_{0}\cdot\mathbf{u)q}
+iG\mathbf{(\omega)(\partial}_{t_{1}}\nabla_{1}A_{0}\cdot\mathbf{q)u}
\nonumber\\
-ic^{2}\nabla_{1}(\nabla_{1}A_{0}\cdot\mathbf{t})\times\mathbf{u}-c^{2}
k\nabla_{2}A_{0}\times\mathbf{t}\nonumber\\
+H\mathbf{(\omega)\partial}_{t_{2}}A_{0}\mathbf{q}-3i\eta\omega^{4}|A_{0}
|^{2}A_{0}\},\label{Eq8.82}\\
\nonumber\\
ik\mathbf{u}\cdot\mathbf{b}=-\{i\nabla_{1}(\nabla_{1}A_{0}\cdot\mathbf{t}
)\cdot\mathbf{u}+k\mathbf{\nabla}_{2}A_{0}\cdot\mathbf{t\}},\label{Eq8.83}\\
\nonumber\\
ikn^{2}\mathbf{u}\cdot\mathbf{a}=-\{in^{2}\nabla_{1}\partial_{t_{1}}A_{0}
\cdot\mathbf{q}+in^{2}\frac{\omega}{k}\nabla_{1}(\nabla_{1}A_{0}
\cdot\mathbf{u})\cdot\mathbf{q}\nonumber\\
+in^{2}\frac{\omega}{k}\nabla_{1}(\nabla_{1}\cdot\mathbf{q})\cdot
\mathbf{u}+\omega n^{2}\nabla_{2}A_{0}\cdot\mathbf{q}\}.\label{Eq8.84}
\end{gather}
We introduce the  longitudinal and transversal vector components for $\mathbf{a}$
and $\mathbf{b}$ like before, and find from (\ref{Eq8.83}) and (\ref{Eq8.84})
that%
\begin{align}
a_{\Vert}  & =(-\frac{1}{k}\nabla_{1}\partial_{t_{1}}A_{0}\cdot\mathbf{q}
-\frac{\omega}{k^{2}}\nabla_{1}(\nabla_{1}A_{0}\cdot\mathbf{u})\cdot
\mathbf{q}\nonumber\\
& -\frac{\omega}{k^{2}}\nabla_{1}(\nabla_{1}A_{0}\cdot\mathbf{q}
)\cdot\mathbf{u}+i\frac{\omega}{k}\nabla_{2}A_{0}\cdot\mathbf{q})\mathbf{u},\label{Eq8.85}\\
b_{\Vert}  & =(i\nabla_{2}A_{0}\cdot\mathbf{t}-\frac{1}{k}\nabla_{1}
(\nabla_{1}\cdot\mathbf{t})\cdot\mathbf{u})\mathbf{u}.\label{Eq8.86}
\end{align}
The longitudinal part of (\ref{Eq8.81}) is
\begin{equation}
\mathbf{u}\cdot\mathbf{b}=\frac{1}{\omega}\{\partial_{t_{1}}\nabla_{1}
A_{0}\cdot\mathbf{t-\nabla}_{1}\partial_{t_{1}}A_{0}\cdot\mathbf{t}
-\frac{\omega}{k}\nabla_{1}(\nabla_{1}A_{0}\cdot\mathbf{u})\cdot
\mathbf{t}+i\omega\mathbf{\nabla}_{2}A_{0}\cdot\mathbf{t}\},\label{Eq8.87}
\end{equation}
and in order for (\ref{Eq8.87}) to be consistent with (\ref{Eq8.86}), we find that
the following solvability condition must hold
\begin{equation}
\partial_{t_{1}}\nabla_{1}A_{0}\cdot\mathbf{t}=\nabla_{1}\partial_{t_{1}}
A_{0}\cdot\mathbf{t}.\label{Eq8.88}
\end{equation}
The longitudinal part of (\ref{Eq8.82}) is
\begin{equation}
\mathbf{u}\cdot\mathbf{a=}\frac{1}{\omega n^{2}}\{G(\omega)\partial_{t_{1}
}\nabla_{1}A_{0}\cdot\mathbf{q}+ic^{2}k\nabla_{2}A_{0}\cdot\mathbf{q}
\},\label{Eq8.89}
\end{equation}
and in order for (\ref{Eq8.89}) to be consistent with (\ref{Eq8.85}) we find,
after a little algebra, that the solvability condition
\begin{gather}
\frac{\omega}{k}n^{2}(\omega)\nabla_{1}\partial_{t_{1}}A_{0}\cdot
\mathbf{q}+G(\omega)\mathbf{\partial}_{t_{1}}\nabla_{1}A_{0}\cdot
\mathbf{q=}\nonumber\\
-c^{2}\nabla_{1}(\nabla_{1}A_{0}\cdot\mathbf{q)\cdot u}-c^{2}\nabla_{1}
(\nabla_{1}A_{0}\cdot\mathbf{u})\cdot\mathbf{q},\label{Eq8.90}
\end{gather}
must hold. 

The transverse parts of (\ref{Eq8.81}) and (\ref{Eq8.82}) are
\begin{gather}
-i\omega b_{\bot}+ik\mathbf{u}\times a_{\bot}=-\{i\nabla_{1}\partial_{t_{1}
}A_{0}\cdot\mathbf{u}+i\frac{\omega}{k}\nabla_{1}(\nabla_{1}A_{0}
\cdot\mathbf{u})\cdot\mathbf{u}\nonumber\\
-i\frac{\omega}{k}\nabla_{1}(\nabla_{1}A_{0}\cdot\mathbf{q})\cdot
\mathbf{q}+k\partial_{t_{2}}A_{0}+\omega\nabla_{2}A_{0}\cdot\mathbf{u\}t}
-\{i\frac{\omega}{k}\nabla_{1}(\nabla_{1}A_{0}\cdot\mathbf{q)\cdot
t\}q},\nonumber\\
\nonumber\\
-i\omega n^{2}a_{\bot}-ic^{2}k\mathbf{u}\times b_{\bot}=-\{iF(\omega
)\partial_{t_{1}t_{1}}A_{0}+iG(\omega)\partial_{t_{1}}\nabla_{1}A_{0}
\cdot\mathbf{u}\nonumber\\
-ic^{2}\nabla_{1}(\nabla_{1}A_{0}\cdot\mathbf{t})\cdot\mathbf{t}+c^{2}
k\nabla_{2}A_{0}\cdot\mathbf{u}+H(\omega)\partial_{t_{2}}A_{0}-3\eta
i\omega^{4}|A_{0}|^{2}A_{0}\}\mathbf{q}\nonumber\\
-\{ic^{2}\nabla_{1}(\nabla_{1}A_{0}\cdot\mathbf{t})\cdot\mathbf{q}
\}\mathbf{t}.\label{Eq8.91}
\end{gather}

\noindent  The matrix for this linear system is the same as for the order $\varepsilon$
case, (\ref{Eq8.66}), so that the two solvability conditions are determined,
through the Fredholm Alternative, by the vectors (\ref{Eq8.69}). The
solvability condition corresponding to the first of the vectors in
(\ref{Eq8.69}) is
\begin{gather}
(-c^{2}k)(-i\frac{\omega}{k}\nabla_{1}(\nabla_{1}A_{0}\cdot\mathbf{q)\cdot
t)}+\omega\mathbf{(-}ic^{2}\nabla_{1}(\nabla_{1}A_{0}\cdot\mathbf{t}
)\cdot\mathbf{q})=0,\nonumber\\
\Updownarrow\nonumber\\
\nabla_{1}(\nabla_{1}A_{0}\cdot\mathbf{q})\cdot\mathbf{t}=\nabla_{1}
(\nabla_{1}\cdot\mathbf{t})\cdot\mathbf{q},\label{Eq8.92}
\end{gather}
and the solvability condition corresponding to the second vector in
(\ref{Eq8.69}) is
\begin{gather}
c^{2}k(-\{i\nabla_{1}\partial_{t_{1}}A_{0}\cdot\mathbf{u}+i\frac{\omega}
{k}\nabla_{1}(\nabla_{1}A_{0}\cdot\mathbf{u})\cdot\mathbf{u}\nonumber\\
-i\frac{\omega}{k}\nabla_{1}(\nabla_{1}A_{0}\cdot\mathbf{q})\cdot
\mathbf{q}+k\partial_{t_{2}}A_{0}+\omega\nabla_{2}A_{0}\cdot\mathbf{u\})}
+\omega(-\{iF(\omega)\partial_{t_{1}t_{1}}A_{0}\nonumber\\
+iG(\omega)\partial_{t_{1}}\nabla_{1}A_{0}\cdot\mathbf{u}-ic^{2}\nabla
_{1}(\nabla_{1}A_{0}\cdot\mathbf{t})\cdot\mathbf{t}+c^{2}k\nabla_{2}A_{0}
\cdot\mathbf{u}\nonumber\\
+H(\omega)\partial_{t_{2}}A_{0}-3\eta i\omega^{4}|A_{0}|^{2}A_{0}
\}\mathbf{q})=0,\nonumber\\
\Updownarrow\nonumber\\
\partial_{t_{2}}A_{0}+\mathbf{v}_{g}\cdot\nabla_{2}A_{0}+i\delta_{1}\nabla
_{1}\partial_{t_{1}}A_{0}\cdot\mathbf{u}+i\delta_{2}\partial_{t_{1}}\nabla
_{1}A_{0}\cdot\mathbf{u}\nonumber\\
-i\beta(\nabla_{1}(\nabla_{1}A_{0}\cdot\mathbf{q)\cdot q}+\nabla_{1}
(\nabla_{1}A_{0}\cdot\mathbf{t})\cdot\mathbf{t-\nabla}_{1}(\nabla_{1}
A_{0}\cdot\mathbf{u})\cdot\mathbf{u})\nonumber\\
+i\alpha\partial_{t_{1}t_{1}}A_{0}-i\gamma|A_{0}|^{2}A_{0}=0,\label{Eq8.93}
\end{gather}
where we have defined
\begin{align}
\alpha & =\frac{\omega^{\prime}F(\omega)}{2c^{2}k},\nonumber\\
\beta & =\frac{\omega^{\prime}}{2k},\nonumber\\
\gamma & =\frac{3\eta\omega^{\prime}\omega^{4}}{2c^{2}k},\nonumber\\
\delta_{1}  & =\frac{\omega^{\prime}}{2\omega},\nonumber\\
\delta_{2}  & =\frac{\omega^{\prime}G(\omega)}{2c^{2}k}.\nonumber
\end{align}
We have now found all solvability conditions. These are (\ref{Eq8.88}
),(\ref{Eq8.90}),(\ref{Eq8.92}) and (\ref{Eq8.93}). 

We now, as usual, define an amplitude
$A(\mathbf{x},t)$ by
\[
A(\mathbf{x},t)=A_{0}(\mathbf{x}_{1},t_{1},...)|\mathbf{x}_{j}=\varepsilon
^{j}\mathbf{x},t_{j}=\varepsilon^{j}t,
\]
and derive the amplitude equations from the solvability conditions in the
usual way. This gives us the following system
\begin{gather}
\partial_{t}\nabla A\cdot\mathbf{t=\nabla}\partial_{t}A\cdot\mathbf{t},
\label{Eq8.95}\\
\nonumber\\
\frac{\omega}{k}n^{2}(\omega)\nabla\partial_{t}A\cdot\mathbf{q}
+G\mathbf{(\omega)\partial}_{t}\nabla A\cdot\mathbf{q=}\nonumber\\
-c^{2}\nabla(\nabla A\cdot\mathbf{q)\cdot u}-c^{2}\nabla(\nabla A\cdot
\mathbf{u})\cdot\mathbf{q},\label{Eq8.96}\\
\nonumber\\
\nabla(\nabla A_{0}\cdot\mathbf{q})\cdot\mathbf{t}=\nabla(\nabla
A\cdot\mathbf{t})\cdot\mathbf{q},\label{Eq8.97}\\
\nonumber\\
\partial_{t}A+\mathbf{v}_{g}\cdot\nabla A+i\delta_{1}\nabla\partial_{t}
A\cdot\mathbf{u}+i\delta_{2}\partial_{t}\nabla A\cdot\mathbf{u}\nonumber
\\
-i\beta(\nabla(\nabla A\cdot\mathbf{q)\cdot q}+\nabla(\nabla A\cdot
\mathbf{t})\cdot\mathbf{t-\nabla}(\nabla A\cdot\mathbf{u})\cdot\mathbf{u}
)\nonumber\\
+i\alpha\partial_{tt}A-i\gamma|A|^{2}A=0,\label{Eq8.98}
\end{gather}
where we as usual have set the formal perturbation parameter equal to $1$.
Equations (\ref{Eq8.95}) and (\ref{Eq8.97}) are automatically satisfied since
$A(\mathbf{x},t)$ is a smooth function of space and time. We know that only
amplitudes such that
\begin{equation}
\partial_{t}A\sim-\mathbf{v}_{g}\cdot\nabla A=\omega^{\prime}\nabla
A\cdot\mathbf{u},\label{Eq8.99}
\end{equation}
can be allowed as solutions. This is assumed by the multiple scale method. If
we insert (\ref{Eq8.99}) into (\ref{Eq8.96}), assume smoothness and use the
dispersion relation, we find that (\ref{Eq8.96}) is automatically satisfied.
The only remaining equation is then (\ref{Eq8.98}) and if we insert the
approximation (\ref{Eq8.99}) for the derivatives with respect to time in the
second and third term of (\ref{Eq8.98}) we get, using the dispersion relation,
that (\ref{Eq8.98}) simplify into
\begin{equation}
\partial_{t}A+\mathbf{v}_{g}\cdot\nabla A-i\beta\nabla^{2}A+i\alpha
\partial_{tt}A-i\gamma|A|^{2}A=0,\label{Eq8.100}
\end{equation}
where we have also used the fact that
\[
\mathbf{qq}+\mathbf{tt}+\mathbf{uu}=I.
\]
The amplitude $A$ determines the electric and magnetic fields through the
identities
\begin{align}
\mathbf{E}(\mathbf{x},t)  & \approx\{(\omega A+i(\frac{\omega}{k}
-\omega^{\prime})\mathbf{u}\cdot\nabla A\mathbf{)q}\nonumber\\
& +i(\frac{\omega}{k}\mathbf{q}\cdot\nabla A)\mathbf{u}\}e^{i(\mathbf{k}
\cdot\mathbf{x}-\omega t)}+(\ast)\nonumber\\
\mathbf{B}(x,t)  & \approx\{kA\mathbf{t}+i(\mathbf{t}\cdot\nabla
A)\mathbf{u}\}e^{i(\mathbf{k}\cdot\mathbf{x}-\omega t)}+(\ast).\label{Eq8.101}
\end{align}
The equations (\ref{Eq8.100}) and (\ref{Eq8.101}) are the key elements in a
fast numerical scheme for linearly polarized wave packet solutions to
Maxwell's equations. Wave packets of circular polarization or arbitrary
polarization can be treated in an entirely similar manner, as can sums of
different polarized wave packets and materials with nontrivial magnetic response.

The reader will of course recognize the amplitude equation (\ref{Eq8.100}) as the 3D
nonlinear Schrodinger equation including group velocity dispersion. As we have
seen before, an equation like this can be solved as an ordinary initial value
problem if we first use (\ref{Eq8.99}) to make the term containing a second
derivative\ with respect to time into one containing only a first derivative with respect to time.

The derivation of the nonlinear Scrodinger equation for linearly polarized
wave packets I have given in this section is certainly not the fastest and
simplest way this can be done. The main aim in this section was to illustrate
how to apply the multiple scale method for vector PDEs in general, not to do
it in the most effective way possible for the particular case of linearly
polarized electromagnetic wave packets in non-magnetic materials. If the material
has a significant magnetic response, a derivation along the lines given is necessary.

All the essential elements we need in order to apply the method of multiple
scales to problems in optics and laser physics, and other areas of science
too, are at this point known. There are no new tricks to learn. Using the
approach described in these lecture notes, amplitude equations can be derived
for most situations of interest. Applying the method is mechanical, but for
realistic systems there can easily be a large amount of algebra involved. This
is unavoidable; solving nonlinear partial differential equations, even
approximately, is hard.

In these lecture notes we have focused on applications of the multiple scale
method for time-propagation problems. The method was originally developed for
these kind of problems and the mechanics of the method is most transparent for
such problems. However the method is by no means limited to time propagation problems.

Many pulse propagation schemes are most naturally formulated as a boundary
value problem where the propagation variable is a space variable. A very
general scheme of this type is the well known UPPE\cite{UPPE} propagation
scheme. More details on how the multiple scale method is applied for these
kind of schemes can be found in \cite{per1} and \cite{newell}.
\pagebreak

\bigskip

\bibliographystyle{abbrv}
\bibliography{PerturbationMethodLectureNotes}

\end{document}